\documentclass[EJP]{ejpecp} 

\usepackage{stackrel}
\usepackage{mathtools}
\usepackage[scr=boondox,  
            cal=cm]   
           {mathalpha}
\usepackage[comma, sort&compress]{natbib}
\usepackage{multirow}
\usepackage{accents}
\usepackage{amssymb}
\usepackage{xcolor}
\usepackage{paralist}

\newtheorem{property}[theorem]{Property}
\newtheorem{configuration}[theorem]{Configuration}

\numberwithin{equation}{section}
\numberwithin{theorem}{section}

\numberwithin{figure}{section}

\newcommand{\Tsn}{T_{SN}}

\newcommand{\barPhi}{\bar{\Phi}}

\newcommand{\frakB}{\mathfrak{B}}
\newcommand{\frakz}{\mathfrak{z}}

\newcommand{\cR}{\mathcal{R}}

\newcommand{\cX}{\mathcal{X}}

\newcommand{\cM}{\mathcal{M}}

\newcommand{\cE}{\mathcal{E}}
\newcommand{\barDelta}{\bar{\Delta}}
\newcommand{\scE}{\mathscr{E}} 

\newcommand{\cN}{\mathcal{N}}
\newcommand{\cZ}{\mathcal{Z}}

\newcommand{\bX}{{\bf X}}

\newcommand{\boldi}{{\bf i}}
\newcommand{\boldj}{{\bf j}}
\newcommand{\boldk}{{\bf k}}

\newcommand{\bR}{\mathbb{R}}

\newcommand{\barD}{{\bar{D}}}

\newcommand{\barW}{{\bar{W}}}
\newcommand{\barxi}{{\bar{\xi}}}
\newcommand{\bareta}{{\bar{\eta}}}
\newcommand{\bark}{{\bar{k}}}

\newcommand{\Inm}{{I_{n,m}}}

\newcommand{\bE}{\mathbb{E}}

\newcommand{\configref}[1]{Configuration~\ref{config:#1}}
\newcommand{\configsref}[1]{Configurations~\ref{config:#1}}
\newcommand{\configssref}[1]{\ref{config:#1}}

\newcommand{\secref}[1]{Section~\ref{sec:#1}}
\newcommand{\secsref}[1]{Sections~\ref{sec:#1}}
\newcommand{\secssref}[1]{\ref{sec:#1}}

\newcommand{\appref}[1]{Appendix~\ref{app:#1}}
\newcommand{\appsref}[1]{Appendices~\ref{app:#1}}
\newcommand{\appssref}[1]{\ref{app:#1}}

\newcommand{\lemref}[1]{Lemma~\ref{lem:#1}}
\newcommand{\lemsref}[1]{Lemmas~\ref{lem:#1}}
\newcommand{\lemssref}[1]{\ref{lem:#1}}

\newcommand{\thmref}[1]{Theorem~\ref{thm:#1}}

\newcommand{\thmsref}[1]{Theorems~\ref{thm:#1}}
\newcommand{\thmssref}[1]{\ref{thm:#1}}

\newcommand{\propertyref}[1]{Property~\ref{property:#1}}

\SHORTTITLE{B-E theorems for incomplete U-statistics with Bernoulli sampling}

\TITLE{Berry-Esseen theorems for the asymptotic normality of incomplete U-statistics with Bernoulli sampling}

\AUTHORS{%
  Dennis~Leung\footnote{University of Melbourne, Australia.
    \EMAIL{dennis.leung@unimelb.edu.au}}\orcid{0000-0002-5924-8161}
    }

\KEYWORDS{U-statistics;   Berry-Esseen theorems; exponential tail bounds; self-normalized limit theory; Stein's method; variable censoring} 

\AMSSUBJ{62E17} 

\SUBMITTED{July 29, 2024} 
\ACCEPTED{Jan 7, 2026} 

\VOLUME{0}
\YEAR{2025}
\PAPERNUM{0}
\DOI{10.1214/YY-TN}

\ABSTRACT{There has been a   resurgence of interest   in  incomplete U-statistics that only sum over a subset of  kernel evaluations, due to their computational efficiency and asymptotic normality which can be leveraged to quantify the uncertainty of ensemble predictions in machine learning. In this paper, we study the weak  convergences to normality of one such construction, the incomplete U-statistic  with Bernoulli sampling, under three different regimes on the relative sizes of the raw sample and the computational budget. Under minimalistic   moment assumptions, we establish accompanying Berry-Esseen  bounds with  the natural rates that characterize the accuracy of these normal approximations. The key ingredients in our proofs  include a variable censoring technique and a methodology for establishing Berry-Esseen bounds for the so-called Studentized nonlinear statistics recently formalized in the Stein's method literature, as well as an exponential lower tail bound for non-negative kernel U-statistics. }

\begin{document}



\section{Introduction}

Suppose $X_1, \dots, X_n \in M$ are  independent and identically distributed random variables 
  taking values in a given measurable space $(M, \cM)$. \citet{hoeffding1948} first introduced U-statistics, which, for a given degree $m \geq 2$, are statistics of the form
\begin{equation} \label{complete_u_stat}
U_n = \frac{1}{{n \choose m}} \sum_{(i_1, \dots, i_m) \in I_{n, m}} h(X_{i_1}, \dots, X_{i_m}),
\end{equation}
where $\Inm \equiv \{(i_1, \dots, i_m): 1 \leq i_1 < \cdots < i_m \leq n \}$ and $h: M^m \rightarrow \bR$ is a given symmetric function in $m$ arguments, i.e. 
\[
h(x_1, \dots, x_m) = h(x_{\pi_1}, \dots, x_{\pi_m})
\]
 for any permutation $(\pi_1, \dots, \pi_m)$ of the indices $(1, \dots, m)$.  It generalizes the sample mean and encompasses many examples such as the sample variance; see \citet{KB1994Ustat} for a comprehensive treatment of the classical theory on U-statistics. 
  To  simplify notation,  for each $m$-tuple  ${\bf i} =  (i_1, \dots, i_m) \in I_{n,m}$, we will  henceforth adopt  the   shorthands $\bX_{\bf i} \equiv (X_{i_1}, \dots, X_{i_m})$  and $h(\bX_{\bf i}) \equiv h(X_{i_1}, \dots, X_{i_m})$.  
  
     $U_n$ is \emph{complete}, in a sense that it is an average of  the kernel $h$ evaluated at \emph{all possible} ordered subsamples of size $m$ from $\{1, \dots, n\}$. 
\citet{blom1976}  suggested alternative constructions dubbed ``incomplete U-statistics''  which only average over  an appropriate subset of $\Inm$, with the intuition that  the strong dependence among   the summands in \eqref{complete_u_stat} should  allow one to only use a well-chosen subset of them without losing too much statistical efficiency. In this work, we study the asymptotic normality of one such construction, known as the \emph{incomplete U-statistic with Bernoulli sampling}, which has the form
\begin{equation} \label{icu_def}
U_{n, N}' \equiv \frac{1}{\hat{N}}  \sum_{{\bf i}\in I_{n, m}} Z_{{\bf i}}h(\bX_{\bf i}),
\end{equation} 
where  the $Z_{\bf i}$'s are ${n \choose m}$ independent and identically distributed Bernoulli random variables with success probability 
\begin{equation}\label{bernoulli_samp_prob}
p = p_{n, N, m}\equiv N/{n \choose m} \text{ for a positive integer } N < {n \choose m},
\end{equation}
i.e., $P(Z_{\bf i} = 1) = p = N {n \choose m}^{-1}$ and $P(Z_{\bf i} = 0) = 1- p = 1 -  N {n \choose m}^{-1}$, 
and $\hat{N}$ is defined as 
\[
\hat{N} \equiv  \sum_{{\bf i} \in I_{n, m}} Z_{\bf i}.
\]
The Bernoulli samplers $(Z_\boldi)_{\boldi \in \Inm}$ are also  assumed to be independent of the raw data $(X_i)_{i \in [n]}$.  We  call $N$ the \emph{computational budget} because it is the expected number of  kernel evaluations  $\hat{N}$ that  $U_{n, N}'$ winds up summing over.  
In the sequel, we will use $\cX \equiv (X_i)_{i \in [n]}$ and $\cZ \equiv(Z_{\bf i})_{{\bf i} \in \Inm}$ to represent the entirety of  the raw data and Bernoulli samplers, respectively. 

Let $\sigma_h^2  \equiv \text{var}[h(X_1, \dots, X_m)]$. To describe the asymptotic normality of $U_{n, N}'$, throughout this paper, we will  assume
\begin{equation} \label{mean0_assumption}
\bE[h(X_1, \dots, X_m)] = 0 \text{ and }  0 < \sigma_h^2  < \infty,
\end{equation} 
without loss of generality; see  \secref{sp_remark_Nhat} below. For each $r = 1, \dots, m$, we will define the $r$-argument function $h_r : M^r \rightarrow \bR$ associated with the kernel $h$ by
\[
h_r (x_1, \dots, x_r) \equiv \bE[h(x_1, \dots, x_r, X_{r+1}, \dots, X_m )],
\]
and let 
\begin{equation*} \label{kernel_rank}
d \equiv \min \{r: \text{var}[ h_r(X_1, \dots, X_r)] > 0\}
\end{equation*}
be the \emph{rank} of  $h$. Moreover, for simplicity, we rename the single-argument function $h_1$ as 
 \[
g(x) \equiv h_1(x),
\] 
and the kernel $h$ is said to be  non-degenerate when  
\begin{equation} \label{non_degenerate_condition}
\sigma_g^2 \equiv \text{var}[g(X_1)] > 0,
\end{equation}
or equivalently, when $d =1$.
 In the sequel, we will also let
\[
\alpha = \alpha_{n, N} \equiv \frac{n}{N}.
\]
and
\[
\sigma^2 \equiv m^2 \sigma_g^2 + \alpha \sigma_h^2.
\]
The following weak convergences of $U_{n, N}'$ to normality under different regimes on the relative sizes between $n$ and $N$ are implied by \citet[Corollary 1]{janson1984asymptotic}; also see \secref{sp_remark_Nhat} below. We use $\cN(0, 1)$ to denote the standard normal distribution.  

\begin{theorem}[Normal convergence of $U_{n, N}'$ for a given degree $m$] 
\label{thm:weak_convergences}
Fix $m \geq 2$. Suppose both $n$ and $N$ tend to $\infty$ and the assumptions in \eqref{mean0_assumption} holds. 
\begin{enumerate}[(i)]
\item   If $\alpha \longrightarrow 0$, and \eqref{non_degenerate_condition} is true, then $\frac{\sqrt{n} U_{n, N}'}{m\sigma_g}$ converges weakly to $\cN(0, 1)$.

\item    If $n^d/N \longrightarrow \infty$, then $\frac{\sqrt{N}U_{n, N}'}{\sigma_h}$ converges weakly to  $\cN(0, 1)$.

\item   If $\alpha$ tends to a positive constant  and \eqref{non_degenerate_condition} is true, then $\frac{\sqrt{n} U_{n, N}'}{\sigma}$ converges weakly to $\cN(0, 1)$.

\end{enumerate}

\end{theorem}

Note that \thmref{weak_convergences}$(ii)$ does not require the non-degeneracy condition \eqref{non_degenerate_condition}. While  these convergent results  have been known for a long time,  the corresponding  \emph{Berry-Esseen} (B-E) bounds  that characterize their  normal approximation accuracy have been missing from the literature. In recent years, 
these normal convergences  have found a  resurgence of interest among the machine learning community,  since \citet{mentch2016quantifying}  pointed out that, by summing over only a subset of $\Inm$,  incomplete U-statistics provide a computationally friendly framework for quantifying the uncertainty of ensemble predictions.  This has led to \citet{peng2022rates} following up on this line of work by proving B-E bounds for the convergences in \thmref{weak_convergences}$(i)$ and $(iii)$; however, their bounds contain  artefacts that clearly suggest they are not of the optimal rate, as we will explain in \secref{main_results}.


%

 In this work, we precisely offers  B-E bounds with the natural rates for the convergences in \thmref{weak_convergences} to fill this literature gap. Broadly speaking, our refined results are enabled by three technical advances: \begin{inparaenum}[(i)] \item  a framework of using Stein's method to prove B-E bounds   for the so-called ``Studentized nonlinear statistics'' formalized in \citet{leung2024another}, \item the handy ``variable censoring'' technique introduced in the same paper, and 
  \item an  exponential  lower tail bound  for the non-negative U-statistic proven in \citet{leungshao2024nonuniform}. 
  \end{inparaenum}  
    In passing, we mention that the statistical community has also been  developing B-E bounds for \emph{high-dimensional} incomplete U-statistics with  vector-valued kernels \citep{SongChenKato2019, chen2019randomized}, but  our focus here is  to provide bounds of the best rate for the one-dimensional convergent results in \thmref{weak_convergences}.    
    
   Our organization is simple: In \secref{heuristics}, we first review the nature of the normal convergences in \thmref{weak_convergences}. \secref{main_results} will then present our main Berry-Esseen bound results (\thmsref{BE_bdd_Nggn}-\thmssref{main}), the proofs of which are split into the ensuing \secsref{pf_Nggn_case} and \secssref{pf_Nlln_Nasympn}. We finish by discussing further extensions of our work in \secref{discuss}.

   In terms of notation,
 $P(\cdot)$, $\bE[\cdot]$ denote generic probability and expectation  operators. For any $q \geq 1$,  we use $\|Y\|_q = ( \mathbb{E}|Y|^q)^{1/q}$ to denote the $L_q$-norm of  any real-valued random variable $Y$. Given $r \in \{1, \dots, n\}$ and a function $t:M^r \rightarrow \bR$ in $r$  arguments in $M$, we  use $\bE[t]$ and $\|t\|_q$ as shorthands for $\bE[t(X_1, \dots, X_r)]$ and $\|t(X_1, \dots, X_r)\|_q$ respectively; in addition, if $\ell \in \bR$ and $t$ is a non-negative function, then $\bE[t^\ell]$ is the shorthand for $\bE[(t(X_1, \dots, X_r))^\ell]$.  
$\phi(\cdot)$ and $\Phi(\cdot)$ are respectively the standard normal density and distribution functions, with $\barPhi(\cdot) \equiv 1 - \Phi(\cdot)$; 
$\cN(\mu, v^2)$ denotes the univariate normal distribution with mean $\mu$ and variance $v^2$. The indicator function is denoted by $I(\cdot)$. Generally speaking,
 $C$  denotes an unspecified absolute positive constant, where ``absolute'' means it is universal for the statement it appears in and doesn't depend on other quantities; 
it generally varies in values at different occurrences. For two real numbers $a$ and $b$, we may also use the shorthands $a\vee b \equiv \max (a, b)$ and $a\wedge b \equiv \min (a, b)$.


%

We make two extra remarks before proceeding:

\subsection{The random normalizer \texorpdfstring{$\hat{N}$}{} of $U'_{n, N}$ } \label{sec:sp_remark_Nhat}

Suppose we do not assume $\bE[h] = 0$ as in  \eqref{mean0_assumption}. Because $U_{n, N}'$ is normalized by $\hat{N}$, one can write
   \[
   U'_{n, N} - \bE[h]  =  \frac{1}{\hat{N}}  \sum_{{\bf i}\in I_{n, m}} Z_{{\bf i}}(h(\bX_{\bf i}) - \bE[h]);
   \]
   which is the same incomplete U-statistic construction in \eqref{icu_def} based on the centered kernel $h(\cdot) - \bE[h]$. Therefore,  the three convergences in \thmref{weak_convergences}  can be restated simply by replacing all the occurrences of $U'_{n, N}$ by $(U'_{n, N} - \bE[h])$; for instance, \thmref{weak_convergences}$(ii)$ can be stated as the weak convergence of $\sqrt{N} (U_{n, N}' - \bE[h])/\sigma_h$ to $\cN(0, 1)$. This also explains why
  assuming $\bE[h] = 0$  does not forgo any generality of the Berry-Esseen-type results developed for these convergences.
 
The random  $\hat{N}$ is also  preferred to the deterministic  $N$ as a normalizer  for statistical power. In fact, \citet[Corollary 1]{janson1984asymptotic} is originally stated for the weak convergence of the deterministically normalized $\breve{U}_{n, N}'  \equiv  N^{-1} \sum_{{\bf i}\in I_{n, m}} Z_{{\bf i}}h(\bX_{\bf i})$: Under exactly the  same asymptotic regimes and  assumptions  as  \thmref{weak_convergences} \emph{except} $\bE[h] = 0$, \citet[Corollary 1]{janson1984asymptotic}  states the three weak convergences in \thmref{weak_convergences}   with  $U_{n, N}'$ replaced by $(\breve{U}_{n, N}' - \bE[h])$ and $\sigma_h$ replaced by $\sqrt{\sigma_h^2 + (\bE[h])^2}$;  specifically, \thmref{weak_convergences}$(ii)$ and $(iii)$ become the weak convergences of 
\[
\frac{\sqrt{N}(\breve{U}_{n, N}' - \bE[h])}{\sqrt{\sigma_h^2 + (\bE[h])^2}}  \quad \text{ and } \quad 
\frac{\sqrt{n} (\breve{U}_{n, N}' - \bE[h])}{\sqrt{m^2 \sigma_g^2 + \alpha( \sigma_h^2 + (\bE[h])^2)}}
\]
 to $\cN(0, 1)$\footnote{Therefore, \citet[Corollary 1]{janson1984asymptotic} implies \thmref{weak_convergences} as $N/\hat{N} \longrightarrow_p 1$.}.  Since $(\bE[h])^2 >0$  when $\bE[h] \neq 0$, by requiring a larger standardizing factor, $\breve{U}_{n, N}'$ is asymptotically less powerful than $U_{n, N}'$ under the two regimes  of  \thmref{weak_convergences}$(ii)$ and $(iii)$ on $\alpha$; \citet[Remarks 2.1 \& 3.3]{chen2019randomized} made a similar remark.

\subsection{Random kernel incomplete U-statistics } \label{sec:random_kernel_ICU}

Our  present work  is partly stimulated by the applications of incomplete U-statistics  in uncertainty quantification for ensemble predictions, where each base learner is built on a subsample of the original data. Some  of these applications may employ base learners of a randomized nature, such as the trees in \citet{breiman2001random}'s random forests,  which randomly select the input covariates to split each node. To accommodate such applications,  the  kernels in \eqref{icu_def} have to be embellished by  additional random elements, to result in a \emph{random kernel} incomplete U-statistic of the form
 \begin{equation} \label{random_kernel_icu}
U_{n, N, \omega}' \equiv \frac{1}{\hat{N}}  \sum_{{\bf i}\in I_{n, m}} Z_{{\bf i}}h(\bX_{\bf i}; \omega_\boldi),
 \end{equation}
 where $\omega \equiv (\omega_\boldi)_{\boldi \in \Inm}$ are independent and identically distributed random parameters that are also independent of  $\cX$ and $\cZ$ \citep{mentch2016quantifying, peng2022rates}. Each    $h(\bX_{\bf i}; \omega_\boldi)$, while still symmetric in  its arguments $X_{i_1}, \dots, X_{i_m}$,   also depends on the realized value of its random parameter $\omega_\boldi$.   The  methodology in this paper can also be readily extended to prove B-E bounds for  \eqref{random_kernel_icu}; see \secref{discuss}.

 \section{The  normal convergences of $U_{n, N}'$} \label{sec:heuristics}

In what follows, we  let
\begin{equation} \label{non_negative_ustat_def}
 U_{h^2} \equiv {n \choose m}^{-1} \sum_{\boldi \in I_{n,m}} h^2 (\bX_\boldi)
 \text{ and } U_{|h|^3} \equiv  {n \choose m}^{-1} \sum_{\boldi \in I_{n,m}} |h|^3 (\bX_\boldi),
 \end{equation}
be the U-statistics whose \emph{non-negative} kernels $h^2$ and $|h|^3$ are induced by taking the square and cube of the absolute value of the original kernel $h$. Moreover,  as in many recent works on B-E bounds for incomplete U-statistics \citep{chen2019randomized, SongChenKato2019, peng2022rates, sturma2024testing}, much of the development in this paper revolves around a basic decomposition of $U_{n, N}'$, which is
\begin{equation} \label{basic_decomp}
U_{n, N}' = \frac{N}{\hat{N}} (\sqrt{(1 - p)}  B_n  +   U_n)  \text{ when } \hat{N} \neq 0,\footnote{ When $\hat{N} = 0$, $U_{n, N}'$ is simply zero as $Z_\boldi = 0$ for all $\boldi \in \Inm$.}
\end{equation}
for
\begin{equation}\label{Bn_def}
B_n \equiv  \frac{1}{N} \sum_{{\bf i} \in I_{n, m}} \frac{Z_{\bf i} - p}{\sqrt{1 - p}} h(\bX_{\bf i}).
\end{equation}
 Therefore, up the  multiplicative scaling factor $N/\hat{N}$,  $U'_{n, N}$ can be understand as a sum of the complete U-statistic $U_n$ and $\sqrt{(1- p)}B_n$.

We   first expand on the nature of the normal convergences in \thmref{weak_convergences}, which  can shed light on  the forms of their corresponding B-E bounds stated in \secref{main_results}.   To that end, it is instructive to  get a sense of the distributional behavior of $U_n$ and $B_n$.  Given \eqref{mean0_assumption} and a fixed $m$, it is well-known that
 \begin{multline} \label{convg_complete_ustat}
  \sqrt{n}U_n \text{ converges weakly to } \cN(0, m^2 \sigma_g^2) \text{ as } n \longrightarrow \infty \\
 \text{ under the non-degeneracy condition \eqref{non_degenerate_condition};}
 \end{multline}
  see \citet[Ch.4]{KB1994Ustat}. On the other hand, when \emph{conditioned on} the raw data $\cX$, $\frac{\sqrt{N}B_n}{\sqrt{U_{h^2}}} = \sum_{\boldi \in \Inm} \zeta_\boldi$ is a sum 
  of the conditionally independent random variables 
  \begin{equation} \label{xi_tilde_def}
\zeta_\boldi \equiv \frac{   (Z_\boldi - p)h(\bX_\boldi)}{\sqrt{U_{h^2}N(1 - p)}} \footnote{For now, we ignore that $U_{h^2}$ can be equal to zero, in which case $\zeta_\boldi$ becomes undefined.}
  \end{equation}
   indexed by $\boldi \in \Inm$, each having  the conditional mean
 \begin{equation}\label{cond_1st_moment_sqrtNBn}
 \bE[\zeta_\boldi  \mid \cX ] = 0, 
\end{equation}
and the conditional second and third absolute moments 
\begin{equation*} 
\bE[ \zeta_\boldi^2  \mid \cX ] =  \frac{ h^2(\bX_\boldi)}{U_{h^2}{n \choose m}}  \quad \text{ and }  \quad 
\bE[ |\zeta_\boldi|^3 \mid \cX ] =  \frac{  (1- 2p + 2p^2)|h|^3(\bX_\boldi)}{ U_{h^2}^{3/2}{n \choose m} \sqrt{N(1 - p)}};
\end{equation*}
we have used the calculation $\bE[|Z_\boldi - p|^3]
= p(1- p)( 1 - 2 p + 2p^2)$ in the latter.  
 Since $\sqrt{N}B_n/\sqrt{U_{h^2}}$ has conditional variance $1$, as in 
 \begin{equation}\label{cond_2nd_moment_sqrtNBn}
 \bE\bigg[\frac{NB_n^2}{U_{h^2}} \mid \cX\bigg] = \sum_{\boldi \in \Inm}\bE[ \zeta_\boldi^2  \mid \cX ]  = \sum_{\boldi \in \Inm}  \frac{ h^2(\bX_\boldi)}{U_{h^2} {n \choose m}} = 1,
 \end{equation}
 and the  conditional absolute third moments  of its summands add up to  
  \begin{equation} \label{cond_3rd_moment_sqrtNBn}
\sum_{\boldi \in \Inm}\bE\Big[ |\zeta_\boldi|^3   \mid \cX \Big]=  \frac{U_{|h|^3} ( 1 - 2 p + 2p^2)}{ U_{h^2}^{3/2}\sqrt{N (1- p)}} ,
  \end{equation}
the  Lyapunov's  central limit theorem \citep[p.362]{billingsley}   
suggests that
\begin{multline} \label{cond_convg_sqrtNBn}
\text{
 $\sqrt{N}B_n$ converges weakly to $\cN(0, U_{h^2})$
  conditional on $\cX$
  }\\
 \text{under the asymptotic regimes of \thmref{weak_convergences}$(ii)$ and $(iii)$,}
 \end{multline}
  because $\frac{1 - 2p + 2p^2}{\sqrt{N(1-p)}} \longrightarrow 0$
  as $n^d/N \longrightarrow \infty$ or $\alpha$ converges to a positive constant, and $U_{|h|}^3$ and $U_{h^2}$ respectively converges to $\bE[|h^3|]$ and $\sigma_h^2$ almost surely by the strong law of large number for U-statistics \citep[Theorem 3.1.2]{KB1994Ustat}\footnote{Here, we have implicitly assumed the absolute third moment $\bE[|h|^3]$ is finite.}. 
The conditional weak convergence  in \eqref{cond_convg_sqrtNBn} suggests that, marginally, 
\begin{multline} \label{convg_Bn}
 \sqrt{N}B_n \text{ is  approximately distributed as  } \cN(0, \sigma_h^2) \text{ for large $n$ and $N$}, \\
  \text{under the asymptotic regimes of \thmref{weak_convergences}$(ii)$ and $(iii)$.}
\end{multline} 

The observations in \eqref{convg_complete_ustat} and \eqref{convg_Bn}, along with the fact that
\begin{equation} \label{N_Nhat_ratio_convg}
\frac{N}{\hat{N}} \text{ converges  to $1$ in probability as } N \longrightarrow \infty, 
\end{equation}
 allow one to either prove or heuristically deduce the weak convergences  in \thmref{weak_convergences}$(i)-(iii)$  under their respective asymptotic regimes:

\subsection{The ``$N \gg n$'' regime, $\alpha \longrightarrow 0$}
\label{sec:Nggn}
By  \eqref{basic_decomp},   write 
\[
\sqrt{n}U_{n, N}'=  \frac{N}{\hat{N}} (  \sqrt{(1 - p)n}  B_n  +   \sqrt{n} U_n).
\] 
Because \eqref{cond_2nd_moment_sqrtNBn} implies $\bE[ NB_n^2] = \sigma_h^2$, $\sqrt{N} B_n$ is bounded in probability, which further implies $\sqrt{(1-p)n} B_n =  ((1-p)n/N)^{1/2} \sqrt{N} B_n$ converges to   $0$ in probability as $\alpha = n/N \longrightarrow 0$.  Therefore,  with \eqref{N_Nhat_ratio_convg}, it is easy to see that 
\begin{equation} \label{Nggn_approximation}
\sqrt{n}U_{n, N}'= \sqrt{n} U_n+ o_p(1),
\end{equation} 
and hence $\sqrt{n}U_{n, N}'$ converges weakly to the same limit as \eqref{convg_complete_ustat}, which establishes \thmref{weak_convergences}$(i)$ .

\subsection{The ``$N \ll n^d$'' regime, $n^d/N \longrightarrow \infty$} \label{sec:Nlln}
From the decomposition \eqref{basic_decomp} one can write 
\[
\frac{\sqrt{N}U_{n, N}'}{\sigma_h} =  \frac{N}{\hat{N}} \cdot  \frac{\sqrt{(1 - p)N}  B_n  +   \sqrt{N} U_n}{\sigma_h}.\]
First, as $n^d/N \longrightarrow \infty$,  $p= N/{n \choose m}$   converges  to $0$.  Second, $\sqrt{N}U_n/ \sigma_h$  converges to $0$ in probability, because the scaling by $\sqrt{N}$ is  smaller than the necessary scaling by $n^{d/2}$ under $n^d/N \longrightarrow \infty$ to render a non-vanishing  limiting variance for $U_n$\footnote{ Under \eqref{mean0_assumption},  $\bE[U_n^2]$ has the  order $O(n^{-d})$; see \citet[Sec 5.2.2, p.185-186, Lemma B]{serfling1980}. Therefore, a minimal scaling of order $O(n^{d/2})$  is required.}.  Therefore, together with \eqref{N_Nhat_ratio_convg}, we have
\begin{equation} \label{Nlln_approximation}
\frac{\sqrt{N}U_{n, N}'}{\sigma_h}  = \frac{\sqrt{N}B_n}{\sigma_h} + o_p(1),
\end{equation}
and the weak convergence in \thmref{weak_convergences}$(ii)$ can be deduced in light of \eqref{convg_Bn}.

\subsection{The ``$N \asymp n$'' regime, $\alpha$ tending to a positive constant}  \label{sec:Nasympn}
From \eqref{basic_decomp} one can write
\begin{equation} \label{U_prime_as_U_star}
\frac{\sqrt{n}U_{n, N}' }{\sigma} =   \frac{N}{\hat{N}} \bigg(\frac{\sqrt{(1 - p)n}  B_n }{\sigma} +   \frac{\sqrt{n}U_n}{\sigma}\bigg)  \approx \frac{ \alpha^{1/2}\sqrt{N}  B_n + \sqrt{n}U_n}{\sigma} ,\end{equation}
where the ``$\approx$'' above comes from  that $p \longrightarrow 0$ as $\alpha$ tends to a positive constant (because $m \geq 2$), as well as \eqref{N_Nhat_ratio_convg}. 
Note that under the zero-mean assumption in \eqref{mean0_assumption}, $U_n$ and $B_n$ are uncorrelated, because
\[
\bE[U_n B_n ] = \bE[U_n  \bE[ B_n \mid \cX]] = \bE\bigg[  \frac{U_n }{N} \sum_{{\bf i} \in I_{n, m}} \frac{\bE[Z_{\bf i}] - p}{\sqrt{1 - p}} h(\bX_{\bf i}) \bigg] =  0
\]
 given that $Z_\boldi$ has success probability $p$. 
Therefore, via \eqref{convg_complete_ustat} and \eqref{convg_Bn}, $\alpha^{1/2}\sqrt{N}  B_n + \sqrt{n}U_n$ is expected to be approximately normally distributed with mean $0$ and variance $\sigma^2 = \alpha \sigma_h^2 + m^2 \sigma_g^2$, which in turn deduces the  the weak convergence in \thmref{weak_convergences}$(iii)$.

In summary, under the three different asymptotic regimes, the relative dominance between $U_n$ and $B_n$ in the decomposition \eqref{basic_decomp}  lead to the approximations in  \eqref{Nggn_approximation}-\eqref{U_prime_as_U_star}, which in turn explain the weak convergences in  \thmref{weak_convergences}.

\section{The main results}  \label{sec:main_results}

 Since the normal convergences of $U_{n, N}'$ under the three regimes in \thmref{weak_convergences} are driven by either or both of the distributions of $U_n$ and $B_n$ as discussed in  \secref{heuristics}, one would intuitively expect their corresponding B-E bounds to capture features of the existing Berry-Esseen-type results for the convergences of $U_n$ and $B_n$ in  \eqref{convg_complete_ustat} and \eqref{cond_convg_sqrtNBn}. Firstly, under the non-degeneracy condition  \eqref{non_degenerate_condition}, the complete U-statistic $U_n$  has an established  B-E bound \citep{chen2007normal, peng2022rates} that  reads as
 \begin{multline} \label{complete_U_stat_BE}
 \sup_{z \in \bR}\bigg| P\bigg( \frac{\sqrt{n} U_n}{m \sigma_g} \leq z\bigg) - \Phi(z)\bigg| \leq \frac{6.1 \bE[|g|^3]}{\sqrt{n} \sigma_g^3} + (1 + \sqrt{2}) \bigg( \frac{m}{n} \bigg(\frac{\sigma_h^2}{m \sigma_g^2} - 1 \bigg) \bigg)^{1/2}.
 \end{multline}
Moreover,  the calculations  \eqref{cond_1st_moment_sqrtNBn}-\eqref{cond_3rd_moment_sqrtNBn} and the classical B-E bound  \citep{berry1941accuracy, esseen1942, shevtsova2010improvement} for a sum of independent variables suggests that
\begin{equation} \label{classical_BE_bdd}
\sup_{z \in \bR} \bigg|P\bigg(\frac{\sqrt{N} B_n }{\sqrt{U_{h^2}}} 
\leq z \mid \cX\bigg) - \Phi(z) \bigg| \leq 
\frac{0.56 U_{|h|^3} (1 - 2p + 2p^2)}{U_{h^2}^{3/2} \sqrt{N(1 - p)}}.
 \end{equation}

 We first state our B-E bound corresponding to the normal convergence in \thmref{weak_convergences}$(i)$ under the ``$N \gg n$'' regime. 

\begin{theorem}[Berry-Esseen theorem for  $U_{n ,N}'$, ``$N \gg n$''] \label{thm:BE_bdd_Nggn}
Let  $U_{n ,N}'$ be constructed as in \eqref{icu_def} with a symmetric kernel $h$ on $M^m$. Suppose that \eqref{mean0_assumption}, \eqref{non_degenerate_condition}  and  $2 \leq m < n/2$ hold. For an absolute constant $C >0$, 
\begin{multline*}
 \sup_{z \in \bR}\bigg| P\bigg(\frac{\sqrt{n} U_{n, N}'}{m\sigma_g} \leq z\bigg) - \Phi(z) \bigg|
 \leq \\
 C \bigg\{
  \frac{ \bE[|g|^3]}{\sqrt{n} \sigma_g^3}  + \bigg( \frac{m}{n} \bigg(\frac{\sigma_h^2}{m \sigma_g^2} - 1 \bigg) \bigg)^{1/2}  
  + \bigg( \frac{n (1-p)  \sigma_h^2 }{Nm \sigma_g^2}   \bigg)^{1/2} + \frac{1}{\sqrt{N}} \bigg\}.
\end{multline*}

\end{theorem}
As discussed in \secref{Nggn},  the asymptotic normality of  $\sqrt{n}U_{n, N}'$ is primarily driven by the complete U-statistic $\sqrt{n}U_n$ when $N \gg n$, so it should come as no surprise that the first two  terms in the bound have their roots in  the B-E bound   \eqref{complete_U_stat_BE}. The third term $( \frac{n (1-p)  \sigma_h^2 }{Nm \sigma_g^2} )^{1/2}$ accounts for the small remainder  $o_p(1)$ figuring in \eqref{Nggn_approximation}, and sensibly reflects that the normal approximation is  accurate only when $N$ is sufficiently large compared to $n$. The last term $1/\sqrt{N}$ stems from normalizing $U_{n, N}'$ with the more powerful $\hat{N}$ instead of $N$ (cf. \secref{sp_remark_Nhat}). In fact, 
\citet[Theorem 4]{peng2022rates} states almost exactly the same B-E bound as  \thmref{BE_bdd_Nggn}, except that our last term $1/\sqrt{N}$ is replaced by the  term   $N^{-1/2+ \epsilon}$ for an   arbitrary fudge factor $\epsilon \in (0, 1/2)$. We avoided such slackness by the variable censoring technique in our proof; see \secref{pf_preps}.

Our  B-E bounds for \thmref{weak_convergences}$(ii)$ and $(iii)$  have to be stated in terms of  projection functions with respect to the non-negative kernel $h^2$ of $U_{h^2}$ in \eqref{non_negative_ustat_def}. For each $r = 1, \dots, m$, first define the function $\tilde{\Psi}_r : M^r \rightarrow \bR$  in $r$ arguments by
\[
\tilde{\Psi}_r (x_1, \dots, x_r) \equiv \bE[h^2(x_1, \dots, x_r, X_{r+1}, \dots, X_m)] - \sigma_h^2,
\]
which is centered in the sense that $\bE[\tilde{\Psi}_r] \equiv \bE[\tilde{\Psi}_r (X_1, \dots, X_r)]= 0$; specifically for $\tilde{\Psi}_1$, we also let  $\Psi_1: M \rightarrow \bR_{\geq 0}$ be its non-negative uncentered version
\begin{equation*}
\Psi_1(x_1) \equiv \bE[h^2(x_1, X_{2}, \dots, X_m)].
\end{equation*}
Next, we  recursively define the \emph{Hoeffding projection kernels} $\pi_r(h^2) : M^r \rightarrow \bR$ of $h^2$ as
\[
\pi_ 1(h^2) (x_1) \equiv     \tilde{\Psi}_1 (x_1) \text{ for } r =1
\]
and 
\begin{multline*}
\pi_r(h^2) (x_1 \dots, x_r) \equiv  \tilde{\Psi}_r(x_1 \dots, x_r) - \sum_{s = 1}^{r-1} \sum_{1 \leq j_1 < \cdots < j_s\leq r} \pi_s(h^2) (x_{ j_1} , \dots, x_{ j_s})\\
\text{ for } r =2, \dots, m.
\end{multline*}
These projection functions give rise to the  \emph{Hoeffding decomposition} \citep[p.137]{de2012decoupling}
of $U_{h^2}$ that, upon re-centering and rescaling with $\sigma_h^2$, can be written as
\begin{equation}   \label{H_decomp_h2_ustat}
\frac{U_{h^2}- \sigma_h^2}{\sigma_h^2} 
= 
\sum_{i =1}^n \eta_i  - m+ \cR,
\end{equation}
where 
\begin{equation} \label{eta_def}
\eta_i \equiv \frac{m \Psi_1(X_i) }{n \sigma_h^2} \text{ for } i = 1, \dots, n
\end{equation}
are the leading terms, and 
\begin{equation} \label{H_decomp_remainder}
\cR \equiv  \sum_{r = 2}^m {m \choose r} {n \choose r}^{-1} \sum_{1 \leq j_1 < \dots < j_r \leq n}\frac{\pi_r(h^2) (X_{j_1}, \dots, X_{j_r})}{\sigma_h^2}
\end{equation}
is the degenerate lower-order remainder term in the decomposition. We  note that  each of the $\eta_i$'s has mean $\bE[\eta_i] = m/n$, giving rise to the centering of $\sum_{i=1}^n \eta_i$ by $m$ in \eqref{H_decomp_h2_ustat}. For our convenience later, we also define the  term
\begin{equation*}
\frakB_1 \equiv  \frac{\bE[|h|^3] (1 - 2p + 2p^2)}{ \sqrt{N(1-p)}\sigma_h^3}  + \exp \bigg( \frac{- n \sigma_h^6}{24 m(\bE[|h|^3])^2} \bigg).
\end{equation*}

We can now state our  B-E bound for  \thmref{weak_convergences}$(ii)$.
\begin{theorem}[Berry-Esseen theorem for $U_{n, N}'$,  ``$N \ll n^d$''] \label{thm:BE_bdd_Nlln}
Let  $U_{n ,N}'$ be constructed as in \eqref{icu_def} with a symmetric kernel $h$ on $M^m$. Suppose that \eqref{mean0_assumption}, $\bE[|h|^3] < \infty$  and  $2 \leq m < n/2$ hold.  Then, for an absolute constant $C >0$ and
\begin{equation*}
K_{n, m, d} \equiv {n \choose m}^{-1}  \sum_{r = d}^m {m -1 \choose r -1} {n - m \choose m - r},
\end{equation*}
\begin{equation}  \label{BE_bdd_Nlln_3rd_moment}
 \sup_{z \in \bR}\bigg| P\bigg(\frac{\sqrt{N} U_{n, N}'}{\sigma_h} \leq z\bigg) - \Phi(z) \bigg|
 \leq 
 C \bigg\{ 
\frakB_1 + 
  \frac{  \sqrt{N K_{n, m, d}}}{\sqrt{1-p}} +
  \frac{m^{3/2}\bE[\Psi_1^{3/2}]}{ \sqrt{n}\sigma_h^3}  + 
   \|\cR\|_{3/2}
\bigg\}.
\end{equation}
 If, in addition,  $\bE[h^4] < \infty$, then the last two terms can be replaced by $\frac{\sqrt{m \bE[(h^2 - \sigma_h^2)^2]}}{\sqrt{n}\sigma_h^2 }$ to give
\begin{equation}  \label{BE_bdd_Nlln_4th_moment}
 \sup_{z \in \bR}\bigg| P\bigg(\frac{\sqrt{N} U_{n, N}'}{\sigma_h} \leq z\bigg) - \Phi(z) \bigg|
 \leq \\
 C \bigg\{ 
\frakB_1  + 
 \frac{  \sqrt{NK_{n, m, d}}}{\sqrt{1-p}} +
\frac{\sqrt{m \bE[(h^2 - \sigma_h^2)^2]}}{ \sqrt{n}\sigma_h^2}
\bigg\}.
\end{equation}
\end{theorem}
 We note that   among  the summands of $K_{n, m, d}$, the term ${n \choose m}^{-1} {m -1 \choose d -1} {n - m \choose m - d}$ has the largest order with respect to $n$, which is $O(n^{- d})$. It is also possible to be more elaborate about the term $\|\cR\|_{2/3}$  in \eqref{BE_bdd_Nlln_3rd_moment}. By standard moment inequalities for  degenerate  U-statistics  \citep[Theorem 2.1.3]{KB1994Ustat}\footnote{Also see  \citet{Ferger1996}.}, we have the slightly cumbersome-looking inequality
\begin{equation} \label{complicated_expression}
\bE[|\cR|^{2/3}]\leq \sum_{r = 2}^m  {m \choose r}^{3/2} {n \choose r}^{-1/2}  \frac{2^{r/2}\bE[|\pi_r(h^2)|^{3/2}]}{\sigma_h^{3}}.
\end{equation}
Therefore, $\|\cR\|_{2/3}$ is seen to be of order at most $O(n^{-2/3})$. Next, we break down the interpretations of  the different parts in the B-E bounds as follows:
\begin{itemize}
\item 
$\frakB_1$: 
It has been seen in   \secref{Nlln} (cf. \eqref{Nlln_approximation}) that $\sqrt{N}B_n$ drives the distribution of $\sqrt{N} U_{n, N}'$ under the regime $n^d/N \longrightarrow \infty$. Therefore, the term $\frakB_1$ in \thmref{BE_bdd_Nlln}  captures some essential features of the B-E bound in \eqref{classical_BE_bdd} for the conditional distribution of $\sqrt{N}B_n$.   Because  \thmref{weak_convergences}$(ii)$ really concerns the  convergence of the \emph{marginal} distribution for $\sqrt{N}U_{n, N}'$,  the ratio $U_{|h|^3}/U_{h^2}^{3/2}$  in \eqref{classical_BE_bdd}, which depends on the data $\cX$,   
is now replaced by $\bE[|h|^3]/\sigma_h^3$, a ratio of 
 population moments that can be loosely thought of as its ``average'' given $\bE[U_{|h|^3}] = \bE[|h|^3]$ and  $\bE[U_{h^2}] = \sigma_h^2$. The extra  term $\exp \big ( \frac{- n \sigma_h^6}{24 m(\bE[|h|^3])^2} \big)$ captures the fact that $U_{|h|^3}/U_{h^2}^{3/2}$ can never be too large due to the exponentially decaying lower tail behavior of $U_{h^2}$ as a \emph{non-negative} U-statistic; see \eqref{exp_lower_bdd_app1} below.  
 
 The  factor $\frac{1 - 2p + 2p^2}{(1-p)^{1/2}}$ should also be understood with some care, because $p$  depends on $N$ and $n$ as per  \eqref{bernoulli_samp_prob}. Indeed, as $N$ approaches ${n \choose m}$, $p$ approaches $1$ and the factor  $\frac{1 - 2p + 2p^2}{(1-p)^{1/2}}$ blows up.
  However, our B-E bound in \thmref{BE_bdd_Nlln} is really  geared for the regime   $n^d/N \longrightarrow \infty$, in which case  $\lim_{n \rightarrow \infty}\frac{1 - 2p + 2p^2}{(1-p)^{1/2}}= 1$  because $p$   converges  to $0$ under the regime $n^d/N \longrightarrow \infty$ as discussed in \secref{Nlln}. 
 
 \item   $\frac{\sqrt{N K_{n, m, d}}}{\sqrt{1-p}}$:  Because $K_{n, m, d}$ has the order $O(n^{-d})$ with respect to $n$ as mentioned before, this term, stemming from the $o_p(1)$ term in \eqref{Nlln_approximation},  captures that as $N$ becomes too large relative to $n^d$, the accuracy of the normal approximation deteriorates.   For the special cases of $d = 1$ and $d = 2$, one can actually replace 
 $K_{n, m, d}$ by a neater term of the same order in $n$  in the B-E bounds of \thmref{BE_bdd_Nlln}. For $d=1$,  it is easy to see that  $K_{n, m, 1} = m/n$ by Vandermonde's identity. For $d = 2$, \citet[Ch. 10, p.285-286]{chen2010normal}
  has demonstrated that
 \begin{equation*}\label{Knm2_chen_shao_bdd}
 K_{n, m, 2} \leq \frac{m(m-1)^2}{ n(n-m+1)}.
 \end{equation*}
 We suspect that a bound in terms of $(n, m, d)$ of the same quality as the latter display   can be proved  for $K_{n, m, d}$ with  $d \geq 3$. However,  at the time of writing, it is unclear to us  how to generalize the arguments in  \citet[Ch. 10, p.285-286]{chen2010normal} and we leave that for future research. 
 
 \item $\frac{m^{3/2}\bE[\Psi_1^{3/2}]}{ \sqrt{n}\sigma_h^3}  + 
   \|\cR\|_{3/2}$  in \eqref{BE_bdd_Nlln_3rd_moment} and $\frac{\sqrt{m \bE[(h^2 - \sigma_h^2)^2]}}{ \sqrt{n}\sigma_h^2}$ in \eqref{BE_bdd_Nlln_4th_moment}: These bounding terms  stem from the  \emph{denominator remainder} $D_2$ defined in \eqref{D2_def}. Heuristically, to bridge between  the \emph{conditional} normal convergence in \eqref{cond_convg_sqrtNBn} and  the \emph{marginal} normal approximation in \eqref{convg_Bn}, one has to account for the estimation of marginal variance $\sigma_h^2$ by the conditional variance $U_{h^2}$, and the term $D_2$ precisely  emerges under our proof strategy for this purpose.

\end{itemize}

Lastly, we state our B-E bound for  \thmref{weak_convergences}$(iii)$. To that end, we also conveniently define the  bounding term 
\[
 \frakB_2 \equiv \frac{ \bE[|g|^3]}{\sqrt{n} \sigma_g^3} + \bigg(\frac{m}{n} \bigg(\frac{\sigma_h^2}{m \sigma_g^2} - 1 \bigg)\bigg)^{1/2}.
\]
\begin{theorem}[Berry-Esseen theorem for $U_{n, N}'$, ``$N \asymp n$'']\label{thm:main}
Let  $U_{n ,N}'$ be constructed as in \eqref{icu_def} with a symmetric kernel $h$ on $M^m$. Suppose that \eqref{mean0_assumption}-\eqref{non_degenerate_condition}, $\bE[|h|^3] < \infty$  and  $2 \leq m < n/2$ hold. Then, for an absolute constant $C >0$,
\begin{multline}  \label{BE_bdd_N_asymp_n_3rd_moment}
 \sup_{z \in \bR}\bigg| P\bigg(\frac{\sqrt{n} U_{n, N}'}{\sigma} \leq z\bigg) - \Phi(z) \bigg|
 \leq \\
  C \bigg\{
  \frakB_1 + \frakB_2
+
 \frac{  \sqrt{m}}{\sqrt{n(1-p)}} + \bigg(1 + \frac{\sqrt{mN}}{\sqrt{n(1-p)}}\bigg)\frac{m^{3/2}\bE[\Psi_1^{3/2}]}{ n^{1/2}\sigma_h^3}  + \|\cR\|_{3/2}
\bigg\}.
\end{multline}
If, in addition,  $\bE[h^4] < \infty$, then the last two terms can be replaced by $\frac{N}{n^2 (1-p)} + \frac{\sqrt{m \bE[(h^2 - \sigma_h^2)^2]}}{\sigma_h^2 \sqrt{n}}$ to give
\begin{multline}  \label{BE_bdd_N_asymp_n_4th_moment}
 \sup_{z \in \bR}\bigg| P\bigg(\frac{\sqrt{n} U_{n, N}'}{\sigma} \leq z\bigg) - \Phi(z) \bigg|
 \leq \\
  C \bigg\{   \frakB_1 + \frakB_2 +
 \frac{  \sqrt{m}}{\sqrt{n(1-p)}} + \frac{N}{n^2 (1-p)} +   \frac{ \sqrt{m \bE[(h^2 - \sigma_h^2)^2]}}{\sqrt{n}\sigma_h^2 }
\bigg\}.
\end{multline}

\end{theorem}
Again, the $\|\cR\|_{3/2}$ in \eqref{BE_bdd_N_asymp_n_3rd_moment} is no larger than the order $O(n^{-2/3})$ due to \eqref{complicated_expression}. Since the asymptotic normality of $U_{n, N}'$
  under the ``$N \asymp n$'' regime   can be viewed as the  distribution resulting from summing the two approximately independent and normal random variables $\sqrt{n}U_n$ and $\sqrt{N}B_n$ as discussed in \secref{Nasympn}, the two B-E bounds in  \eqref{complete_U_stat_BE} and \eqref{classical_BE_bdd}   intuitively suggest that a ``correct'' B-E bound 
for the marginal cumulative distribution of $U_{n, N}'$ should be of the rate $n^{-1/2}$, when $n$ and $N$ are of comparable order in such a way that  
\begin{equation} \label{n_asymp_N}
 \max\Big(\frac{n}{N} , \frac{N}{n}  \Big) < C
\end{equation}
for an absolute constant $C>0$. Indeed, when \eqref{n_asymp_N} holds, both  \eqref{BE_bdd_N_asymp_n_3rd_moment} and \eqref{BE_bdd_N_asymp_n_4th_moment} attain this natural rate. This  is clearly superior to the B-E bound stated  in \citet[Theorem 5]{peng2022rates} for  the convergence in \thmref{weak_convergences}$(iii)$, which doesn't achieve the rate $n^{-1/2}$ under \eqref{n_asymp_N}; their  bound contains some slackening  terms such as  $(m/n)^{1/3}$ and $N^{-1/2 + \epsilon}$ for an arbitrary $\epsilon \in (0, 1/2)$,  which arise from a less-than-refined proof methodology. In addition, they assume a finite sixth moment on the kernel $h$, while we only assume a finite third moment, a holy grail for establishing B-E bounds,  for \eqref{BE_bdd_N_asymp_n_3rd_moment} and a finite fourth moment for the cleaner-looking \eqref{BE_bdd_N_asymp_n_4th_moment}.

We also summarize the different parts of the bounds in \thmref{main} below:
\begin{itemize}
\item $\frakB_1 + \frakB_2$:  $\frakB_2$ stems from the B-E bound \eqref{complete_U_stat_BE} for the complete U-statistic $\sqrt{n}U_n$, and similarly to what is explained  after \thmref{BE_bdd_Nlln}, $\frakB_1$ can be considered a B-E bound for the marginal distribution of $\sqrt{N} B_n$. Together, they account for the combined sum-of-independent-terms and   complete-U-statistic behaviors of $\sqrt{n}U_{n, N}'$ in light of the decomposition  \eqref{U_prime_as_U_star} under the ``$N \asymp n$'' regime.

\item $\Big(1 + \frac{\sqrt{mN}}{\sqrt{n(1-p)}}\Big)\frac{m^{3/2}\bE[\Psi_1^{3/2}]}{ n^{1/2}\sigma_h^3}  + \|\cR\|_{3/2}
$ in \eqref{BE_bdd_N_asymp_n_3rd_moment} 
and $\frac{ \sqrt{m \bE[(h^2 - \sigma_h^2)^2]}}{\sqrt{n}\sigma_h^2 }$ in \eqref{BE_bdd_N_asymp_n_4th_moment}: Analogous to $\frac{m^{3/2}\bE[\Psi_1^{3/2}]}{ \sqrt{n}\sigma_h^3}  + 
   \|\cR\|_{3/2}$  in \eqref{BE_bdd_Nlln_3rd_moment}  and $\frac{\sqrt{m \bE[(h^2 - \sigma_h^2)^2]}}{ \sqrt{n}\sigma_h^2}$ in \eqref{BE_bdd_Nlln_4th_moment}, these terms  also stem from the \emph{denominator remainder} $D_2$ defined in \eqref{D2_def} and account for  the estimation of $\sigma_h^2$ by $U_{h^2}$.

\item  Other terms, such as $ \sqrt{m}/\sqrt{n(1-p)}$,  stem from various error terms that emerge during the proof of \thmref{main}.
\end{itemize}

\subsection{Preparation for the proofs}  \label{sec:pf_preps}

The proofs of \thmsref{BE_bdd_Nggn}-\thmssref{main} all adopt a handy technique called \emph{variable censoring}, recently introduced in  \citet{leung2024another} as a refinement of  the variable truncation technique for proving B-E bounds with Stein's method. For a given real-valued random variable $Y$, its censored version with respect to two  constants $a, b \in \bR \cup \{-\infty, \infty\}$, where $a \leq b$, is defined as 
\[
 \bar{Y} \equiv a I(Y < a) + Y I(a \leq Y \leq b)+bI(Y>b),
\] 
i.e. once $Y$ goes beyond the range $[a, b]$, its censored version $\bar{Y}$ can only be either $a$ or $b$ depending on which direction the value of $Y$ goes; that $a$ or/and $b$ are allowed to be infinity means that $Y$ may only be  censored in  one direction, or not be censored at all.
This technique has the following apparent but   useful properties for establishing bounds:

\begin{property}[Properties of variable censoring] \label{property:censoring_property}
Let $Y$ and $Z$ be any two real-value variables. The following facts hold:
 
 \begin{enumerate}[(i)]
 \item  Suppose, for some  $a, b \in \bR \cup \{-\infty, \infty\}$ with $a \leq b$, we define
\[
\bar{Y} \equiv a I(Y < a) + Y I(a \leq Y \leq b)+bI(Y>b)
\] and 
\[
\bar{Z}  \equiv a I(Z < a) + Z I(a \leq Z \leq b)+bI(Z>b);
\]
then it must be that 
$
|\bar{Y} - \bar{Z} | \leq |Y- Z|.
$
\item
If $Y$ is a non-negative random variable, then it must also be true that 
\[
Y I(0 \leq Y \leq b)+bI(Y>b) \leq Y \text{ for any } b \in (0, \infty), 
\]
i.e., the upper-censored version of $Y$ is always no larger than $Y$ itself.
 \end{enumerate}
\end{property}

Throughout, we will use $\Omega$ to denote the entire underlying probability space. Moreover, we will define
the  (upper-and-lower) censored version of $\hat{N}$, 
\begin{equation}\label{N_censored}
\hat{N}_c \equiv \hat{N} I\Big(\hat{N} \in \Big[\frac{N}{2} ,\frac{3N}{2} \Big]\Big) + \frac{3N}{2}I\Big(\hat{N} > \frac{3N}{2}\Big) + \frac{N}{2} I\Big(  \hat{N}< \frac{N}{2}\Big),
\end{equation}
and 
 the event 
\begin{equation*} 
\cE_\cZ \equiv \bigg \{ \hat{N} \in \Big[ \frac{N}{2}, \frac{3N}{2} \Big]\bigg\}
\end{equation*}
that only depends on the  Bernoulli samplers in $\cZ$; note that  
\begin{equation} \label{censored_N_equal_N}
\hat{N}_c  = \hat{N} \text{ on } \cE_\cZ. 
\end{equation}
In addition,  regarding $\hat{N}$, we have the inequalities
 \begin{equation} \label{Ncminus_N_bdd}
 |\hat{N}_c - N| \leq  |\hat{N} - N|,
 \end{equation}
which  is true because it  must be that $ |\hat{N}_c - N| = \frac{N}{2}$ and $  |\hat{N} - N| \geq \frac{N}{2}$ if $\hat{N} \not \in [\frac{N}{2}, \frac{3N}{2}]$,
\begin{equation}\label{Nhat_1_norm_bddb}
\bE[| \hat{N} - N| ] 
\leq \|\hat{N} - N\|_2 = \sqrt{ \sum_{{\bf i} \in I_{n, m}}p(1 -p)  } = \sqrt{N (1 - p)} \leq \sqrt{N},
\end{equation}
as well as Bernstein's inequality \citep[Lemma 2.2.9]{VDVW} 
\begin{equation} \label{bernstein1}
P ( \Omega \backslash \cE_\cZ ) = P\Big(\Big|\frac{\hat{N}}{N} - 1\bigg|> \frac{1}{2}\Big)
\leq 2 \exp \Big( - \frac{3N}{28}\Big).\end{equation}
Later, for the proofs of  \thmsref{BE_bdd_Nlln} and \thmssref{main} in \secref{pf_Nlln_Nasympn}, for each $\boldi \in \Inm$, we shall define 
\begin{equation} \label{Nhat_leave_one_out}
\hat{N}^{(\boldi)} \equiv \sum_{\boldj \in \Inm ,  \boldj \neq \boldi} Z_\boldj,
\end{equation} 
the the leave-one-out version of $\hat{N}$, and  
censor it as
\begin{equation}\label{Nhatc_leave_one_out}
\hat{N}_c^{(\boldi)} \equiv \hat{N}^{(\boldi)}I\bigg(\hat{N}^{(\boldi)} \in \Big[\frac{N}{2} ,\frac{3N}{2} \Big]\bigg) + \frac{3N}{2}I\bigg(\hat{N}^{(\boldi)} > \frac{3N}{2}\bigg) + \frac{N}{2} I\bigg(  \hat{N}^{(\boldi)}< \frac{N}{2}\bigg).
\end{equation}

\section{Proof of \thmref{BE_bdd_Nggn}, the case of ``$N \gg n$''} \label{sec:pf_Nggn_case}
In this section, we will use  $\sum_{\boldj \in I_{n,m}, i \not\in \boldj}$ to denote a summation over the ${n-1  \choose m}$ index vectors $\boldj = (j_1, \dots, j_m)$ in $\Inm$ for which $j_k \neq i$ for any $k =1, \dots, m$. 
In the same spirit, we will use $\sum_{\boldj \in I_{n,m}, i \in \boldj}$ to denote a summation over the ${n -1 \choose m-1}$ index vectors $\boldj = (j_1, \dots, j_m)$ in $\Inm$ for which $j_k = i$ for some $k =1, \dots, m$. By letting
\[
\bar{h}_m (x_1, \dots, x_m) \equiv h(x_1, \dots, x_m) - \sum_{i = 1}^m g(x_i)
\]
we first define the degenerate U-statistic
\[
U_{\bar{h}_m} \equiv  {n \choose  m}^{-1} \sum_{\boldi \in \Inm}\bar{h}_m({\bf X}_\boldi).
\]
For each $i \in [n]$, 
 we also define the ``leave-one-out'' version of $U_{\bar{h}_m}$ and $B_n$, 
\[
U_{\bar{h}_m}^{(i)} \equiv  {n \choose  m}^{-1}\sum_{\boldj \in I_{n,m}, i \not\in \boldj} \bar{h}_m({\bf X}_\boldj)
\]
and
\[
B_n^{(i)} \equiv \frac{1}{N} \sum_{\boldj \in I_{n,m}, i \not\in \boldj} \frac{Z_\boldj - p}{\sqrt{1 - p}} h(\bX_\boldj).
\]
Using the arguments in \citet[p.284-287]{chen2010normal}\footnote{Combining $(10.74)$, $(10.75)$, $(10.80)$, $(10.81)$ in \citet[Ch. 10, Appendix]{chen2010normal} gives the first inequality in \eqref{Ubarhm_sq_estimate}; combining slightly modified forms of $(10.74)$, $(10.75)$,$(10.80)$ and $(10.82)$ in \citet[Ch. 10, Appendix]{chen2010normal} gives the first inequality in \eqref{Ubarhm_leave_one_out_sq_estimate}. }  and $m < n/2$, we have the following second moment estimates
\begin{equation} \label{Ubarhm_sq_estimate}
\frac{n\bE[U_{\bar{h}_m}^2] }{m^2 \sigma_g^2}\leq \frac{(m-1)^2}{m(n-m+1)} \bigg(\frac{\sigma_h^2}{m \sigma_g^2} -1\bigg) \leq \frac{2m}{n} \bigg(\frac{\sigma_h^2}{m \sigma_g^2} -1\bigg)
\end{equation}
and
\begin{equation} \label{Ubarhm_leave_one_out_sq_estimate}
\frac{n \bE[(U_{\bar{h}_m} - U_{\bar{h}_m}^{(i)} )^2]}{m^2 \sigma_g^2} \leq \frac{2(m-1)^2}{n m(n - m+1)}  \bigg(\frac{\sigma_h^2}{m \sigma_g^2} -1\bigg) \leq \frac{ 4m}{n^2}\bigg(\frac{\sigma_h^2}{m \sigma_g^2} -1\bigg) .
\end{equation}
Moreover,   we have
\begin{equation}\label{exp_Bn2}
\bE[B_n^2] = \frac{ \sigma_h^2 }{N}
\end{equation}
from \eqref{cond_2nd_moment_sqrtNBn} and
 \begin{multline} \label{Bn_minus_Bni_2nd_moment}
\bE[(B_n - B_n^{(i)})^2] = \bE\bigg[  \bigg(\frac{1}{N} \sum_{\boldj \in I_{n,m}, i \in \boldj} \frac{Z_\boldj - p}{\sqrt{1 - p}} h(\bX_\boldj) \bigg)^2 \bigg] \\
= \frac{1}{N^2 (1-p)} \sum_{\boldj \in I_{n,m}, i \in \boldj}  \bE[(Z_\boldj - p)^2 h^2(\bX_\boldj)  ] = \frac{{n -1 \choose m-1}\sigma_h^2}{N {n \choose m}} = \frac{m \sigma_h^2}{Nn}.
 \end{multline}

To prove  \thmref{BE_bdd_Nggn}, we write
\begin{multline*}
\bigg|P\bigg(\frac{\sqrt{n} U_{n, N}'}{m\sigma_g} \leq z \bigg) - \Phi(z) \bigg| \leq 
  \bigg|P\bigg(\frac{N}{\hat{N}_c} \cdot \frac{\sqrt{n} (U_n + \sqrt{(1-p)} B_n)}{m\sigma_g }  \leq z \bigg) - \Phi(z) \bigg| + P(\Omega \backslash \cE_\cZ)
\end{multline*}
by virtue of the decomposition in \eqref{basic_decomp}. 
So,  it suffice to prove the B-E bound  
\begin{multline} \label{BE_bound_our_nonlinear_stat}
\bigg|P\bigg(\frac{N}{\hat{N}_c} \cdot \frac{\sqrt{n} (U_n + \sqrt{(1-p)} B_n)}{m\sigma_g }  \leq z \bigg) - \Phi(z) \bigg| \leq \\
C \bigg\{\frac{ \bE[|g|^3]}{\sqrt{n} \sigma_g^3}  + \bigg( \frac{m}{n} \bigg(\frac{\sigma_h^2}{m \sigma_g^2} - 1 \bigg) \bigg)^{1/2}  
  + \bigg( \frac{n (1-p)  \sigma_h^2 }{Nm \sigma_g^2}   \bigg)^{1/2} + \frac{1}{\sqrt{N}} \bigg\}
\end{multline}
in light of \eqref{bernstein1}, as $\exp ( - 3N/28) \leq C/\sqrt{N}$. To that end, we will invoke an  existing result on Berry-Esseen bounds for so-called nonlinear statistics from  \citet[Theorem 10.1]{chen2010normal}:

\begin{lemma}[B-E bound for nonlinear statistics]\label{lem:chen_shao_2011_thm10_1} Let $\xi_1, \dots, \xi_n$ be independent random variables satisfying $\bE[\xi_i] = 0$ for all $i = 1, \dots, n$ and $\sum_{i=1}^n \bE[\xi_i^2] = 1$, and define $S \equiv \sum_{i=1}^n \xi_i$. Moreover, let $\Delta$ be any random variable on the same probability space as $\xi_1, \dots, \xi_n$. Then
\[
\sup_{z \in \bR} |P(S + \Delta \leq z) - P(S \leq z)| \leq 2 \sum_{i=1}^n \bE[\xi_i^3] + \bE[|S \Delta|] + \sum_{i=1}^n \bE[|\xi_i(\Delta - \Delta^{(i)})|],
\]
where $\Delta^{(i)}$ is any random variable on the same probability space such that $\xi_i$ and $(S - \xi_i, \Delta^{(i)})$ are independent, for each $i = 1, \dots n$. 
\end{lemma}
Now, we will define 
\begin{multline}
\Delta_1  \equiv \frac{\sqrt{n}(U_{\bar{h}_m}+  \sqrt{(1-p)}B_n )}{m \sigma_g} \text{ and }  \\
\Delta_2  \equiv  \bigg(\frac{N - \hat{N}_c}{\hat{N}_c} \bigg) \cdot \frac{\sqrt{n} (U_n + \sqrt{(1-p)} B_n)}{m\sigma_g },\\
\text{ and } \Delta \equiv \Delta_1 + \Delta_2
\end{multline}
as well as writing
\begin{equation} \label{write_in_terms_of_Delta1_and_Delta2}
\frac{N}{\hat{N}_c} \cdot \frac{\sqrt{n} (U_n + \sqrt{(1-p)} B_n)}{m\sigma_g }  
= \frac{\sum_{i=1}^n g(X_i)}{\sigma_g\sqrt{n}} +  \Delta_1+ \Delta_2.
\end{equation}
via the  decomposition
\begin{equation} \label{light_H_decomp}
U_n = \frac{m \sum_{i=1}^n g(X_i)}{n} +U_{\bar{h}_m}.
\end{equation}
For each $i \in [n]$, we also define the leave-one-out versions of $\Delta_1$ and $\Delta_2$
\begin{multline*}
\Delta_1^{(i)} \equiv \frac{\sqrt{n}(U_{\bar{h}_m}^{(i)}+  \sqrt{(1-p)}B_n^{(i)} )}{m \sigma_g} \text{ and }\\
\Delta_2^{(i)} \equiv \bigg(\frac{N - \hat{N}_c}{\hat{N}_c} \bigg) \cdot \frac{\sqrt{n} (\frac{m \sum_{j \neq i}g(X_j)}{n} + U_{\bar{h}_m}^{(i)} + \sqrt{(1-p)} B_n^{(i)})}{m\sigma_g }
\end{multline*}
that leave out any terms involving $X_i$, 
where we have used \eqref{light_H_decomp} in defining $\Delta_2^{(i)}$. 
To prove a bound on \eqref{BE_bound_our_nonlinear_stat},  under the non-degeneracy condition \eqref{non_degenerate_condition}, we will now apply \lemref{chen_shao_2011_thm10_1} by taking
\begin{equation} \label{using_chen_shao_lemma}
\xi_i = \frac{g(X_i)}{\sigma_g \sqrt{n}}, \quad \Delta = \Delta_1 + \Delta_2, \quad \Delta^{(i)} = \Delta^{(i)}_1 + \Delta^{(i)}_2,
\end{equation}
as well as using the following estimates regarding the remainder terms proven in \appref{pf_delta_estimate_Nggn}:
\begin{lemma}[Estimates regarding $\Delta_1$ and $\Delta_2$]\label{lem:Delta_estimates_Nggn}
The following estimates are true:
\begin{enumerate}[(i)]
\item
 $ \bE\Big[ \Big| \frac{\sum_{i=1}^n g(X_i)}{\sigma_g\sqrt{n}}  \cdot \Delta_1 \Big|\Big] \leq \Big(\frac{2m}{n} \Big( \frac{\sigma_h^2}{m \sigma_g^2} - 1\Big) \Big)^{1/2} +  \Big(\frac{n(1 - p)\sigma_h^2}{Nm \sigma_g^2}\Big)^{1/2}$;
\item 
$
\bE\Big[ \Big| \frac{g(X_i)}{\sigma_g\sqrt{n}} \cdot( \Delta_1- \Delta_1^{(i)}) \Big| \Big]
\leq \frac{1}{\sqrt{n}} \Big(  \frac{ 2\sqrt{m}}{n} \Big(\frac{\sigma_h^2}{m \sigma_g^2} -1\Big)^{1/2} + 
  \frac{\sigma_h}{\sigma_g}  \sqrt{\frac{(1-p)  }{Nm}}\Big)
$;
\item 
$\bE\Big[ \Big| \frac{\sum_{i=1}^n g(X_i)}{\sigma_g\sqrt{n}}  \cdot \Delta_2 \Big|\Big] 
\leq \frac{2}{\sqrt{N}} \Big(1 + \Big(\frac{2m}{n} \Big( \frac{\sigma_h^2}{m \sigma_g^2} - 1\Big) \Big)^{1/2} +  \Big(\frac{n (1 - p)\sigma_h^2}{Nm \sigma_g^2}\Big)^{1/2} \Big)$; 
\item 
$\bE\Big[ \Big| \frac{g(X_i)}{\sigma_g\sqrt{n}} \cdot( \Delta_2- \Delta_2^{(i)}) \Big| \Big]
 \leq \frac{2}{\sqrt{nN}} \cdot
\Big(  \frac{1}{\sqrt{n}}+     \frac{ 2\sqrt{m}}{n} \Big(\frac{\sigma_h^2}{m \sigma_g^2} -1\Big)^{1/2}+ \frac{\sigma_h}{\sigma_g}  \sqrt{\frac{(1-p)  }{Nm}}  \Big)$.
\end{enumerate}

\end{lemma}
%
 
Combining \eqref{using_chen_shao_lemma} and  \lemsref{chen_shao_2011_thm10_1}-\lemssref{Delta_estimates_Nggn}  proves \eqref{BE_bound_our_nonlinear_stat}, and hence \thmref{BE_bdd_Nggn}.

\section{Proofs of \thmsref{BE_bdd_Nlln} and \thmssref{main}, the cases of  ``$N \ll n^d$'' and ``$N \asymp n$''} \label{sec:pf_Nlln_Nasympn}

The proofs of \thmsref{BE_bdd_Nlln} and \thmssref{main} are placed under the same umbrella, because  both  employ two crucial ingredients: \begin{inparaenum}[(i)]  \item an  exponential  lower tail bound  for the non-negative U-statistic $U_{h^2}$ in \eqref{non_negative_ustat_def}, and \item a framework of using Stein's method to prove B-E bounds   for the so-called ``Studentized nonlinear statistics''. \end{inparaenum}

Our first ingredient is the lower tail probability bound 
\begin{equation} \label{exp_lower_bdd_app1}
 P(U_{h^2}  \leq \sigma^2/2 ) \leq 1.05 \exp \Big ( \frac{- n \sigma_h^6}{24 m(\bE[|h|^3])^2} \Big) \text{ when }  2m < n,
\end{equation}
which can be established by  \lemref{chernoff_lower_tail_bdd_U_stat} below.
\begin{lemma}[Exponential lower tail bound for U-statistics with non-negative kernels, \citet{leungshao2024nonuniform}] \label{lem:chernoff_lower_tail_bdd_U_stat}
Suppose  $U_\kappa = {n \choose m}^{-1} \sum_{{\bf i} \in \Inm} \kappa(X_{\bf i})$ is a U-statistic, and $\kappa : M^m \longrightarrow \bR_{\geq 0}$ is a symmetric kernel  of degree $m$ that can only take non-negative values, with the property that $0 <\bE[\kappa^\ell] < \infty$ for some $\ell \in (1, 2]$. Then for $0 < t \leq  \bE[\kappa]$,
\[
P(U_\kappa \leq t) \leq \exp\left(\frac{- [ n/m ] (\ell-1)  (\bE[\kappa]  -t)^{\ell/(\ell-1)}}{\ell (\bE[\kappa^\ell])^{1/(\ell-1)} }\right),
\]
where $[n/m]$ is defined as the greatest integer no larger  than $n/m$.
\end{lemma}

Note that \lemref{chernoff_lower_tail_bdd_U_stat} immediately gives $P(U_{h^2}  \leq \sigma^2/2 ) \leq \exp \Big ( \frac{- [n/m] \sigma_h^6}{24(\bE[|h|^3])^2} \Big)$ by taking $\ell = 3/2$ and $t = \sigma_h^2/2$, but when $2m< n$, one can further simplify it with
\[
\exp \Big ( \frac{- [n/m] \sigma_h^6}{24(\bE[|h|^3])^2} \Big) \leq \exp \Big ( \frac{- (n/m-1) \sigma_h^6}{24(\bE[|h|^3])^2} \Big) < 1.05\exp \Big ( \frac{- n \sigma_h^6}{24 m(\bE[|h|^3])^2} \Big),
\]
to give \eqref{exp_lower_bdd_app1}, where the last inequality uses $\sigma_h^6 \leq (\bE[|h|^3])^2$ and $\exp(1/24) < 1.05$.
The exponential decay in $n$ of the bound \eqref{exp_lower_bdd_app1}  implies that the event 
\begin{equation}\label{cX_events_def}
\cE_{1, \cX} \equiv \bigg\{ U_{h^2} > \frac{\sigma_h^2 }{2} \bigg\}
\end{equation}
has a dominating probability, a fact we will  leverage.  For the rest of this paper, with a slight abuse of notation we will use $\cX \in \cE_{1, \cX}$ to indicate  the realized value of $\cX$ is such that  $U_{h^2} > \sigma_h^2/2$. Moreover, we will lower censor $U_{h^2}$  as
\begin{equation} \label{Uhsq_lower_censored}
\bar{U}_{h^2} \equiv  \frac{\sigma_h^2}{2}I \Big(U_{h^2} \leq \frac{\sigma_h^2}{2}\Big) + U_{h^2} I\Big(U_{h^2} > \frac{\sigma_h^2}{2}\Big),
\end{equation}
which has the property that
\begin{equation} \label{Uh2_equal}
\bar{U}_{h^2}= U_{h^2} \text{ on }\cE_{1, \cX}.
\end{equation}

 Now we describe the second  ingredient. Let 
\begin{equation} \label{xi_def}
\xi_\boldi \equiv \frac{ (Z_\boldi - p)h(\bX_\boldi)}{(\bar{U}_{h^2}N(1 - p))^{1/2}},
\end{equation}
whose first moments conditional on the raw data $\cX$ can be  computed as
\begin{equation} \label{cX_cond_moment_properties_of_xi}
\bE[\xi_\boldi \mid \cX] = 0, \quad 
\bE[ \xi_\boldi ^2 \mid \cX ] =  \frac{ h^2(\bX_\boldi)}{\bar{U}_{h^2} {n \choose m}} \text{ and } \bE [ |\xi_\boldi|^3 \mid \cX ] =  \frac{  (1- 2p + 2p^2)|h(\bX_\boldi)|^3}{ {n \choose m}\bar{U}_{h^2}^{3/2} \sqrt{N(1 - p)}}.
\end{equation}
Note the similarity between the forms of $\xi_\boldi$ and $\zeta_\boldi$  in \eqref{xi_tilde_def}; however, $\xi_\boldi$ is always well-defined because, unlike $U_{h^2}$, the $\bar{U}_{h^2}$ in its denominator is always greater than $1/2$, so division by zero is not an issue. Now, for any  \emph{arbitrary} random variable 
\begin{equation} \label{generic_D1}
D_1  =  D_1 ( \cZ; \cX)
\end{equation}
that is a function of both $\cX$ and $\cZ$ 
 and 
\begin{equation} \label{D2_def}
D_2\equiv - \frac{\sigma_h^2}{\bar{U}_{h^2}} \Big(\frac{U_{h^2}}{\sigma_h^2} - 1 \Big) ,
\end{equation}
we 
 define
 the generic random variable  
\begin{equation}\label{SN_nonlinear_form}
T_{SN}
\equiv \frac{\sum_{\boldi \in I_{n, m}} \xi_\boldi + D_1}{\sqrt{1 + D_2}},
\end{equation}
which  will serve as a proof device for establishing \thmsref{BE_bdd_Nlln} and \thmssref{main}. In particular, for a generic threshold 
\begin{equation} \label{generic_frak_z}
\frakz_\cX = \frakz(\cX) \in \bR
\end{equation}
 which may possibly depend on the data $\cX$, we will 
leverage the following general Berry-Esseen bound for the   absolute difference $\big| P( T_{SN}  \leq \frakz_\cX)- 
\bE[\Phi (\frakz_\cX)] \big|$.
\begin{lemma}[B-E bound for  $|P( T_{SN}  \leq \frakz_\cX) - \mathbb{E}  \mbox{[}\Phi(\frakz_\cX)  \mbox{]}  |$ ] \label{lem:prelim_BE}
Suppose that \eqref{mean0_assumption}, $\bE[|h|^3] < \infty$  and  $2 \leq m < n/2$ hold. Then  for any $ \frakz_\cX \in \bR$ that possibly depends on the data $\cX$, there exists  an absolute constant $C >0$ such that
\begin{equation*} \label{prelim_BE_bdd}
|P( T_{SN} \leq \frakz_\cX ) - \bE[\Phi(\frakz_\cX)]|
\leq C\bigg( T(D_1) + T(D_2) + P(\Omega \backslash \cE_{1, \cX}  ) + \frac{\bE[|h|^3] (1 - 2p + 2p^2)}{\sigma_h^3 \sqrt{N(1 - p)}}  \bigg)  ,
\end{equation*}
for 
\[
T(D_1) \equiv P(|D_1|>1/2)+ \bE\Bigg[ \sqrt{\bE[\bar{D}_1^2 | \cX]} +\frac{\sum_{\boldi \in I_{n, m}} |h(\bX_\boldi)|\sqrt{\bE[ (\bar{D}_1- \bar{D}_1^{(\boldi)})^2 | \cX] }}{\sigma_h\sqrt{ {n \choose m}}} \Bigg]
\]
and
\[
T(D_2) \equiv P\Big(|D_2|>1/2\Big)+ \bE[\bar{D}_2^2 ] +  \Big| \bE\Big[\frakz_\cX \bar{D}_2  f_{\frakz_\cX} \Big(\sum_{\boldi \in I_{n, m}} \barxi_\boldi \Big) \Big]\Big|,
\]
where  we have defined
\begin{equation*}\label{Dk_censored}
\bar{D}_k \equiv  D_k I\Big(-\frac{1}{2}\leq D_k \leq \frac{1}{2}\Big)- \frac{1}{2} I\Big(D_k  < -\frac{1}{2}\Big)+ \frac{1}{2 }I\Big(D_k >\frac{1}{2}\Big)   \text{ for } k = 1, 2,
\end{equation*}
\begin{equation}\label{xi_censored}
\barxi_\boldi \equiv \xi_\boldi I(|\xi_\boldi| \leq 1) + 1 I(\xi_\boldi > 1) - 1 I(\xi_\boldi <- 1),
\end{equation}
and
\begin{equation} \label{Stein_solution_random_threshold}
f_{\frakz_\cX} (w) \equiv
  \begin{cases} 
  \sqrt{2 \pi} e^{w^2/2} \Phi(w)\bar{\Phi}(\frakz_\cX) &  w \leq \frakz_\cX\\
 \sqrt{2 \pi} e^{w^2/2} \Phi(\frakz_\cX) \bar{\Phi}(w)        &  w > \frakz_\cX
     \end{cases};
\end{equation}
moreover, 
\begin{equation} \label{D1_leave_one_out}
D_1^{(\boldi)}= D_1^{(\boldi)}(  \cZ\backslash \{Z_\boldi\}; \cX) \text{ is any function in  } \cX \text{ and }\cZ \text{ except }Z_\boldi
\end{equation}
and
\begin{equation} \label{barD1_leave_one_out}
\bar{D}_1^{(\boldi)} \equiv D_1^{(\boldi)} I\Big(-\frac{1}{2}\leq D_1^{(\boldi)} \leq \frac{1}{2}\Big)- \frac{1}{2} I\Big(D_1^{(\boldi)}  < -\frac{1}{2}\Big)+ \frac{1}{2 }I\Big(D_1^{(\boldi)} >\frac{1}{2}\Big),
\end{equation}
is its censored version.

\end{lemma}

We remark that  \eqref{Stein_solution_random_threshold} is  simply the typical solution  \citep[p.14]{chen2010normal}
\begin{equation} \label{Stein_solution}
f_z (w) \equiv
  \begin{cases} 
  \sqrt{2 \pi} e^{w^2/2} \Phi(w)\bar{\Phi}(z) &  w \leq z\\
 \sqrt{2 \pi} e^{w^2/2} \Phi(z) \bar{\Phi}(w)        &  w > z
     \end{cases}
\end{equation}
to the normal Stein equation  \citep{stein1972bound}
\begin{equation*} \label{steineqt}
I(w \leq z) - \Phi(z) = \frac{\partial f_z}{\partial w}(w) - wf_z(w),
\end{equation*}
but with the fixed $z \in \bR$  replaced by the possibly random $\frakz_\cX$ to accommodate our proof later. 
\lemref{prelim_BE}  is  established at some length in \appref{pf_prelim_BE} via  adapting the  arguments from  \citet{leung2024another} while conditioning on a given data realization $\cX \in \cE_{1, \cX}$.  Importantly,  \citet{leung2024another} fully  formalizes a paradigm of using Stein's method to prove B-E bounds  for \emph{Studentized nonlinear statistics}: For $n$ real-valued deterministic functions $t_1(\cdot), \dots, t_n(\cdot)$ and a sample of independent random variables $G_1, \dots, G_n \in \bR$ such that $\bE[t_i(G_i)] = 0$ for all $i$ and $\sum_{i=1}^n \bE[t^2(G_i)] = 1$, a Studentized nonlinear statistic is a statistic that can be written into the generic form 
  \begin{equation} \label{sn_stat_generic_form}
 \frac{ \sum_{i =1}^n t_i(G_i)+ \mathcal{D}_1}{\sqrt{1 +\mathcal{D}_2}},
 \end{equation}
 where $\mathcal{D}_1 = \mathcal{D}_1( G_1, \dots, G_n)$ and $\mathcal{D}_2 =  \mathcal{D}_2( G_1, \dots, G_n)$ are (usually small)  remainder terms that may depend on the data $G_1, \dots, G_n$, with $\mathcal{D}_2 > -1$ almost surely; in particular, the generic statistic in  \eqref{sn_stat_generic_form} is ``Studentized'' by the data-driven variance estimate $1 +\mathcal{D}_2$. Indeed, 
 \begin{multline} \label{obs_as_sn_stat}
 \text{$T_{SN}$ can be understood as an instance of \eqref{sn_stat_generic_form} }\\
 \text{ when conditioning on a given data realization $\cX \in \cE_{1, \cX}$}:
 \end{multline}
From  the conditional moment calculations in \eqref{cX_cond_moment_properties_of_xi}, we have
  \begin{equation} \label{SN_condition_satisfied_1}
 \bE[\xi_\boldi \mid \cX ] = 0 \text{ and }
 \sum_{\boldi \in \Inm} \bE[\xi_\boldi^2 \mid \cX ] = \frac{ U_{h^2}}{\bar{U}_{h^2} }  = 1\text{ for }  \cX \in \cE_{1, \cX},
 \end{equation}
 where the last equality comes from the fact in \eqref{Uh2_equal}.
 Moreover, only $D_1$ depends on $\cZ$, and $D_2$ is no smaller than $-1$ because, by \eqref{Uh2_equal} again, 
  \begin{equation} \label{1plusD2_equal_cEx}
  1 + D_2  
 =1 - \frac{\sigma_h^2}{\bar{U}_{h^2}}\Big( \frac{{U}_{h^2}}{\sigma_h^2} -1 \Big) = \frac{\sigma_h^2}{\bar{U}_{h^2}} >0 \text{ for }  \cX \in \cE_{1, \cX}.
  \end{equation}
 These verify the claim in \eqref{obs_as_sn_stat}.

Note that all of $D_1$, $\{D_1^{(\boldi)}\}_{\boldi \in \Inm}$ and $\frakz_\cX$ are arbitrary and not   yet defined at this point. Turns out, we can finish our  proofs of \thmsref{BE_bdd_Nlln} and \thmssref{main} by  picking appropriate configurations for them; recall the definition of $\hat{N}^{(\boldi)} $ and $\hat{N}^{(\boldi)}_c $ in \eqref{Nhat_leave_one_out} and \eqref{Nhatc_leave_one_out} in  \secref{pf_preps}.

 \begin{configuration}[Configuration of $D_1, D_1^{(\boldi)}, \frakz_\cX$ for  ``$N \ll n^d$''] \label{config:Nlln}
 For proving  \thmref{BE_bdd_Nlln}, we will set
 \begin{enumerate}[(i)]
   \item $\frakz_\cX = z$;
 \item $ D_1 = \frac{N - \hat{N}_c}{\hat{N}_c}   \sum_{{\bf i} \in I_{n, m}}  \xi_\boldi  +  \frac{ N^{3/2}U_n}{\hat{N}_c \sqrt{ \bar{U}_{h^2} (1 - p)}  } $;
 \item    $D_1^{(\boldi)} =  \frac{N - \hat{N}_c^{(\boldi)}}{\hat{N}_c^{(\boldi)}} \sum_{{\bf j} \in I_{n, m}, \boldj \neq \boldi}\xi_\boldj  +  \frac{ N^{3/2}U_n}{\hat{N}_c^{(\boldi)}\sqrt{ \bar{U}_{h^2} (1 - p)}  }$ \text{ for each } $\boldi \in \Inm$. 
 \end{enumerate}

 \end{configuration}
 
  \begin{configuration}[Construction of $D_1, D_1^{(\boldi)}, \frakz_\cX$ for  ``$N \asymp n$'' ]  \label{config:Nasympn}

  For proving  \thmref{main}, we will set
  \begin{enumerate}[(i)]
   \item  $\frakz_\cX \equiv \frac{z \sigma}{\sigma_h\sqrt{ \alpha(1 -p)}} -   \frac{\sqrt{N}U_n}{\sigma_h \sqrt{1 - p}}$;
   \item $D_1  \equiv  \frac{N - \hat{N}_c}{\hat{N}_c} \Big( \sum_{\boldi \in \Inm} \xi_\boldi  +   \frac{ \sqrt{N}U_n}{ \sqrt{\bar{U}_{h^2} (1 - p)}}\Big)$;
   \item $D_1^{(\boldi)} 
\equiv   \frac{N - \hat{N}_c^{(\boldi)}}{\hat{N}_c^{(\boldi)}} \Big(\sum_{{\bf j} \in I_{n, m}, \boldj \neq \boldi}\xi_\boldj +   \frac{ \sqrt{N}U_n}{ \sqrt{\bar{U}_{h^2} (1 - p)}}\Big)$ \text{ for each } $\boldi \in \Inm$.
 \end{enumerate}
 
 \end{configuration}

Under these two configurations, we have  the following  bounds on the terms $T(D_1)$ and $T(D_2)$ defined in \lemref{prelim_BE}:

\begin{lemma}[Bounds on $T(D_1)$] \label{lem:TD1_bdd} 
Suppose that \eqref{mean0_assumption}, $\bE[|h|^3] < \infty$  and  $2 \leq m < n/2$ hold, and  $T(D_1)$ is defined as in \lemref{prelim_BE}.
\begin{enumerate}[(i)]
\item   If $D_1$ and $D_1^{(\boldi)}$ for each ${\boldi \in \Inm}$ are defined as in \configref{Nlln}, then
 \[
 T(D_1) \leq C\bigg( \frac{1}{\sqrt{N}} + \frac{  \sqrt{N K_{n, m, d}}}{\sqrt{1-p}} \bigg)
 \]
for an absolute constant $C >0$, where $K_{n, m, d}$ is  defined as in \thmref{BE_bdd_Nlln}.
\item If $D_1$ and $D_1^{(\boldi)}$ for each ${\boldi \in \Inm}$  are defined as in \configref{Nasympn}, then  \[
 T(D_1) \leq C\bigg( \frac{1}{\sqrt{N}} + \frac{  \sqrt{m}}{\sqrt{n(1-p)}} \bigg)
 \]
 for an absolute constant $C >0$.  
 \end{enumerate}
\end{lemma}

\begin{lemma}[Bounds on $T(D_2)$]  \label{lem:TD2_bdd} 
Suppose that  \eqref{mean0_assumption}, $\bE[|h|^3] < \infty$  and  $2 \leq m < n/2$ hold, and  $T(D_2)$ is defined as in \lemref{prelim_BE}.
\begin{enumerate}[(i)]
\item If $\frakz_\cX$ is defined as in \configref{Nlln}$(i)$, then
 \begin{multline*}
T(D_2) \leq 
C \bigg[ \frac{\bE[|h|^3] (1 - 2p + 2p^2)}{\sigma_h^3 (N(1-p))^{1/2}} +   \frac{m^{3/2}\bE[\Psi_1^{3/2}]}{ n^{1/2}\sigma_h^3}   +
   \|\cR\|_{3/2}  \bigg]+ P(\Omega \backslash \cE_{1, \cX} ) .
 \end{multline*}
 for an absolute constant $C >0$.
\item  If $\frakz_\cX$ is defined as in \configref{Nasympn}$(i)$, then
 \begin{multline*}
T(D_2) \leq 
C \bigg[ \frac{\bE[|h|^3] (1 - 2p + 2p^2)}{\sigma_h^3 (N(1-p))^{1/2}} +  \bigg(1 + \frac{\sqrt{Nm}}{\sqrt{n(1-p)}}\bigg)  \frac{m^{3/2}\bE[\Psi_1^{3/2}]}{ n^{1/2}\sigma_h^3}   +
  \|\cR\|_{3/2} \bigg]+ P(\Omega \backslash \cE_{1, \cX} ) 
 \end{multline*}
  for an absolute constant $C >0$.
 \item If, in addition, $ \bE[h^4] < \infty$, then
\[
 T(D_2) \leq \frac{8 \sqrt{m \bE[(h^2 - \sigma_h^2)^2]}}{\sigma_h^2 \sqrt{n}},
 \]
  regardless of whether $\frakz_\cX$ takes the form in \configref{Nlln}$(i)$ or \configssref{Nasympn}$(i)$.
\end{enumerate}
\end{lemma}
These two lemmas are proved in  \appsref{bounding_T_D1} and  \appssref{bounding_T_D2}  respectively; the proof of \lemref{TD1_bdd}  in \appref{bounding_T_D1} entails  tight moment calculations that  exploits conditioning  on the raw data $\cX$  wherever necessary, and the proof of \lemref{TD2_bdd} in \appref{bounding_T_D2} entails   leave-one-out arguments very unique to this problem. We advise the reader to skip these proofs on first reading.

In  the next two subsections, with basic arguments, we will finish the proofs of \thmsref{BE_bdd_Nlln} and \thmssref{main} using  the ingredients  above, as well as the simple facts that 
\begin{equation} \label{simple_fact_1}
\frac{1}{\sqrt{N}} \leq  \frac{14 \sqrt{7}}{25} \cdot \frac{\bE[|h|^3] (1 - 2p + 2p^2)}{\sigma_h^3 (N(1-p))^{1/2}}
\end{equation}
(by virtue of $ \min_{p \in [0,1]}\frac{1 - 2p + 2p^2}{(1 - p)^{1/2}} =  \frac{25}{14 \sqrt{7}}$ and  $\bE[|h|^3] \geq \sigma_h^3$) and 
\begin{equation} \label{simple_fact_2}
|\Phi(az) - \Phi(z)| \leq \frac{ \exp(- 0.5)}{\sqrt{ 2\pi}} |a - 1| \text{ for any } a \geq 1 \text{ and } z  \in \bR
\end{equation}
(by virtue of  Taylor's theorem and $\sup_{z \in \bR} |z \phi(z)| = \frac{\exp(- 0.5)}{\sqrt{2 \pi}} $).

 \subsection{Proof of \thmref{BE_bdd_Nlln} }


Throughout this subsection, we let  $\frakz_\cX$, $D_1$ and $D_1^{\boldi}$ for each $\boldi \in \Inm$ be as  in \configref{Nlln}. From the basic decomposition in  \eqref{basic_decomp} and  \eqref{censored_N_equal_N}, write
 \begin{multline} \label{pf_Nlln_N_censored_1st_step}
\bigg|P\bigg(\frac{\sqrt{N} U_{n, N}'}{\sigma_h}\leq z \bigg) - \Phi(z) \bigg|\\
 \leq  \bigg|P\bigg( \frac{N^{3/2} B_n }{\hat{N}_c\sigma_h   } +  \frac{ N^{3/2}U_n}{\hat{N}_c\sigma_h\sqrt{(1 - p)}  } \leq \frac{z}{\sqrt{(1 - p)} } \bigg) - \Phi(z) \bigg| + P(\Omega \backslash \cE_\cZ).
 \end{multline} 
Because $1 + D_2 = \sigma_h^2/\bar{U}_{h^2}$ on $\cE_{1, \cX }$ (cf. \eqref{1plusD2_equal_cEx}), one can simply check from the definition of $T_{SN}$  that, when $D_1$ is defined as in \configref{Nlln}$(ii)$, 
\[
 T_{SN} = \frac{N^{3/2} B_n }{\hat{N}_c\sigma_h   } +  \frac{ N^{3/2}U_n}{\hat{N}_c\sigma_h\sqrt{(1 - p)}  }  \text{ on }  \cE_{1, \cX}.
\]
Hence,  continuing from \eqref{pf_Nlln_N_censored_1st_step}, we can further write
\begin{align}
&\bigg|P\bigg(\frac{\sqrt{N} U_{n, N}'}{\sigma_h}\leq z \bigg) - \Phi(z) \bigg| \notag\\
&\leq  \bigg|P\Big( T_{SN}\leq \frac{z}{\sqrt{(1-p)}} \Big) - \Phi(z) \bigg| + P(\Omega \backslash \cE_{1, \cX}) + P(\Omega \backslash \cE_\cZ) \notag\\
&\leq \big|P( T_{SN}\leq z) - \Phi(z) \big| +\bigg| \Phi\Big(\frac{z}{\sqrt{1-p}}\Big) - \Phi(z) \bigg| + P(\Omega \backslash \cE_{1, \cX}) + P(\Omega \backslash \cE_\cZ) \notag \\
&\leq  \big|P( T_{SN}\leq z) - \Phi(z) \big| +\frac{ \exp(- 0.5)}{\sqrt{ 2\pi}} \Big|\frac{\sqrt{p}}{\sqrt{1-p}} \Big| + P(\Omega \backslash \cE_{1, \cX}) + P(\Omega \backslash \cE_\cZ), \notag
\end{align}
where the last inequality uses  \eqref{simple_fact_2}.
Applying \eqref{bernstein1}, \eqref{exp_lower_bdd_app1}, \eqref{simple_fact_1},  \lemref{prelim_BE}, \lemref{TD1_bdd}$(i)$, \lemref{TD2_bdd}$(i)$\&$(iii)$ and   the fact that
\[
\frac{\sqrt{p}}{\sqrt{1-p}} =  \frac{\sqrt{N } {n \choose m}^{-1/2}}{\sqrt{1 -p}} \leq  \frac{  \sqrt{N K_{n, m, d}}}{\sqrt{1-p}}
\]
(because ${n \choose m}^{-1} \leq K_{n , m , d}$) to the last inequality, we have proven  \thmref{BE_bdd_Nlln}.

\subsection{Proof of \thmref{main}} \label{sec:finishing_main_pf}
Throughout this subsection, we let  $\frakz_\cX$, $D_1$ and $D_1^{\boldi}$ for each $\boldi \in \Inm$  be as in \configref{Nasympn}.
The proof of \thmref{main} can be separated into three broad steps:

\subsubsection{Bounding $|P(\sqrt{n}U_{n, N}' /\sigma \leq z ) - \Phi(z) |$ in terms of $ | P (  T_{SN} \leq  \frakz_\cX ) -  \bE[\Phi(\frakz_\cX)] |$}

For a given $z \in \bR$ and $\frakz_\cX$ defined in \configref{Nasympn}$(i)$ (which \emph{does} depend on the data $\cX$),
we will first follow iterative arguments  similar to those pioneered by \citet{chen2019randomized} and establish the bound
\begin{multline} \label{BE_bdd_after_iterative_argument}
\bigg|P\bigg(\frac{\sqrt{n}U_{n, N}' }{\sigma} \leq z \bigg) - \Phi(z) \bigg| \leq   \sup_{z' \in \bR}  \bigg| P\bigg(\frac{\sqrt{n} U_n }{m\sigma_g} \leq  z' \bigg) - \Phi (z')  \bigg| +    \\
 \varepsilon_z +  P(\Omega \backslash \cE_\cZ) + P(\Omega \backslash \cE_{1, \cX}) +  \big| P( T_{SN} \leq  \frakz_\cX ) -  \bE[\Phi(\frakz_\cX)] \big|  ,
\end{multline}
where
\begin{equation}\label{var_epsilon_def}
\varepsilon_z \equiv \Big| \Phi \Big(\frac{\sigma z}{\tilde{\sigma}}  \Big)  - \Phi(z)\Big| 
\quad \text{ for } \quad  \tilde{\sigma}^2 \equiv  m^2  \sigma_g^2+  \alpha (1 - p) \sigma_h^2.
\end{equation}
By defining 
\[
T \equiv \frac{N^{3/2}  B_n }{\hat{N}_c\sigma_h}  +  \bigg(  \frac{N - \hat{N}_c}{\hat{N}_c}\bigg) \frac{ \sqrt{N}U_n}{\sigma_h \sqrt{1 - p}},
\]
one can easily check that, when $D_1$ is defined as in \configref{Nasympn}$(ii)$, 
\[
 T_{SN} = T\text{ on }  \cE_{1, \cX};
\]
therefore, to prove \eqref{BE_bdd_after_iterative_argument}, it suffices to establish
\begin{multline} \label{BE_bdd_immediately_after_iterative_argument}
\bigg|P\bigg(\frac{\sqrt{n}U_{n, N}' }{\sigma} \leq z \bigg) - \Phi(z) \bigg| \leq     \sup_{z' \in \bR}  \bigg| P\bigg(\frac{\sqrt{n} U_n }{m\sigma_g} \leq  z' \bigg) - \Phi (z')  \bigg| + \\
 \varepsilon_z+ P(\Omega \backslash \cE_\cZ) +  | P( T \leq  \frakz_\cX ) -  \bE[\Phi(\frakz_\cX)] |.
\end{multline}

To prove \eqref{BE_bdd_immediately_after_iterative_argument},
let $Y_1$ and $Y_2$ be two independent standard normal random variables that are also independent of $\cX$ and $\cZ$. 
With the equality in \eqref{U_prime_as_U_star}, one can first write the inequality 
\begin{align} \label{1st_iter_conditioning}
P\bigg(\frac{\sqrt{n}U_{n, N}' }{\sigma} \leq z \bigg)  &\leq 
 P\bigg(\frac{N\sqrt{n}  B_n }{\hat{N}_c\sigma}  +   \frac{(N - \hat{N}_c)\sqrt{n}U_n}{\hat{N}_c\sigma \sqrt{1 - p}}  \leq \frac{z}{\sqrt{1 -p}} -   \frac{\sqrt{n}U_n}{\sigma \sqrt{1 - p}}\bigg)  + P(\Omega \backslash \cE_\cZ) \notag\\
 &= P(T \leq \frakz_\cX) + P(\Omega \backslash \cE_\cZ).
\end{align}
The first term on the right hand side of  \eqref{1st_iter_conditioning} can be further bounded as
\begin{equation}
P(T   \leq  \frakz_\cX )  \leq P\left( Y_1 \leq   \frakz_\cX  \right) +  \big| P\big( T \leq  \frakz_\cX\big) -  \bE[\Phi(\frakz_\cX)]\big|.
 \label{2nd_iter_conditioning}
 \end{equation}
 Continuing, we can bound the first term on the right hand side of \eqref{2nd_iter_conditioning} as 
  \begin{align}
P\left( Y_1 \leq   \frakz_\cX  \right)&= P\bigg( \frac{\sqrt{N}U_n}{\sigma_h \sqrt{1 - p}} \leq \frac{z \sigma}{\sigma_h\sqrt{\alpha(1 -p)}} - Y_1 \bigg)  \notag\\
&= \bE \bigg[ P\bigg( \frac{\sqrt{n} U_n }{m \sigma_g} 
\leq \frac{ z \sigma}{  m \sigma_g} - \frac{ \sigma_h\sqrt{ \alpha (1- p) } }{m \sigma_g} Y_1 \mid Y_1 \bigg) \bigg]
\notag\\
&\leq  P\bigg(Y_2 \leq  \frac{z \sigma}{ m \sigma_g} - \frac{ \sigma_h\sqrt{ \alpha (1- p) } }{m \sigma_g} Y_1\bigg) + \sup_{z' \in \bR}  \Big| P\Big(\frac{\sqrt{n} U_n }{m\sigma_g} \leq  z' \Big) - \Phi (z')  \Big| \notag\\
&= \Phi \Big(\frac{\sigma z}{\tilde{\sigma}}  \Big) 
+ \sup_{z' \in \bR}  \bigg| P\bigg(\frac{\sqrt{n} U_n }{m\sigma_g} \leq  z' \bigg) - \Phi (z')  \bigg| \notag\\
&\leq \Phi(z)  + \sup_{z' \in \bR}  \bigg| P\bigg(\frac{\sqrt{n} U_n }{m\sigma_g} \leq  z' \bigg) - \Phi (z')  \bigg| + \varepsilon_z \label{3rd_iter_conditioning},
\end{align}
where $\tilde{\sigma}$ and $\varepsilon_z$ are as defined in \eqref{var_epsilon_def};
the  equality above \eqref{3rd_iter_conditioning} uses that  $m \sigma_g Y_2 + \sigma_h \sqrt{\alpha (1- p)} Y_1$ is distributed as a normal random variable with mean $0$ and variance $\tilde{\sigma}^2$.
Combining \eqref{1st_iter_conditioning}-\eqref{3rd_iter_conditioning}, we have
\begin{multline} \label{1st_to_3rd_iter_conditioning}
P\Big(\frac{\sqrt{n}U_{n, N}' }{\sigma} \leq z \Big) -  \Phi(z) \leq  \sup_{z' \in \bR}  \Big| P\Big(\frac{\sqrt{n} U_n }{m\sigma_g} \leq  z' \Big) - \Phi (z')  \Big| + \\
 \varepsilon_z + 
 P(\Omega \backslash \cE_\cZ)+ | P ( T \leq  \frakz_\cX ) - \bE[ \Phi(\frakz_\cX)]| 
\end{multline}
 By reverse arguments, in parallel to \eqref{1st_iter_conditioning}-\eqref{3rd_iter_conditioning}, we can get
 \begin{equation*} 
P\Big(\frac{\sqrt{n}U_{n, N}' }{\sigma} \leq z \Big)  \geq 
P(T \leq \frakz_\cX)
 - P(\Omega \backslash \cE_\cZ),
\end{equation*}
\begin{equation*}
P(T \leq \frakz_\cX)  
 \geq P\left( Y_1 \leq   \frakz_\cX  \right) -  | P( T  \leq  \frakz_\cX) -  \bE[\Phi(\frakz_\cX)]|
 \end{equation*}
 and
 \begin{equation*}   
 P\left( Y_1 \leq   \frakz_\cX  \right) \geq   \Phi(z)  - \sup_{z' \in \bR}  \bigg| P\bigg(\frac{\sqrt{n} U_n }{m\sigma_g} \leq  z' \bigg) - \Phi (z')  \bigg| - \varepsilon_z,
 \end{equation*}
 which can be combined to give 
 \begin{multline}  \label{4th_to_6th_iter_conditioning}
P\bigg(\frac{\sqrt{n}U_{n, N}' }{\sigma} \leq z \bigg) - \Phi(z) \geq  - \sup_{z' \in \bR}  \bigg| P\bigg(\frac{\sqrt{n} U_n }{m\sigma_g} \leq  z' \bigg) - \Phi (z')  \bigg| - \\
 \varepsilon_z -
 P(\Omega \backslash \cE_\cZ) - \big| P\big( T  \leq  \frakz_\cX\big) -  \bE[\Phi(\frakz_\cX)]\big|.
\end{multline}
Together, \eqref{1st_to_3rd_iter_conditioning} and \eqref{4th_to_6th_iter_conditioning} render \eqref{BE_bdd_immediately_after_iterative_argument}.

\subsubsection{Bounding $\varepsilon_z$ in \eqref{var_epsilon_def}} We will show that, for an absolute constant $C > 0$,
\begin{equation} \label{var_ep_z_bdd}
\varepsilon_z \leq \frac{0.31 \cdot  N}{n^2 (1-p)} \text{ under the assumptions of \thmref{main}}
\end{equation}
irrespective of the value of $z$.  Since  $0< \tilde{\sigma} \leq \sigma$,  by \eqref{simple_fact_2}, 
we get  that
\begin{align} 
\varepsilon_z
&\leq  \frac{ \exp(- 0.5)}{\sqrt{ 2\pi}} \Bigg|\frac{\sigma  - \tilde{\sigma} }{\tilde{\sigma}}   \Bigg|  
\notag \\
&=     \frac{\exp(- 0.5)}{\sqrt{ 2\pi}}  \Bigg(\frac{  n {n \choose m}^{-1} \sigma_h^2 }{\tilde{\sigma} ( \sigma+ \tilde{\sigma})} \Bigg) \text{ by } \frac{\sigma  - \tilde{\sigma} }{\tilde{\sigma}} = \frac{\sigma^2  - \tilde{\sigma}^2 }{\tilde{\sigma} (\sigma + \tilde{\sigma})} \notag \\
&\leq \frac{\exp(- 0.5)}{2 \sqrt{ 2\pi}}  \Bigg(\frac{  n {n \choose m}^{-1} \sigma_h^2 }{\tilde{\sigma}^2} \Bigg)    \text{ by } \tilde{\sigma} \leq \sigma \notag\\
&\leq  \frac{\exp(- 0.5) N}{\sqrt{ 2\pi}  (n-1)n (1-p)}  \label{varep2_ultimate_bdd} \notag,
\end{align}
where the last inequality uses $\alpha (1 - p) \sigma_h^2 \leq \tilde{\sigma}^2 $ and 
\begin{equation*}
\frac{n}{{n \choose m} } = \frac{m \cdots 1}{ (n-1)\cdots (n-(m -1)) }  
\leq \frac{2}{n-1} \text{ as a consequence of } 2 \leq m < n/2.
\end{equation*}
The last inequality gives \eqref{var_ep_z_bdd}, as $\frac{n}{n-1} \leq \frac{5}{4}$ under $2 \leq m < \frac{n}{2}$ and $\frac{5\exp(- 0.5)}{4\sqrt{ 2\pi}} < 0.31$.
\subsubsection{Wrapping up} With  \eqref{BE_bdd_after_iterative_argument}  and \eqref{var_ep_z_bdd}, along with the bounds \eqref{complete_U_stat_BE},  \lemref{prelim_BE} and  \eqref{exp_lower_bdd_app1} for $ \sup_{z' \in \bR}  \big| P\big(\frac{\sqrt{n} U_n }{m\sigma_g} \leq  z' \big) - \Phi (z')  \big|$, $| P( T_{SN} \leq  \frakz_\cX ) -  \bE[\Phi(\frakz_\cX)]|$ and $P(\Omega \backslash \cE_{1, \cX})$, 
 we can first write
\begin{equation*}  
\bigg|P\bigg(\frac{\sqrt{n}U_{n, N}' }{\sigma} \leq z \bigg) - \Phi(z) \bigg| \leq   \frac{0.31 \cdot  N}{n^2 (1-p)}  +      
  C\bigg(\frakB_1+\frakB_2+ T(D_1) + T(D_2)   \bigg),
\end{equation*}
where  $P(\Omega \backslash \cE_\cZ)$ has been absorbed into $\frac{\bE[|h|^3] (1 - 2p + 2p^2)}{\sigma_h^3 (N(1-p))^{1/2}}$ of  $\frakB_1$ via the exponential inequality   in \eqref{bernstein1}, $\exp ( - 3N/28) \leq C/\sqrt{N}$ and  \eqref{simple_fact_1}. Applying the bound for $T(D_1)$ in \lemref{TD1_bdd}$(ii)$ to the  above display, we further get  
\begin{equation}  
\label{right_before_final_main_thm_pf}
\bigg|P\bigg(\frac{\sqrt{n}U_{n, N}' }{\sigma} \leq z \bigg) - \Phi(z) \bigg| \leq  \frac{0.31 \cdot  N}{n^2 (1-p)} +    
    C\bigg( \frakB_1+ \frakB_2+ \frac{  \sqrt{m}}{\sqrt{n(1-p)}}  + T(D_2)   \bigg),
\end{equation}
where $1/\sqrt{N}$ is absorbed into $\frac{\bE[|h|^3] (1 - 2p + 2p^2)}{\sigma_h^3 (N(1-p))^{1/2}}$ of  $\frakB_1$ by  \eqref{simple_fact_1}.

Finally, to get the B-E bound in \thmref{main} under the stronger finite fourth moment assumption $\bE[h^4] < \infty$, we apply \lemref{TD2_bdd}$(iii)$ to bound $T(D_2)$ in \eqref{right_before_final_main_thm_pf}. To get the B-E bound in \thmref{main} under the weaker finite third moment assumption $\bE[|h|^3] < \infty$ from    \eqref{right_before_final_main_thm_pf}, 
we apply the bound for $T(D_2)$ in \lemref{TD2_bdd}$(ii)$, the bound for $P(\Omega \backslash \cE_{1, \cX})$ in \eqref{exp_lower_bdd_app1} and the fact that
\begin{equation} \label{last_absorption}
\frac{N}{n^2 (1-p)} \leq \frac{\sqrt{N}}{n \sqrt{(1-p)}} \leq \frac{\sqrt{N}}{ \sqrt{n(1-p)}} \frac{\bE[\Psi^{3/2}]}{\sqrt{n} \sigma_h^3},
\end{equation}
to absorb $\frac{N}{n^2 (1-p)} $ into  $\big(1 + \frac{\sqrt{Nm}}{\sqrt{n(1-p)}}\big)  \frac{m^{3/2}\bE[\Psi_1^{3/2}]}{ n^{1/2}\sigma_h^3} $; the terms $\frac{\bE[|h|^3] (1 - 2p + 2p^2)}{\sigma_h^3 (N(1-p))^{1/2}}$ and $\exp \Big ( \frac{- n \sigma_h^6}{24 m(\bE[|h|^3])^2} \Big)$ are absorbed into $\frakB_1$.

\section{Discussion} \label{sec:discuss}

With more careful book keeping, it is  possible to produce numerical estimates of the absolute constants $C$ appearing in \thmsref{BE_bdd_Nggn}-\thmssref{main}. The most important contribution of this work nevertheless lies in our methodology of adapting the latest techniques in the Stein's method literature to establish B-E bounds  of  natural rates in $N$ and $n$ for $U_{n, N}'$, which we believe to be optimal. That being said, we don't exclude the possibility that comparable bounds with simpler forms could be  proved, potentially with alternative methods. At present, it is unclear to us how the terms in our  bounds can be further simplified, such as by  the ``absorption trick'' illustrated in \eqref{last_absorption}. We also speculate that some  techniques might  be combined with the recent advances in the Berry-Esseen-type theory for sums of independent random vectors  \citep{cckk2022, ccky_2023} to improve upon the existing Berry-Esseen-type theory for  vector-valued incomplete U-statistics \citep{SongChenKato2019, chen2019randomized}.

 It  takes only slight adjustments to our present methodology to establish B-E bounds for the random kernel incomplete U-statistic $U_{n, N, \omega}'$  in \eqref{random_kernel_icu} as well. For instance, the Hoeffding decomposition is also available to random kernel U-statistics \citep[Definition 2]{peng2022rates}, and the conditioning on $\cX \in \cE_{1, \cX}$ in \eqref{obs_as_sn_stat} simply has to be changed to conditioning on both $\cX \in \cE_{1, \cX}$ and $\omega$ altogether. To keep our presentation the most streamlined possible, we leave these details to the reader. Moreover, in the real applications for uncertainty quantification of ensemble predictions mentioned in \secref{random_kernel_ICU}, one almost always has to \emph{studentize} their incomplete U-statistics by data-driven estimates of the standardizing factors such as $\sigma_g$ and $\sigma_h$ in \thmref{weak_convergences}; we refer the reader to \citet{zhou_mentch_hooker_2021} and \citet{peng2021bias} for related work in this vein. Therefore, the more imperative  and challenging task is to develop a Berry-Esseen theory for \emph{studentized} incomplete U-statistics that is directly applicable; we are currently pursuing a deep investigation into this.

\newpage
\appendix 


\section{Two technical lemmas} \label{app:two_technical_lemmas}

We need  the a Bennett's inequality from  \citet[Appendix A]{leung2024another}, as well as some bounds on the solution to the Stein equation in \eqref{Stein_solution}.
\begin{lemma} [Bennett's inequality for a sum of censored random variables]  \label{lem:Bennett} 
Let $\xi_1, \dots, \xi_n$ be independent random variables with $\bE[\xi_i] = 0$ for all $i = 1, \dots, n$ and $\sum_{i=1}^n \bE[\xi_i^2] \leq 1$, and define $\barxi_i = \xi_i I(|\xi_i| \leq 1) + 1 I(\xi_i >1) - 1  I(\xi_i <-1)$. 
Then  
\[
\bE[e^{\sum_{i=1}^n \barxi_i}] \leq  \exp\left( 4^{-1} (e^2+1) \right) \leq 8.15.
\]
\end{lemma}

 \begin{lemma}[Bounds for $f_z$ and $\frac{\partial f_z}{\partial w}$] \label{lem:fxfxprime_bdd}

The following bounds are true for $f_z$  in \eqref{Stein_solution} and its derivative $\frac{\partial f_z}{\partial w}$:
\begin{enumerate}[(i)]
\item For all $w, z \in \bR$, $0 < f_z(w) \leq 0.63$ and $ |\frac{\partial f_z}{\partial w}(w)| \leq 1$ .

\item For $z\geq 1$, 

  \begin{equation*} \label{fxbdd}
 f_z (w) \leq  
  \begin{cases} 
  1.7 e^{-z} &  \text{ for }\quad w \leq z - 1\\
  1/z     &    \text{ for }\quad z  - 1 < w \leq z \\
   1/w     &   \text{ for }\quad z  < w 
     \end{cases}
     \end{equation*}
     and
    \begin{equation*} \label{fx'bdd}
      \Big| \frac{\partial f_z}{\partial w}(w) \Big| \leq  
  \begin{cases} 
    e^{1/2 -z } &   \text{ for }\quad w \leq z - 1\\
   1       &    \text{ for }\quad z  - 1 < w \leq z \\
      (1 + z^2)^{-1}    &  \text{ for }\quad w > z  
     \end{cases}.
\end{equation*}
\item For $z \leq -1$, 
\begin{equation*} \label{fxbddzleqminus1}
 f_z (w) \leq  
  \begin{cases} 
  1.7 e^{-|z|} &  \text{ if }\quad w \geq z+1\\
  1/|z|     &    \text{ if }\quad  z \leq w < z+1 \\
   1/|w|     &   \text{ if }\quad w  < z 
     \end{cases} 
     \end{equation*}
     and
    \begin{equation*} \label{fx'bddzleqminus1}
      \Big| \frac{\partial f_z}{\partial w}(w) \Big| \leq  
  \begin{cases} 
    e^{1/2 -|z| } &   \text{ if }\quad w \geq  z+1\\
   1       &    \text{ if }\quad z \leq w < z +1\\
      (1 + |z|^2)^{-1}    &  \text{ if }\quad w  < z 
     \end{cases}.
     \end{equation*}
     \item $\sup_{z , w \in \bR}|zf_z(w)| \leq 1$.
\end{enumerate}
\end{lemma}

\begin{proof}[Proof of \lemref{fxfxprime_bdd}]
\citet[Appendix A]{leung2024another} has already established 
 $(i)$ and  $(ii)$.
 To prove $(iii)$, let $\tilde{z} = - z$ and $\tilde{w} = -w$. By the definition in \eqref{Stein_solution}, it is easy to see that
 \begin{equation*} \label{flip}
 f_z(w) = f_{\tilde{z}} (\tilde{w});
 \end{equation*}
 since $\tilde{z} \geq 1$, we can apply $(ii)$ to $f_{\tilde{z}} (\tilde{w})$ and get $(iii)$. $(iv)$ is a simple consequence of the bounds for $f_z$ in $(i)$-$(iii)$.

\end{proof}

 \section{Proof of \lemref{Delta_estimates_Nggn}, estimates regarding the remainder terms $\Delta_1$ and $\Delta_2$} \label{app:pf_delta_estimate_Nggn}

 For the first term, we use \eqref{Ubarhm_sq_estimate} and \eqref{exp_Bn2} to get
 \begin{align}
  \bE\bigg[ \bigg| \frac{\sum_{i=1}^n g(X_i)}{\sigma_g\sqrt{n}}  \cdot \Delta_1 \bigg|\bigg] 
  &\leq \Big\| \frac{\sum_{i=1}^n g(X_i)}{\sqrt{n} \sigma_g}\Big \|_2 \bigg(\Big\|\frac{\sqrt{n}U_{\bar{h}_m}}{ m \sigma_g}\Big \|_2 + \Big\|\frac{\sqrt{n(1 - p)}B_n}{ m \sigma_g}  \Big \|_2\bigg)  \notag\\
&\leq   \sqrt{\frac{2m}{n} \Big( \frac{\sigma_h^2}{m \sigma_g^2} - 1\Big) } +  \sqrt{\frac{n}{N}(1 - p)\frac{\sigma_h^2}{m \sigma_g^2}} \label{1st_bdd_chen_shao_lemma}
 \end{align}
 For the second term, we use \eqref{Ubarhm_leave_one_out_sq_estimate} and \eqref{Bn_minus_Bni_2nd_moment} to get 
 \begin{align}
  \bE\bigg[ \bigg| \frac{g(X_i)}{\sigma_g\sqrt{n}} \cdot( \Delta_1- \Delta_1^{(i)}) \bigg| \bigg] 
  &\leq 
  \bigg\| \frac{ g(X_i)}{\sqrt{n} \sigma_g}\bigg \|_2 \bigg\|  \frac{\sqrt{n}(U_{\bar{h}_m}- {U_{\bar{h}_m}^{(i)}} + \sqrt{(1 - p)} (B_n - B_n^{(i)}))}{m \sigma_g}\bigg\|_2 \notag\\
  &\leq \frac{1}{\sqrt{n}} \bigg(  \frac{ 2\sqrt{m}}{n} \bigg(\frac{\sigma_h^2}{m \sigma_g^2} -1\bigg)^{1/2} + 
  \frac{\sigma_h}{\sigma_g}  \sqrt{\frac{(1-p)  }{Nm}}\bigg)  \label{2nd_bdd_chen_shao_lemma}
 \end{align}
For the third term, we use
\begin{align}
 &\bE\Bigg[ \bigg| \frac{\sum_{i=1}^n g(X_i)}{\sigma_g\sqrt{n}}  \cdot \Delta_2 \bigg|\Bigg] \notag\\
 &\leq \bigg\| \bigg( \frac{N - \hat{N}_c}{\hat{N}_c} \bigg) \frac{\sum_{i=1}^n g(X_i)}{\sigma_g\sqrt{n}}  \bigg\|_2  \bigg\|   \frac{\sqrt{n} (U_n + \sqrt{(1-p)} B_n)}{m\sigma_g }\bigg\|_2 \notag\\
 &= \bigg\| \bigg( \frac{N - \hat{N}_c}{\hat{N}_c} \bigg) \bigg\|_2 \bigg\| \frac{\sum_{i=1}^n g(X_i)}{\sigma_g\sqrt{n}}  \bigg\|_2 \bigg\|   \frac{\sqrt{n} (U_n + \sqrt{(1-p)} B_n)}{m\sigma_g }\bigg\|_2 \text{ by independence of $\cX$ and $\cZ$} \notag\\
 &=\bigg\| \bigg( \frac{N - \hat{N}_c}{\hat{N}_c} \bigg) \bigg\|_2 \bigg\|   \frac{\sqrt{n} (U_n + \sqrt{(1-p)} B_n)}{m\sigma_g }\bigg\|_2 \notag\\
 &\leq \frac{2}{\sqrt{N}}  \bigg( \Big\| \frac{\sum_{i=1}^n g(X_i)}{\sqrt{n} \sigma_g}\Big \|_2 + \Big\|\frac{\sqrt{n}U_{\bar{h}_m}}{ m \sigma_g}\Big \|_2 + \Big\|\frac{\sqrt{n(1 - p)}B_n}{ m \sigma_g}  \Big \|_2\bigg) \text{ by \eqref{Nhat_1_norm_bddb}, \eqref{light_H_decomp} and $N/2 \leq \hat{N}_c$} \notag\\
 &\leq \frac{2}{\sqrt{N}} \bigg(1 + \sqrt{\frac{2m}{n} \Big( \frac{\sigma_h^2}{m \sigma_g^2} - 1\Big) } +  \sqrt{\frac{n}{Nm}(1 - p)\frac{\sigma_h^2}{ \sigma_g^2}} \bigg) \text{ by \eqref{Ubarhm_sq_estimate} and \eqref{exp_Bn2}}.  \label{3rd_bdd_chen_shao_lemma}
\end{align}
For the fourth term, we can do
\begin{align}
&\bE\bigg[ \Big| \frac{g(X_i)}{\sigma_g\sqrt{n}} \cdot( \Delta_2- \Delta_2^{(i)}) \Big| \bigg] \notag\\
&\leq \bigg\|   \frac{g(X_i)}{\sigma_g\sqrt{n}}   \cdot  \frac{\hat{N}_c - N}{\hat{N}_c}\bigg\|_2 
\bigg\| \frac{\sqrt{n}}{m\sigma_g} \Big( \frac{m g(X_i)}{n} + U_{\bar{h}_m}- {U_{\bar{h}_m}^{(i)}} + \sqrt{(1 - p)} (B_n - B_n^{(i)})\Big) \bigg\|_2 \notag\\
&= \frac{1}{\sqrt{n}}\cdot \bigg\| \frac{\hat{N}_c - N}{\hat{N}_c} \bigg\|_2 \cdot
\bigg\| \frac{\sqrt{n}}{m\sigma_g} \Big( \frac{m g(X_i)}{n} + U_{\bar{h}_m}- {U_{\bar{h}_m}^{(i)}} + \sqrt{(1 - p)} (B_n - B_n^{(i)})\Big) \bigg\|_2 \notag\\
&\leq  
\frac{2}{\sqrt{nN}} \cdot
\bigg(  \frac{1}{\sqrt{n}}+     \frac{ 2\sqrt{m}}{n} \bigg(\frac{\sigma_h^2}{m \sigma_g^2} -1\bigg)^{1/2}+ \frac{\sigma_h}{\sigma_g}  \sqrt{\frac{(1-p)  }{Nm}}  \bigg), \label{4th_bdd_chen_shao_lemma}
\end{align}
where we have used \eqref{Nhat_1_norm_bddb}, $N/2 \leq \hat{N}_c$,  \eqref{Ubarhm_leave_one_out_sq_estimate} and \eqref{Bn_minus_Bni_2nd_moment} at the end.

\section{Proof of \lemref{prelim_BE}, B-E bound for  $|P( T_{SN}  \leq \frakz_\cX) - \mathbb{E} [\Phi(\frakz_\cX)]   |$} \label{app:pf_prelim_BE}

Recalling $\xi_\boldi$ and  $\barxi_\boldi$  in \eqref{xi_def}  and \eqref{xi_censored}, we first let
\begin{equation}\label{Ws_def}
W \equiv \sum_{\boldi \in \Inm}  \xi_\boldi  \text{ and } 
 \barW \equiv \sum_{\boldi \in \Inm} \barxi_\boldi.
\end{equation}
The algebraic inequalities
\begin{equation} \label{alg_inequalities}
1 + s/2 - s^2/2 \leq (1 + s)^{1/2} \leq 1 + s/2 \text{ for all } s \geq -1, 
\end{equation}
lead to the following two sets of event inclusions depending on the sign of $\frakz_\cX$:
\subsubsection{If $\frakz_\cX \geq 0$:}  \label{sec:case1}In this case we have 
\[
\frakz_\cX \Big(1+\frac{D_2}{2}- \frac{D_2^2}{2} \Big) \leq \frakz_\cX (1+D_2)^{1/2} \leq \frakz_\cX \Big(1 + \frac{D_2}{2}\Big)
\]
  via  \eqref{alg_inequalities}, which in turn implies
\begin{multline*}
\{T_{SN} >\frakz_\cX, \frakz_\cX \geq 0\} \subset \\
\bigg\{ W   + D_1  > \frakz_\cX \Big(1 + \frac{D_2}{2}\Big), \frakz_\cX \geq 0 \bigg\} \cup \bigg\{  \frakz_\cX \Big(1+\frac{D_2}{2} - \frac{D_2^2}{2}\Big) < W   + D_1 \leq \frakz_\cX \Big( 1+ \frac{D_2}{2}\Big), \frakz_\cX \geq 0 \bigg\}
\end{multline*}
and 
\[
\{T_{SN} > \frakz_\cX, \frakz_\cX \geq 0\} \supset \bigg\{ W  + D_1  > \frakz_\cX \Big(1 + \frac{D_2}{2}\Big), \frakz_\cX \geq 0\bigg\}. 
\]
\subsubsection{If $\frakz_\cX < 0$:} \label{sec:case2} In this case we  have 
\[
\frakz_\cX \Big(1 + \frac{D_2}{2}\Big) \leq \frakz_\cX (1+D_2)^{1/2} \leq  \frakz_\cX \Big(1+\frac{D_2}{2} - \frac{D_2^2}{2}\Big)
\]
  via  \eqref{alg_inequalities}, which in turn implies
\begin{multline*}
\{T_{SN} \leq \frakz_\cX, \frakz_\cX < 0\}   \subset \\
\bigg\{ W   + D_1  \leq \frakz_\cX \Big(1 + \frac{D_2}{2}\Big), \frakz_\cX < 0\bigg\} \cup \bigg\{  \frakz_\cX \Big( 1+ \frac{D_2}{2}\Big) < W   + D_1 \leq    \frakz_\cX \Big(1+\frac{D_2}{2} - \frac{D_2^2}{2}\Big), \frakz_\cX < 0  \bigg\}
\end{multline*}
and 
\[
\{T_{SN} \leq \frakz_\cX, \frakz_\cX < 0\} \supset \bigg\{ W  + D_1  \leq \frakz_\cX \Big(1 + \frac{D_2}{2}\Big), \frakz_\cX < 0\bigg\}.
\]

\vspace{1cm}

For the rest of this section, we will use the sign function 
\[
\text{sgn}(z) \equiv  I(z \geq 0)-I(z < 0)
\]
 for any $z \in \mathbb{R}$. The conclusions in \secsref{case1} and \secssref{case2}, along with the equality
\begin{equation*} 
I\big(W  + D_1  \leq \frakz_\cX (1 + D_2/2) \big) - \Phi(\frakz_\cX) = \barPhi(\frakz_\cX) - I\big( W   + D_1 > \frakz_\cX (1 + D_2/2)  \big) , 
\end{equation*}
allow us to write
\begin{multline} \label{Tsn_BE_bound_template}
\big|P\big(T_{SN} \leq \frakz_\cX  \big) - \bE[ \Phi(\frakz_\cX)]\big| \\
= \bigg|\bE\Big[\Big(I\big(T_{SN} \leq \frakz_\cX   \big) -  \Phi(\frakz_\cX) \Big)\cdot I(\frakz_\cX <0)\Big] +\bE\Big[\Big( \barPhi(\frakz_\cX)  - I\big(T_{SN} > \frakz_\cX   \big) \Big)\cdot I(\frakz_\cX \geq 0)\Big] \bigg|\\
\leq \sum_{k=1}^2P(|D_k|>1/2 ) + \frac{ 2^{3/2}\bE[|h|^3](1 - 2p + 2p^2)}{\sigma_h^3 \sqrt{N(1 - p)}}+ 
 \big|P\big( \barW - \bar{\Delta}_{\frakz_\cX} \leq \frakz_\cX  \big) - \bE[\Phi(\frakz_\cX)]\big| +  \\
 P\bigg(\text{sgn}(\frakz_\cX) \cdot \frakz_\cX \cdot \bigg(1+\frac{D_2}{2} - \frac{D_2^2}{2}\bigg)  \leq \text{sgn}(\frakz_\cX) \cdot (W  + D_1) \leq  \text{sgn}(\frakz_\cX) \cdot  \frakz_\cX \cdot \bigg(1 + \frac{D_2}{2}\bigg)  \bigg),
\end{multline}
where we have also defined
\begin{equation} \label{barDeltaz_def}
\bar{\Delta}_{\frakz_\cX} \equiv  \frakz_\cX \bar{D}_2/2 -   \bar{D}_1
\end{equation}
and employed the fact 
\begin{equation*}
P\big(\max_{\boldi \in \Inm}|\xi_\boldi| >1 \big)
 \leq  \sum_{\boldi \in \Inm} \bE[\bE[ |\xi_\boldi|^3 \mid \cX]] \leq  \frac{ 2^{3/2}\bE[|h|^3](1 - 2p + 2p^2)}{\sigma_h^3 \sqrt{N(1 - p)}}
\end{equation*}
 by the conditional third moment calculations in  \eqref{cX_cond_moment_properties_of_xi} and $\sigma_h^2/2 \leq \bar{U}_{h^2}$.
On the other hand, since $f_{\frakz_\cX}$ solves the Stein equation 
\[
I(w \leq \frakz_\cX) - \Phi(\frakz_\cX) = \frac{\partial f_{\frakz_\cX}}{\partial w}(w) - wf_{\frakz_\cX}(w),
\]
 one can write
\begin{equation}  \label{stein_equation_used}
P\big(\barW   - \barDelta_{\frakz_\cX}\leq \frakz_\cX  \big) - \bE[ \Phi(\frakz_\cX)]= F + \bE\Big[\bar{\Delta}_{\frakz_\cX} f_{\frakz_\cX} (\barW)\Big],
\end{equation}
where
\begin{align}
F &\equiv \bE\big[
I(\barW  - \barDelta_{\frakz_\cX} \leq \frakz_\cX )  - \Phi({\frakz_\cX})- \bar{\Delta}_{\frakz_\cX} f_{\frakz_\cX} (\barW )
\big]  \notag\\
&=
\bE\Big[
\frac{\partial f_{\frakz_\cX}}{\partial w} (\barW - \bar{\Delta}_{\frakz_\cX})
- \barW f_{\frakz_\cX} (\barW- \barDelta_{\frakz_\cX})
+ \barDelta_{\frakz_\cX} \Big( f_{\frakz_\cX} (\barW - \barDelta_{\frakz_\cX}) - f_{\frakz_\cX} (\barW ) \Big)
\Big].  \label{R_def}
\end{align}
Hence, with the definition of $\barDelta_{\frakz_\cX}$ in \eqref{barDeltaz_def}  and the property that $|f_{\frakz_\cX}| \leq 0.63$ from  \lemref{fxfxprime_bdd}$(i)$ below,  continuing from \eqref{stein_equation_used} we  obtain 
\begin{equation} \label{conseq_stein_eqt}
\big|P\big(\barW  - \barDelta_{\frakz_\cX}\leq \frakz_\cX  \big) - \bE[\Phi(\frakz_\cX)] \big| \leq |F| + 0.63 \bE[|\bar{D}_1|] + \frac{1}{2} \bigg| \bE\big[ \frakz_\cX \bar{D}_2   f_{\frakz_\cX} (\barW )   \big] \bigg|.
\end{equation}
Hence, to prove \lemref{prelim_BE}, it remain to establish  \eqref{RCI_sqeeuze_application} and \eqref{F_bdd_borrow} below:
\begin{multline} \label{RCI_sqeeuze_application}
 P\bigg(\text{sgn}(\frakz_\cX) \cdot \frakz_\cX \bigg(1+\frac{D_2}{2} - \frac{D_2^2}{2}\bigg)  \leq \text{sgn}(\frakz_\cX) \cdot (W  + D_1) \leq  \text{sgn}(\frakz_\cX) \cdot  \frakz_\cX \bigg(1 + \frac{D_2}{2}\bigg)  \bigg)  \\
  \leq  C \Bigg\{\bE[\bar{D}_2^2] +  \frac{\bE[|h|^3](1 - 2p + 2p^2)}{\sigma_h^3 \sqrt{N(1 - p)}}+ \frac{\sum_{\boldi \in I_{n, m}} \bE\Big[|h(\bX_\boldi)|\sqrt{\bE[ (\bar{D}_1- \bar{D}_1^{(\boldi)})^2 | \cX] }\Big]}{\sigma_h\sqrt{ {n \choose m}}} \Bigg\}\\
+ \sum_{k=1}^2 P(|D_k| > 1/2 ) +P(\Omega\backslash \cE_{1, \cX}  ).
\end{multline}
and
\begin{multline} \label{F_bdd_borrow}
|F| \leq C \Bigg\{ \frac{\bE[|h|^3](1 - 2p + 2p^2)}{\sigma_h^3 \sqrt{N(1 - p)}} + \bE[\bar{D}^2_2] +\\
 \bE\Bigg[\sqrt{E[\bar{D}_1^2 \mid \cX]} +
\frac{\sum_{\boldi \in I_{n, m}} |h(\bX_\boldi)|\sqrt{\bE[ (\bar{D}_1- \bar{D}_1^{(\boldi)})^2 | \cX] }}{\sigma_h\sqrt{{n \choose m}}} \Bigg] + P(\Omega\backslash \cE_{1, \cX}  )\Bigg\}
\end{multline}
Combining \eqref{Tsn_BE_bound_template}, \eqref{conseq_stein_eqt}-\eqref{F_bdd_borrow} and $\bE[|\bar{D}_1|] \leq  \bE[\sqrt{\bE[ \bar{D}_1^2\mid \cX]}]$,  \lemref{prelim_BE} is then established.

In what follows,
we will respectively prove  \eqref{RCI_sqeeuze_application}  and \eqref{F_bdd_borrow} but with some details responsibly skipped. This is because, 
as mentioned in \eqref{obs_as_sn_stat},  $\Tsn$ is an example of a Studentized nonlinear statistic in \eqref{sn_stat_generic_form} when conditioned on $\cX$ for any 
\begin{equation} \label{cx_event_condition}
\cX \in \cE_{1, \cX} ,
\end{equation}
and many  mathematical results already developed in \citet{leung2024another} can be directly  borrowed;  the readers will be given exact pointers to the corresponding sections in \citet{leung2024another} for these skipped details. We will focus on working out derivations that are not covered by \citet{leung2024another}. In the framework of \citet{leung2024another}, for each $\boldi \in \Inm$, one also has to define the ``leave-one-out'' version of the remainder $D_2$ in the normalizer of $\Tsn$ akin to $D_1^{(\boldi)}$ and $\barD_1^{(\boldi)}$ in \eqref{D1_leave_one_out} and \eqref{barD1_leave_one_out}. However,  
  Since $D_2$  doesn't depend on the Bernoulli samplers in $\cZ$ indexed by ${\bf i} \in \Inm$,  unlike $D_1^{(\boldi)}$ and $\barD_1^{(\boldi)}$, we can trivially take the ``leave-one-out'' versions of $D_2$ and $\barD_2$ to  be 
\begin{equation} \label{leave_one_out_D2s}
D_2^{(\boldi)}\equiv D_2 \text{ and } \barD_2^{(\boldi)}\equiv \barD_2.
\end{equation}
We will also define the leave-one-out version of $\barW$ in \eqref{Ws_def} for any $\boldi \in \Inm$ as
\begin{equation} \label{barW_leave_one_out}
\barW^{(\boldi)} = \sum_{\boldj \in \Inm, \boldj \neq \boldi} \barxi_\boldj.
\end{equation}
Moreover, the following four  inequalities will come in handy to supplement the developed results in \citet{leung2024another}. The first is 
\begin{multline}\label{beta2beta3les3}
\sum_{\boldi \in \Inm} \bE[ \xi_\boldi^2 I(| \xi_\boldi| > 1) \mid \cX] 
+ \sum_{\boldi \in \Inm} \bE[ |\xi_\boldi|^3 I(| \xi_\boldi| \leq 1) \mid \cX] \\
 \leq \sum_{\boldi \in \Inm}  \bE[ |\xi_\boldi|^3 \mid \cX] 
 =  \frac{U_{|h|^3} (1 - 2p + 2p^2)}{\bar{U}_{h^2}^{3/2} \sqrt{N(1 - p)}} \leq \frac{2^{3/2}U_{|h|^3} (1 - 2p + 2p^2)}{\sigma_h^3 \sqrt{N(1 - p)}} ,
\end{multline}
coming from  the third moment calculation in \eqref{cX_cond_moment_properties_of_xi} and $\sigma_h^2/2 \leq \bar{U}_{h^2}$. Next is the Cauchy inequality 
\begin{multline} \label{supplementi_borrow_2}
\bE\Big[| \barxi_\boldi e^{\pm\bar{W}^{(\boldi)}/2}  (\bar{D}_1 - \bar{D}_1^{(\boldi)})| \mid \cX\Big] 
\leq   \sqrt{\bE[e^{\pm \barW^{(\boldi)}}\mid \cX] \bE [ \barxi_\boldi^2 \mid \cX ] } \sqrt{\bE[ (\bar{D}_1- \bar{D}_1^{(\boldi)})^2 | \cX] } \\
\leq C \frac{|h(\bX_\boldi)| \sqrt{\bE[ (\bar{D}_1- \bar{D}_1^{(\boldi)})^2 | \cX] }}{\sigma_h \sqrt{{n \choose m}}}
 \text{ for } \cX \in \cE_{1, \cX},
\end{multline}
where in the last inequality we have used the second moment calculation in  \eqref{cX_cond_moment_properties_of_xi}, $\sigma_h^2/2 \leq \bar{U}_{h^2}$  and  \lemref{Bennett}  in light of  \eqref{SN_condition_satisfied_1}. The third is 
\begin{equation} \label{supplementi_borrow_3}
\bE[e^{\pm \barW} | \cX] \leq 8.15 \text{ for } \cX \in \cE_{1, \cX},
\end{equation}
by virtue of  \lemref{Bennett} in light of  \eqref{SN_condition_satisfied_1}.
The last one is 
\begin{align} \label{supplementi_borrow_4}
\bE\Big[ (1+ e^{\pm \bar{W}^{(\boldi)}} )|\bar{D}_1 - \bar{D}_1^{(\boldi)}| \mid \cX\Big] 
&\leq \Big( 1+\sqrt{\bE [ e^{\pm 2\barW^{(\boldi)}}\mid \cX]}\Big) \sqrt{\bE[ (\bar{D}_1- \bar{D}_1^{(\boldi)})^2 | \cX] } \notag\\
&\leq C \sqrt{\bE[ (\bar{D}_1- \bar{D}_1^{(\boldi)})^2 | \cX] } \text{ for } \cX \in \cE_{1, \cX} ,
\end{align}
where we have again applied  \lemref{Bennett}  in light of \eqref{SN_condition_satisfied_1}.

\subsection{Proof of \eqref{RCI_sqeeuze_application}} \label{app:RCI_application}
Given how  $D_2^{(\boldi)}$ and $\bar{D}_2^{(\boldi)}$ are defined in \eqref{leave_one_out_D2s} and  $\text{sgn}(\frakz_\cX) \cdot \frakz_\cX \geq 0$, 
one can 
follow \emph{exactly the same} arguments in  \citet[Appendix C.1]{leung2024another} to show that, under condition \eqref{cx_event_condition},
\begin{multline*}
P\bigg(\text{sgn}(\frakz_\cX) \cdot \frakz_\cX \cdot \bigg(1+\frac{D_2}{2} - \frac{D_2^2}{2}\bigg)  \leq \text{sgn}(\frakz_\cX) \cdot (W  + D_1) \leq  \text{sgn}(\frakz_\cX) \cdot  \frakz_\cX \cdot \bigg(1 + \frac{D_2}{2}\bigg)   \mid \cX  \bigg) \\
\leq  \sum_{k=1}^2 P(|D_k| > 1/2 \mid \cX ) + 
C \Bigg\{\sum_{\boldi \in \Inm} \bE[ \xi_\boldi^2 I(| \xi_\boldi| > 1) \mid \cX] 
 \\
+ \sum_{\boldi \in \Inm} \bE[ |\xi_\boldi|^3 I(| \xi_\boldi| \leq 1) \mid \cX]  + \bE\Big[\big(1 + e^{\text{sgn}(\frakz_\cX)\barW} \big) \bar{D}_2^2 \mid  \cX\Big] 
+ \sum_{\boldi \in \Inm}\bE\big[| \barxi_\boldi e^{\text{sgn}(\frakz_\cX) \bar{W}^{(\boldi)}/2}  (\bar{D}_1 - \bar{D}_1^{(\boldi)})| \mid \cX\big] \Bigg\},
\end{multline*}
which can be further simplified by    \eqref{beta2beta3les3}-\eqref{supplementi_borrow_3} as 
\begin{multline}   \label{sandwich_prob_ineq}
P\bigg(\text{sgn}(\frakz_\cX) \cdot \frakz_\cX \cdot \bigg(1+\frac{D_2}{2} - \frac{D_2^2}{2}\bigg)  \leq \text{sgn}(\frakz_\cX) \cdot (W  + D_1) \leq  \text{sgn}(\frakz_\cX) \cdot  \frakz_\cX \cdot \bigg(1 + \frac{D_2}{2}\bigg)   \mid \cX  \bigg) \\ 
\leq   \sum_{k=1}^2 P\bigg(|D_k| > \frac{1}{2}\mid \cX \bigg) + 
C\Bigg\{ \bar{D}_2^2 + \frac{U_{|h|^3} (1 - 2p + 2p^2)}{\sigma_h^3 \sqrt{N(1 - p)}} +   \frac{\sum_{\boldi \in I_{n, m}} |h(\bX_\boldi)|\sqrt{\bE[ (\bar{D}_1- \bar{D}_1^{(\boldi)})^2 | \cX] }}{\sigma_h\sqrt{ {n \choose m}}}  \Bigg\}\\
 \text{ for } \cX \in \cE_{1, \cX}.
\end{multline}
 To form the bound on the marginal probability, we write
\begin{multline*}
P\bigg(\text{sgn}(\frakz_\cX) \frakz_\cX  \bigg(1+\frac{D_2}{2} - \frac{D_2^2}{2}\bigg)  \leq \text{sgn}(\frakz_\cX)  (W  + D_1) \leq  \text{sgn}(\frakz_\cX)   \frakz_\cX  \bigg(1 + \frac{D_2}{2}\bigg)     \bigg) \leq  P(\Omega\backslash \cE_{1, \cX}  ) + \\
 \bE\bigg[ P\bigg(\text{sgn}(\frakz_\cX) \frakz_\cX  \bigg(1+\frac{D_2}{2} - \frac{D_2^2}{2}\bigg)  \leq \text{sgn}(\frakz_\cX)  (W  + D_1) \leq  \text{sgn}(\frakz_\cX)   \frakz_\cX  \bigg(1 + \frac{D_2}{2}\bigg)    \mid \cX \bigg)I( \cE_{1, \cX} )\bigg] ;
\end{multline*}
applying \eqref{sandwich_prob_ineq} to the prior display  yields  \eqref{RCI_sqeeuze_application}.

\subsection{Proof of  \eqref{F_bdd_borrow}}  \label{app:R_bdd_sec}
Let $\{Z_\boldi^*\}_{\boldi \in \Inm}$ be independent copies of $\{Z_\boldi\}_{\boldi \in \Inm}$. 
For each ${\bf i} \in \Inm$, define
\[
D_{1, \boldi^*} = D_1\Big(\{Z_\boldj\}_{\boldj \in \Inm, \boldj \neq \boldi}, Z_{\boldi}^* ; \cX\Big),
\]
which is constructed in exactly same way as $D_1$ in \eqref{generic_D1}, except that $Z_\boldi^*$  has replaced $Z_\boldi$  as the input, as well as its censored version
\[
\bar{D}_{1, \boldi^*}  \equiv  \bar{D}_{1, \boldi^*} I\bigg(-\frac{1}{2}\leq \bar{D}_{1, \boldi^*} \leq \frac{1}{2}\bigg)- \frac{1}{2} I\bigg(\bar{D}_{1, \boldi^*} < -\frac{1}{2}\bigg)+ \frac{1}{2 }I\bigg(\bar{D}_{1, \boldi^*}>\frac{1}{2}\bigg);
\]
note that $D_{1, \boldi^*} $ is not to be confused with the ``leave-one-out'' $D_1^{(\boldi)}$ defined in \eqref{D1_leave_one_out}.  From these, we also define
\[
\bar{\Delta}_{\frakz_\cX, \boldi^*} \equiv  \frakz_\cX \bar{D}_2/2 -   \bar{D}_{1, \boldi^*},
\]
which is similar to $\bar{\Delta}_{\frakz_\cX}$ in \eqref{barDeltaz_def} but has $\bar{D}_1$ replaced by $\bar{D}_{1, \boldi^*}$. 
Moreover, with the censored $\barxi_\boldi $  from \eqref{xi_censored},
we  define the 
 the conditional \emph{K function}
\begin{align*}
\bark_\boldi(t) &\equiv \bE\Big[\barxi_\boldi \big(I(0 \leq t \leq \barxi_\boldi) - I( \barxi_\boldi\leq t <0)\big)  \mid \bX_\boldi\Big] \\
&= \bE\Big[\barxi_\boldi \big(I(0 \leq t \leq \barxi_\boldi) - I( \barxi_\boldi\leq t <0)\big)  \mid \cX\Big],
\end{align*}
which has the properties
\begin{multline} \label{k_fn_properties}
 \int_{-\infty}^\infty \bark_\boldi(t)  dt =   \int_{-1}^1 \bark_\boldi(t)  dt = \bE[ \barxi_\boldi^2 \mid \cX] \text{ and }\\
  \int_{-\infty}^\infty |t|\bark_\boldi(t)  dt =    \int_{-1}^1 |t|\bark_\boldi(t)  dt = \frac{\bE[ |\barxi_\boldi|^3\mid \cX]}{2}
\end{multline}

Now we begin to establish \eqref{F_bdd_borrow}. 
From the definition of $F$ in \eqref{R_def}, we  write
\begin{multline} \label{F_break_down}
F= \\
\underbrace{\bE\bigg[\bigg\{\frac{\partial f_{\frakz_\cX}}{\partial w}(\barW - \bar{\Delta}_{\frakz_\cX})
- \barW f_{\frakz_\cX} (\barW- \barDelta_{\frakz_\cX})
+ \barDelta_{\frakz_\cX} \Big( f_{\frakz_\cX} (\barW - \barDelta_{\frakz_\cX}) - f_{\frakz_\cX} (\barW ) \Big)\bigg\} \cdot I(\cE_{1, \cX})\bigg]}_{\equiv F_1}\\
+ \underbrace{\bE\bigg[ \Big\{I(\barW  - \barDelta_{\frakz_\cX} \leq \frakz_\cX )  - \Phi({\frakz_\cX})- \bar{\Delta}_{\frakz_\cX} f_{\frakz_\cX} (\barW )\Big\} \cdot I(\Omega \backslash \cE_{1, \cX} )\bigg] }_{\equiv F_2}.
\end{multline}
Next, we further break down $F_1$ by noticing two facts: By independence and the fundamental theorem of calculus, one has
\begin{multline}\label{FTC_consequence}
\bE\bigg[ \int_{-1}^1 \frac{\partial f_{\frakz_\cX}}{\partial w} \Big(\barW^{(\boldi)}  - \bar{\Delta}_{\frakz_\cX, \boldi^*} +t  \Big) \bark_\boldi(t) dt\mid \cX\bigg] = \\
\bE\bigg[ 
\barxi_\boldi \bigg\{ f_{\frakz_\cX} \Big( \barW - \bar{\Delta}_{\frakz_\cX, \boldi^*} \Big) -  f_{\frakz_\cX} \Big( \barW^{(\boldi)}  - \bar{\Delta}_{\frakz_\cX, \boldi^*} \Big)\bigg\}
\mid \cX\bigg] 
\end{multline}
On the other hand, we have 
\begin{equation} \label{sum_xi_sq_to_1}
\sum_{\bold i \in \Inm} \bE[ \xi_\boldi^2 \mid \cX] = 1 \text{ under condition } \eqref{cx_event_condition} 
\end{equation}
 as per the conditional moment calculations in  \eqref{SN_condition_satisfied_1}. Combining the two observations in  \eqref{FTC_consequence} and \eqref{sum_xi_sq_to_1}, as well as the fact $\int_{-\infty}^\infty \bark_\boldi(t)  dt = \bE[ \barxi_\boldi^2 \mid \cX] $ in \eqref{k_fn_properties}, one can write $F_1$ as
\begin{equation} \label{four_R_sum_new}
F_1 = \sum_{k= 1}^4 \bE[ F_{1k}  I(\cE_{1, \cX} )]   
\end{equation}
where
\[
F_{11} \equiv \sum_{\boldi \in \Inm}\bE\Bigg[\int_{-1}^1 \bigg\{\frac{\partial f_{\frakz_\cX}}{\partial w} \bigg( \barW-  \bar{\Delta}_{\frakz_\cX}\bigg) - \frac{\partial f_{\frakz_\cX}} {\partial w}\bigg( \barW^{(\boldi)}  -  \bar{\Delta}_{\frakz_\cX, \boldi^*}  + t\bigg)\bigg\} \bark_\boldi(t) dt  \mid \cX\Bigg], 
\]
\begin{multline*} \label{R2_def}
F_{12} \equiv  \sum_{\boldi \in \Inm} \bE \bigg[ (\xi_\boldi^2- 1)I(|\xi_\boldi|> 1) \mid \cX\bigg]  \cdot
\bE \bigg[\frac{\partial f_{\frakz_\cX}}{\partial w} (\barW-  \bar{\Delta}_{\frakz_\cX})\mid \cX\bigg]  \\
- \sum_{\boldi \in \Inm}\bE \bigg[\barxi_{\boldi} f_{\frakz_\cX} \bigg(\barW^{(\boldi)}   - \bar{\Delta}_{\frakz_\cX, \boldi^*} \bigg) \mid \cX \bigg],
\end{multline*}
\[
F_{13}\equiv\Bigg\{ - \sum_{\boldi \in \Inm}\bE\Bigg[\barxi_{\boldi}  \Bigg\{ f_{\frakz_\cX}\bigg(\barW -  \bar{\Delta}_{\frakz_\cX}\bigg) -
 f_{\frakz_\cX} \bigg(\barW- \bar{\Delta}_{\frakz_\cX, \boldi^*} \bigg) \Bigg\} \mid \cX\Bigg]\Bigg\}
\]
and
\[
F_{14} \equiv\bE\Bigg[ \bar{\Delta}_{\frakz_\cX} \int_0^{-  \bar{\Delta}_{\frakz_\cX}} \frac{\partial f_{\frakz_\cX} }{\partial w} \bigg(\barW+ t\bigg) dt \mid \cX \Bigg].
\]
We will now establish bounds for $F_{11}$, $F_{12}$, $F_{13}$ and $F_{14}$.

Given how we define $D_2^{(\boldi)}$ and $\bar{D}_2^{(\boldi)}$ in \eqref{leave_one_out_D2s},  borrowing  essentially  the exact  arguments in \citet[Appendix C.2]{leung2024another},\footnote{The terms $F_{11}$, $F_{13}$ and $F_{14}$ are completely analogous to the terms $R_1$, $R_3$ and $R_4$ respectively bounded in Appendices C.2.1, C.2.3 and C.2.4 of  \citet{leung2024another}. Under \eqref{cx_event_condition}, all the assumptions for the (fairly laborious) proof arguments there are met, except that our threshold $\frakz_\cX$ may have a negative sign, whereas the threshold appearing in \citet[Appendix C]{leung2024another} is assumed to be non-negative  in a ``without-loss-of-generality'' manner. However, the bound in \lemref{fxfxprime_bdd}$(iii)$ for $f_z$ when $z < -1$ can be leveraged to fully reproduce those arguments; we leave the details to the interested readers who may revisit \citet[Appendix C.2]{leung2024another}.} one can establish the bounds
\begin{multline*} \label{F11_abs_bdd}
|F_{11}| \leq C \Bigg\{ \sum_{\boldi \in \Inm} \bE[ \xi_\boldi^2 I(| \xi_\boldi| > 1) \mid \cX] 
+ \sum_{\boldi \in \Inm} \bE[ |\xi_\boldi|^3 I(| \xi_\boldi| \leq 1) \mid \cX] +\\
 \sum_{\boldi \in \Inm} \bE\big[ \barxi_\boldi^2 \mid \cX\big] \cdot \bE\Big[ (1+ e^{\text{sgn}(\frakz_\cX)\bar{W}^{(\boldi)}} )|\bar{D}_1 - \bar{D}_1^{(\boldi)}| \mid \cX\Big]
+ \sum_{\boldi \in \Inm}  \bE\Big[| \barxi_\boldi e^{\text{sgn}(\frakz_\cX)\bar{W}^{(\boldi)}/2}  (\bar{D}_1 - \bar{D}_1^{(\boldi)})| \mid \cX\Big]  \Bigg\},
\end{multline*}
\begin{equation*} \label{F13_abs_bdd}
|F_{13}| \leq  C   \sum_{\boldi \in \Inm}  \bE\Big[| \barxi_\boldi  (\bar{D}_1 - \bar{D}_1^{(\boldi)})| \mid \cX\Big] 
\end{equation*}
and
\begin{equation*} \label{F14_abs_bdd}
|F_{14}| \leq  C  \Bigg\{ \sqrt{E[\bar{D}_1^2 \mid \cX]} + \bar{D}_2^2 + \bE[e^{\text{sgn}(\frakz_\cX) \barW} \bar{D}_2^2 \mid \cX ] \Bigg\}
, \text{ all under }\eqref{cx_event_condition}.
\end{equation*}
Applying  $\bar{U}_{h^2} \geq \sigma_h^2/2$, \eqref{beta2beta3les3}-\eqref{supplementi_borrow_4} and the  second moment calculation in  \eqref{cX_cond_moment_properties_of_xi}   to the above bounds for $F_{11}$, $F_{13}$ and $F_{14}$ further gives
\begin{equation} \label{F11_abs_bdd}
|F_{11}| \leq C \Bigg\{\frac{U_{|h|^3} (1 - 2p + 2p^2)}{\sigma_h^3 \sqrt{N(1 - p)}}+ \sum_{\boldi \in \Inm} 
\frac{|h(\bX_\boldi)| \sqrt{\bE[ (\bar{D}_1- \bar{D}_1^{(\boldi)})^2 | \cX] }}{\sigma_h \sqrt{{n \choose m}}}
 \Bigg\},
\end{equation}
\begin{equation} \label{F13_abs_bdd}
|F_{13}| \leq  C \sum_{\boldi \in \Inm}  \frac{|h(\bX_\boldi)| \sqrt{\bE[ (\bar{D}_1- \bar{D}_1^{(\boldi)})^2 | \cX] }}{\sigma_h \sqrt{{n \choose m}}}
\end{equation}
and
\begin{equation} \label{F14_abs_bdd}
|F_{14}| \leq  C  \Bigg\{ \sqrt{E[\bar{D}_1^2 \mid \cX]} + \bar{D}_2^2  \Bigg\}
, \text{ all under }\eqref{cx_event_condition},
\end{equation}
where we have also used  $\bE\big[ \barxi_\boldi^2 \mid \cX\big] \leq \sqrt{\bE\big[ \barxi_\boldi^2 \mid \cX\big]}$ for \eqref{F11_abs_bdd}, because $|\barxi_\boldi| \leq 1$.

For $F_{12}$, 
the following bound 
\begin{equation} \label{F12_abs_bdd}
|F_{12}| \leq  C\frac{U_{|h|^3} (1 - 2p + 2p^2)}{\sigma_h^3 \sqrt{N(1 - p)}}.
\end{equation}
can also be established as follows: 
 First, using   $|\frac{\partial f_{\frakz_\cX}}{\partial w}| \leq 1$ (\lemref{fxfxprime_bdd}$(i)$) and  \eqref{beta2beta3les3},  we have
\begin{multline}
 \bigg| \sum_{\boldi \in \Inm} \bE [ (\xi_\boldi^2- 1)I(|\xi_\boldi|> 1) \mid \cX]  \cdot
\bE \bigg[\frac{\partial f_{\frakz_\cX}}{\partial w}(\barW-  \bar{\Delta}_{\frakz_\cX})\mid \cX \bigg]\bigg| \\
\leq  \sum_{\boldi \in \Inm}\bE \Big[ \xi_\boldi^2I(|\xi_\boldi|> 1) \mid \cX\Big]  
\leq \frac{2^{3/2}U_{|h|^3} (1 - 2p + 2p^2)}{\sigma_h^3 \sqrt{N(1 - p)}} \label{R2_bdd_first}.
\end{multline}
Second, using that $|f_{\frakz_\cX}| \leq 0.63$  from \lemref{fxfxprime_bdd}$(i)$,  we have
\begin{align}
& \Bigg| \sum_{\boldi \in \Inm}\bE \Big[\barxi_{\boldi} f_{\frakz_\cX} (\barW^{(\boldi)}   - \bar{\Delta}_{\frakz_\cX, \boldi^*} ) \mid \cX \Big] \Bigg| \notag\\
&=  \Bigg| \sum_{\boldi \in \Inm}\bE [\barxi_{\boldi} \mid \cX ] \bE\Big[ f_{\frakz_\cX} (\barW^{(\boldi)}    - \bar{\Delta}_{{\frakz_\cX}, \boldi^*} ) \mid \cX \Big] \Bigg|  \notag\\
&\leq 0.63 \sum_{\boldi \in \Inm} \Big| \bE \big[\barxi_{\boldi} \mid \cX \big]  \Big| \notag\\
&= 0.63 \sum_{\boldi \in \Inm} \Big| \bE \big[(\xi_\boldi - 1)I(\xi_\boldi >1) + (\xi_\boldi +1) I(\xi_\boldi <-1)\mid \cX \big]  \Big| \text{ by }\bE[\xi_\boldi \mid \cX] = 0 \text{ in }\eqref{cX_cond_moment_properties_of_xi}\notag\\
& \leq 0.63 \sum_{\boldi \in \Inm}  \bE \Big[|\xi_{\boldi}| I(|\xi_{\boldi}| >1) \mid \cX \Big] \notag\\
&\leq 0.63 \sum_{\boldi \in \Inm}\bE [ |\xi_{\boldi}|^3\mid \cX] \notag\\
& = 0.63 \frac{U_{|h|^3} (1 - 2p + 2p^2)}{\bar{U}_{h^2}^{3/2} \sqrt{N(1 - p)}} \leq 1.26 \frac{\sqrt{2}U_{|h|^3} (1 - 2p + 2p^2)}{\sigma_h^3\sqrt{N(1 - p)}}, \label{R2_bdd_second}
\end{align}
where \eqref{cX_cond_moment_properties_of_xi} is used in the last equality. 
So \eqref{R2_bdd_first} and \eqref{R2_bdd_second} together give \eqref{F12_abs_bdd}. Collecting \eqref{four_R_sum_new}-\eqref{F12_abs_bdd}, we can get
\begin{multline} \label{F1_bdd}
|F_1| \leq C \Bigg\{ \frac{\bE[|h|^3](1 - 2p + 2p^2)}{\sigma_h^3 \sqrt{N(1 - p)}} + \bE[\bar{D}^2_2] +\\
 \bE\Bigg[\sqrt{E[\bar{D}_1^2 \mid \cX]} +
\frac{\sum_{\boldi \in I_{n, m}} |h(\bX_\boldi)|\sqrt{\bE[ (\bar{D}_1- \bar{D}_1^{(\boldi)})^2 | \cX] }}{\sigma_h\sqrt{{n \choose m}}} \Bigg] \Bigg\}.
\end{multline}
For $F_2$ in \eqref{F_break_down}, because  $|f_{\frakz_\cX}| \leq 0.63$  from \lemref{fxfxprime_bdd}$(i)$ and $\sup_{z, w \in \bR}|zf_z(w)| \leq 1$ in  \lemref{fxfxprime_bdd}$(iv)$, it is easy to see that
\begin{multline*}
I(\barW  - \barDelta_{\frakz_\cX} \leq \frakz_\cX )  - \Phi({\frakz_\cX})- \bar{\Delta}_{\frakz_\cX} f_{\frakz_\cX}(\barW ) \\
= 
 I(\barW  - \barDelta_{\frakz_\cX}\leq {\frakz_\cX} )  - \Phi(\frakz_\cX) +  f_{\frakz_\cX} (\barW ) \bar{D}_1- \frac{\frakz_\cX f_{\frakz_\cX} (\barW )\bar{D}_2}{2}
\end{multline*}
 can be bounded in absolute value by $3$ to give that
\begin{equation} \label{F2_bdd}
 |F_2| \leq 3 P(\Omega\backslash \cE_{1, \cX} ).
\end{equation}
Lastly, combining \eqref{F_break_down}, \eqref{F1_bdd} and \eqref{F2_bdd} we obtain \eqref{F_bdd_borrow}.

 \section{Proof of  \lemref{TD1_bdd}, the bounds on  $T(D_1)$ } \label{app:bounding_T_D1}
 In this section, we will define the U-statistic with kernel $|h|$
\[
 U_{|h|}  \equiv {n \choose m}^{-1}\sum_{\boldj \in I_{n,m}} |h| (X_\boldj).
\]
 As will be seen later, 
 for the proofs of both  \lemref{TD1_bdd}$(i)$ and $(ii)$, we have to use the conditional expectation 
 estimate 
 \begin{equation}\label{moment_est_for_D1_bdd}
\bE\Big[ (\sum_{\boldj \in I_{n, m}} \xi_\boldj  )^2   (\sum_{\boldk \in I_{n, m}} (Z_\boldk - p) )^2 \mid \cX \Big] \leq \frac{ (1 + 3N)U_{h^2}}{\bar{U}_{h^2}},
\end{equation}
 which can be proved as follows: By the independence among $\{Z_\boldj\}_{\boldj \in \Inm}$ and  that $\bE[ \xi_\boldj| \cX] = 0$ from \eqref{cX_cond_moment_properties_of_xi}, 
\begin{align*}
&\bE\bigg[ \bigg(\sum_{\boldj \in I_{n, m}} \xi_\boldj \bigg)^2   \bigg(\sum_{\boldk \in I_{n, m}} (Z_\boldk - p) \bigg)^2 \mid \cX\bigg]\\
& = \bE\Big[ \sum_{\boldj \in I_{n, m}} \xi_\boldj^2 (Z_\boldj - p)^2 
+  \sum_{\boldj \in I_{n, m}} \xi_\boldj^2  \sum_{\boldk \neq \boldj}  (Z_\boldk - p)^2
+ 2 \sum_{\boldj \in I_{n, m}} \xi_\boldj   (Z_\boldj - p) \bigg(\sum_{\boldk \neq \boldj}  \xi_\boldk (Z_\boldk - p) \bigg) \mid \cX\Big]\\
&=  \sum_{\boldj \in I_{n, m}} \frac{p(1- p)( 1 + 3p^2 - 3p) h^2(\bX_\boldj)}{\bar{U}_{h^2}N(1- p)} 
+ \bigg({n \choose m} - 1\bigg)\frac{ p^2(1- p)^2  \sum_{\boldj \in I_{n, m}} h^2(\bX_\boldj)}{\bar{U}_{h^2}N(1-p)}\\
& \hspace{8cm}
+ \frac{ 2 p^2(1-p)^2}{\bar{U}_{h^2} N (1 - p)}  \sum_{\boldj \in I_{n, m}}   h(\bX_\boldj) \bigg(\sum_{\boldk \neq \boldj } h(\bX_\boldk)\bigg)\\
&= \frac{U_{h^2}}{\bar{U}_{h^2}}\bigg( 1 + 3p^2 - 3p
+ \bigg({n \choose m} - 1\bigg)p(1- p) \bigg)
+ \frac{2 N(1-p)}{\bar{U}_{h^2} {n \choose m}^2 }  \sum_{\boldj \in I_{n, m}}   h(\bX_\boldj) \bigg(\sum_{\boldk \neq \boldj } h(\bX_\boldk)\bigg)\\
&\leq \frac{U_{h^2}}{\bar{U}_{h^2}}\Big( 1 + 3p^2 - 3p + N(1- p)\Big) + \frac{2 N(1- p) U_{|h|}^2 }{\bar{U}_{h^2}} \\
&\leq \frac{U_{h^2} (  1 + 3p^2 - 3p + 3N (1-p) )}{\bar{U}_{h^2}} \leq \frac{ (1 + 3N)U_{h^2}}{\bar{U}_{h^2}},
\end{align*}
where the second last inequality uses 
Jensen's inequality $U^2_{|h|} \leq U_{h^2}$. 


 \subsection{Proof of \lemref{TD1_bdd}$(i)$} We separate the proof into three steps:
 
 \subsubsection{Bounding $\sqrt{\bE[D_1^2 | \cX]} 
$}
We will first show the bound
\begin{equation} \label{1st_setp_D1_bdd}
\sqrt{\bE[D_1^2 | \cX]} \leq  \frac{2 \sqrt{ (1 + 3N)U_{h^2}}}{\sqrt{\bar{U}_{h^2}}N} +   \frac{2\sqrt{2N} \big|U_n  \big|  }{ \sigma_h\sqrt{(1 -p)}}.
\end{equation}
 First, by \eqref{Ncminus_N_bdd}, $N/2\leq \hat{N}_c$ and $\sigma_h^2/2 \leq \bar{U}_{h^2} $, with the form of $D_1$ in \configref{Nlln}$(ii)$,  we have
\[
|D_1| \leq   \frac{2 \big|\big(\sum_{\boldj \in I_{n, m}} \xi_\boldj \big) \big(\sum_{\boldk \in I_{n, m}} (Z_\boldk - p) \big)\big|}{N} 
+  \frac{2\sqrt{2N} \big|U_n  \big|  }{ \sigma_h\sqrt{(1 -p)}}.
\]
By the moment estimate \eqref{moment_est_for_D1_bdd}, we get \eqref{1st_setp_D1_bdd}.

 \subsubsection{Bounding $\sqrt{\bE[ (D_1  - D_1^{(\boldi)})^2 |\cX]}$}
We now  claim that the following conditional $\ell_2$-norm bound holds:
\begin{equation} \label{D1minusD1i_norm}
\sqrt{\bE[ (D_1  - D_1^{(\boldi)})^2 |\cX]} \leq \frac{2|h(\bX_\boldi)| \sqrt{1 + N}}{N\sqrt{ \bar{U}_{h^2} {n \choose m} } }  +  \frac{4\sqrt{ U_{h^2}}}{\sqrt{ \bar{U}_{h^2}  N {n \choose m}} }
+ \frac{4 |U_n| }{\sqrt{\bar{U}_{h^2}(1-p){n \choose m}}}.
\end{equation}
First, from  \configref{Nlln}$(ii)$ and $(iii)$, write
\begin{align}
D_1  - D_1^{(\boldi)} 
&= \xi_\boldi \frac{N - \hat{N}_c}{\hat{N}_c} +\sum_{\boldj \neq \boldi} \xi_\boldj \Bigg( \frac{N - \hat{N}_c}{\hat{N}_c} -  \frac{N - \hat{N}_c^{(\boldi)}}{\hat{N}_c^{(\boldi)}} \Bigg) +  \bigg( \frac{1}{\hat{N}_c} - \frac{1}{\hat{N}_c^{\boldi}} \bigg) \frac{ N^{3/2}U_n}{ \sqrt{ \bar{U}_{h^2} (1 - p)}  } \notag \\
&=  \xi_\boldi \frac{N - \hat{N}_c}{\hat{N}_c} + \bigg( \sum_{\boldj \neq \boldi} \xi_\boldj  + \frac{ \sqrt{N}U_n}{ \sqrt{ \bar{U}_{h^2} (1 - p)}  } \bigg)  \Bigg( \frac{N(\hat{N}_c^{(\boldi)} - \hat{N}_c)}{ \hat{N}_c\hat{N}_c^{(\boldi)} } \Bigg) \label{D1minusD1boldi_form}.
\end{align}
Hence, by the triangular inequality, we have
\begin{multline} \label{D1minusD1i_in_two_parts}
\sqrt{\bE[(D_1  - D_1^{(\boldi)})^2 \mid \cX]} \leq  \sqrt{ \bE\bigg[  \bigg(\xi_\boldi \frac{N - \hat{N}_c}{\hat{N}_c}\bigg)^2 \mid \cX\bigg]}  
+ \\
 \sqrt{ \bE\bigg[\Big(\sum_{\boldj \neq \boldi} \xi_\boldj\Big)^2 \bigg( \frac{N(\hat{N}_c^{(\boldi)} - \hat{N}_c)}{ \hat{N}_c\hat{N}_c^{(\boldi)} } \bigg)^2 \mid \cX \bigg]} 
+  \sqrt{ \bE\bigg[ \frac{N U_n^2}{ \bar{U}_{h^2}(1-p)} \bigg( \frac{N(\hat{N}_c^{(\boldi)} - \hat{N}_c)}{ \hat{N}_c\hat{N}_c^{(\boldi)} } \bigg)^2 \mid \cX \bigg]}.
\end{multline}
We now bound the three  conditional second moments on the right hand side of \eqref{D1minusD1i_in_two_parts} in order. Since 
\begin{equation} \label{1st_2nd_moment_bdd}
\bigg|\xi_\boldi \frac{N - \hat{N}_c}{\hat{N}_c}\bigg| \leq  \frac{  |(Z_\boldi - p)( \hat{N}- N)h(\bX_\boldi)|}{\hat{N}_c\sqrt{\bar{U}_{h^2}N(1 - p)}} \leq  \frac{ 2 |(Z_\boldi - p)( \hat{N}- N)h(\bX_\boldi)|}{N^{3/2} \sqrt{ \bar{U}_{h^2}(1 - p)}}
\end{equation}
by \eqref{Ncminus_N_bdd}, we will compute the second moment of $(Z_\boldi - p)( \hat{N}- N)$. Upon expansion, 
\begin{align*}
&(Z_\boldi - p)^2 (\hat{N}- N)^2\\
 &= (Z_\boldi - p)^2 (Z_\boldi + \sum_{\boldj \neq \boldi} Z_\boldj - N)^2\\
&=  (Z_\boldi - p)^2 \bigg(Z_\boldi^2 + 2 Z_\boldi \bigg(\sum_{\boldj \neq \boldi} Z_\boldj - N\bigg) +  \bigg(\sum_{\boldj \neq \boldi} Z_\boldj - N\bigg)^2 \bigg) \\
&= (Z_\boldi - p)^2 \bigg\{Z_\boldi \bigg( 1 + 2  \bigg(\sum_{\boldj \neq \boldi} Z_\boldj - N\bigg) \bigg) +  \bigg(\sum_{\boldj \neq \boldi} Z_\boldj - N\bigg)^2 \bigg\} \\
&= Z_\boldi (1-p)^2 \bigg( 1 + 2  \bigg(\sum_{\boldj \neq \boldi} Z_\boldj - N\bigg) \bigg)  + (Z_\boldi - p)^2 \bigg(\underbrace{\sum_{\boldj \neq \boldi} Z_\boldj - N}_{=\sum_{\boldj \neq \boldi} Z_\boldj  - \bE[\sum_{\boldj \neq \boldi} Z_\boldj ]  - p }\bigg)^2.
\end{align*}
Taking expectations on both sides, we get
\begin{align}
& \bE[(Z_\boldi - p)^2 (N - \hat{N})^2] \notag\\
&= p(1- p)^2 \bigg(1 + 2(N - p - N) \bigg) + p(1- p) \bigg(\bigg({n \choose m} - 1\bigg) p(1- p) + p^2\bigg) \notag\\
&= p(1 - p) \bigg( (1 - p) (1- 2p) + N (1- p) - p(1- p) + p^2\bigg) \notag \\
&= p(1- p) (1 + 4p^2 - 4p + N(1- p))\notag\\
&=  p(1- p)((2p - 1)^2 + N(1- p)). \label{2nd_moment_blahblahblah}
\end{align}
Hence, from \eqref{1st_2nd_moment_bdd} and \eqref{2nd_moment_blahblahblah}, we get
\begin{align}
\bE\bigg[\bigg(\xi_\boldi \frac{N - \hat{N}_c}{\hat{N}_c}\bigg)^2 \mid \cX \bigg]&\leq  \frac{ 4 \bE[(Z_\boldi - p)^2 (\hat{N} - N)^2h^2(\bX_\boldi) \mid \cX]}{N^3 \bar{U}_{h^2}(1 - p)}\notag\\
&\leq \frac{4 h^2(\bX_\boldi)  ((2p -1)^2 + N(1- p))}{N^2 {n \choose m}  \bar{U}_{h^2}} \notag\\
&\leq \frac{4 h^2(\bX_\boldi)  (1 + N)}{N^2 {n \choose m}  \bar{U}_{h^2}}  \label{1st_2nd_moment_bdd_final}
\end{align}
Now we bound $ \bE\bigg[\Big(\sum_{\boldj \neq \boldi} \xi_\boldj\Big)^2 \bigg( \frac{N(\hat{N}_c^{(\boldi)} - \hat{N}_c)}{ \hat{N}_c\hat{N}_c^{(\boldi)} } \bigg)^2 \mid \cX\bigg]$. Note that, by the definitions of $\hat{N}_c$ and $\hat{N}_c^{(\boldi)}$ in \eqref{N_censored} and \eqref{Nhatc_leave_one_out}, as well as \propertyref{censoring_property}$(i)$,
\begin{equation*}
\Big(\sum_{\boldj \neq \boldi} \xi_\boldj\Big)^2 \bigg( \frac{N(\hat{N}_c^{(\boldi)} - \hat{N}_c)}{ \hat{N}_c\hat{N}_c^{(\boldi)} } \bigg)^2 
\leq \Big(\sum_{\boldj \neq \boldi} \xi_\boldj\Big)^2 \bigg( \frac{4(\hat{N}^{(\boldi)} - \hat{N})}{ N } \bigg)^2 =   \Big(\sum_{\boldj \neq \boldi} \frac{ (Z_\boldj - p)h(\bX_\boldj)}{\sqrt{\bar{U}_{h^2}N(1 - p)}}\Big)^2 \bigg( \frac{4Z_\boldi}{ N } \bigg)^2,
\end{equation*}
we can then take conditional expectation on both sides and get
\begin{align}
 \bE\bigg[\Big(\sum_{\boldj \neq \boldi} \xi_\boldj\Big)^2 \bigg( \frac{N(\hat{N}_c^{(\boldi)} - \hat{N}_c)}{ \hat{N}_c\hat{N}_c^{(\boldi)} } \bigg)^2 \mid \cX\bigg] 
 &\leq \frac{16 p}{N^2} \sum_{\boldj \neq \boldi} \frac{p(1- p) h^2(\bX_\boldj)}{\bar{U}_{h^2}N(1 - p)}  \notag \\
 &=  \frac{16  p}{N^2\bar{U}_{h^2}} \sum_{\boldj \neq \boldi} \frac{ h^2(\bX_\boldj)}{{n \choose m}} \leq \frac{16 U_{h^2}}{N {n \choose m} \bar{U}_{h^2}}  \label{2nd_2nd_moment_bdd_final}
\end{align}
Lastly, for $\bE\Big[ \frac{N U_n^2}{ \bar{U}_{h^2}(1-p)} \Big( \frac{N(\hat{N}_c^{(\boldi)} - \hat{N}_c)}{ \hat{N}_c\hat{N}_c^{(\boldi)} } \Big)^2 \mid \cX \Big]$, we have
\begin{align}
\bE\Bigg[ \frac{N U_n^2}{ \bar{U}_{h^2}(1-p)} \Bigg( \frac{N(\hat{N}_c^{(\boldi)} - \hat{N}_c)}{ \hat{N}_c\hat{N}_c^{(\boldi)} } \Bigg)^2 \mid \cX \Bigg] 
&\leq \frac{16 U_n^2 \bE[ (\hat{N}_c^{(\boldi)} - \hat{N}_c)^2]}{ \bar{U}_{h^2}(1-p)N}  \notag\\
&\leq \frac{16 U_n^2 \bE[ Z_\boldi^2]}{ \bar{U}_{h^2}(1-p)N} = \frac{16 U_n^2 }{ \bar{U}_{h^2}(1-p){n \choose m}},
\label{3rd_2nd_moment_bdd_final}
\end{align}
where the last inequality comes from  \propertyref{censoring_property}$(i)$. Putting \eqref{1st_2nd_moment_bdd_final}-\eqref{3rd_2nd_moment_bdd_final} back into \eqref{D1minusD1i_in_two_parts}, we obtain \eqref{D1minusD1i_norm}.

\subsubsection{Wrapping up}
Since $\bar{D}_1^2$ amounts to the non-negative  $D_1^2$ upper-censored at $1/4$, \propertyref{censoring_property}$(i)$ and $(ii)$ give that $|\bar{D}_1- \bar{D}_1^{(\boldi)}| \leq |D_1- D_1^{(\boldi)}|$ and $\bar{D}_1^2 \leq D_1^2$. These, together with 
 the Markov inequality
\[
 P(|D_1| > 1/2) \leq 2\bE \bigg[\sqrt{\bE[ D_1^2|\cX]}  \bigg],
 \]
 allow us to bound $T(D_1)$ with
\eqref{1st_setp_D1_bdd} and \eqref{D1minusD1i_norm}  as
\begin{align}
T(D_1) 
&\leq \bE \Bigg[ 
 \frac{6 \sqrt{ (1 + 3N)U_{h^2}}}{\sqrt{\bar{U}_{h^2}}N} +   \frac{6\sqrt{2N} \big|U_n  \big|  }{ \sigma_h\sqrt{(1 -p)}}
+ \frac{2 U_{h^2} \sqrt{1 + N}}{\sigma_h N \sqrt{\bar{U}_{h^2}}} 
+ \frac{4 U_{|h|} \sqrt{U_{h^2}}}{\sigma_h\sqrt{N \bar{U}_{h^2}}} + \frac{4 U_{|h|} |U_n| }{\sigma_h\sqrt{\bar{U}_{h^2}(1-p)}}\Bigg]\notag\\
&\leq \bE \Bigg[ \frac{6\sqrt{(1+ 3N)     }}{N  }
+   \frac{6\sqrt{2N} \big|U_n  \big|  }{ \sigma_h\sqrt{(1 -p)}}
+ \frac{2 \sqrt{U_{h^2} (1 + N)}}{\sigma_h N } 
+ \frac{4 \sqrt{U_{h^2}}}{\sigma_h\sqrt{N}} 
+ \frac{4  |U_n| }{\sigma_h\sqrt{(1-p)}}\Bigg] \label{TD1_prelim_bdd}
\end{align}
where we have used  the inequality   $U_{|h|} \leq \sqrt{U_{h^2}} \leq \sqrt{\bar{U}_{h^2}}$ in  \eqref{TD1_prelim_bdd}.   Since $\bE[\sqrt{U_{h^2}}] \leq \sqrt{\bE[U_{h^2}  ]} = \sigma_h$ and 
\begin{equation*} \label{d_2nd_moment_bound}
\bE[|U_n|] \leq \sqrt{\bE[U_n^2]} \leq \sigma_h \sqrt{K_{n, m, d}}, 
\end{equation*}
by the   bound on the variance of a U-statistic of rank $d$ in \citet[Theorem 5.2]{hoeffding1948}, continuing from \eqref{TD1_prelim_bdd} we further get
\begin{equation}\label{TD1_prelim_bddb}
T(D_1) 
\leq 
 \frac{6\sqrt{(1+ 3N)     }}{N  }
+ \frac{6\sqrt{2 N K_{n, m, d}} }{ \sqrt{(1-p)}}
+ \frac{2 \sqrt{(1 + N)}}{ N } 
+ \frac{4 }{\sqrt{N}} 
+ \frac{4   \sqrt{K_{n, m ,d}}}{\sqrt{(1-p)}},
\end{equation}
which implies
\lemref{TD1_bdd}$(i)$.


\subsection{Proof of \lemref{TD1_bdd}$(ii)$} \label{sec:D1_2_norm_pf} Like the proof for \lemref{TD1_bdd}$(i)$, there are three similar steps:
 \subsubsection{Bounding $\sqrt{\bE[D_1^2 | \cX]} 
$}
 First, by the definition of  $D_1$ in \configref{Nasympn}$(ii)$, \eqref{Ncminus_N_bdd}, $N/2\leq \hat{N}_c$ and $\sigma_h^2/2 \leq \bar{U}_{h^2} $
\begin{equation*}
|D_1| \leq  \frac{2 \big|\big(\sum_{\boldj \in I_{n, m}} \xi_\boldj \big) \big(\sum_{\boldk \in I_{n, m}} (Z_\boldk - p) \big)\big|}{N} 
+ \frac{2\sqrt{2} \big|U_n \cdot (\hat{N} - N) \big|  }{ \sigma_h\sqrt{N(1 -p)}}.
\end{equation*}
Hence, by the moment estimate \eqref{moment_est_for_D1_bdd} and \eqref{Nhat_1_norm_bddb},
\begin{equation*} \label{D1_2_norm}
\sqrt{\bE[ D_1^2|\cX]}\leq    \frac{2 \sqrt{ (1 + 3N)U_{h^2}}}{\sqrt{\bar{U}_{h^2}}N} + \frac{2\sqrt{2} \big|U_n  \big|  }{ \sigma_h\sqrt{(1 -p)}}.
\end{equation*}

 \subsubsection{Bounding $\sqrt{\bE[ (D_1  - D_1^{(\boldi)})^2 |\cX]}$}
From the definitions of $D_1$ and $D_1^{(\boldi)}$ in \configref{Nasympn}, write 
\begin{align*}
D_1  - D_1^{(\boldi)} 
&= \xi_\boldi \frac{N - \hat{N}_c}{\hat{N}_c} + \Bigg(\sum_{\boldj \neq \boldi} \xi_\boldj +\frac{ \sqrt{N}U_n}{ \sqrt{\bar{U}_{h^2} (1 - p)}}\Bigg) \Bigg( \frac{N - \hat{N}_c}{\hat{N}_c} -  \frac{N - \hat{N}_c^{(\boldi)}}{\hat{N}_c^{(\boldi)}} \Bigg)\\
&= \xi_\boldi \frac{N - \hat{N}_c}{\hat{N}_c} + \Bigg(\sum_{\boldj \neq \boldi} \xi_\boldj + \frac{ \sqrt{N}U_n}{ \sqrt{\bar{U}_{h^2} (1 - p)}}\Bigg) \Bigg( \frac{N(\hat{N}_c^{(\boldi)} - \hat{N}_c)}{ \hat{N}_c\hat{N}_c^{(\boldi)} } \Bigg).
\end{align*}
Since it has exactly the same form as \eqref{D1minusD1boldi_form} for \configref{Nlln},  using the proof there we get 
\begin{equation*} \label{D1minusD1i_norm_for_2nd_config}
\sqrt{\bE[ (D_1  - D_1^{(\boldi)})^2 |\cX]} \leq \frac{2|h(\bX_\boldi)| \sqrt{1 + N}}{N\sqrt{ \bar{U}_{h^2} {n \choose m} } }  +  \frac{4\sqrt{ U_{h^2}}}{\sqrt{ \bar{U}_{h^2}  N {n \choose m}} }
+ \frac{4 |U_n| }{\sqrt{\bar{U}_{h^2}(1-p){n \choose m}}}.
\end{equation*}
\subsubsection{Wrapping up} Similarly to the proof for \configref{Nlln}, using $|\bar{D}_1- \bar{D}_1^{(\boldi)}| \leq |D_1- D_1^{(\boldi)}|$ and $\bar{D}_1^2 \leq D_1^2$ and $U_{|h|} \leq \sqrt{U_{h^2}} \leq \sqrt{\bar{U}_{h^2}}$, we can first write
\[
T(D_1) 
\leq \bE \Bigg[ \frac{6\sqrt{(1+ 3N)     }}{N  }
+ \frac{6\sqrt{2} |U_n|}{\sigma_h \sqrt{1-p}}
+ \frac{2 \sqrt{U_{h^2} (1 + N)}}{\sigma_h N } 
+ \frac{4 \sqrt{U_{h^2}}}{\sigma_h\sqrt{N}} 
+ \frac{4  |U_n| }{\sigma_h\sqrt{(1-p)}}\Bigg],
\]
which only differs from \eqref{TD1_prelim_bdd} in the second term on the right hand side. Applying similar techniques that lead to \eqref{TD1_prelim_bddb}   but  with the standard bound
\[
\bE[U_n^2] \leq \frac{m\sigma_h^2}{n}
\] 
 for the variance of a  U-statistic \citep[Theorem 5.2]{hoeffding1948}\footnote{This bound is sharp  with respect to $n$ under the non-degeneracy condition \eqref{non_degenerate_condition}.}, we can prove \lemref{TD1_bdd}$(ii)$.

 \section{Proof of \lemref{TD2_bdd}, the bounds on $T(D_2)$} \label{app:bounding_T_D2}

We will first prove \lemref{TD2_bdd}$(iii)$ under the  stronger finite fourth moment assumption $\bE[h^4] < \infty$, which is simple: 

\begin{proof} [Proof of \lemref{TD2_bdd}$(iii)$]

Given  $\bar{D}_2^2 \leq |\bar{D}_2|$, $|\bar{D}_2| \leq |D_2|$  by   \propertyref{censoring_property}$(ii)$ and \lemref{fxfxprime_bdd}$(iv)$, in light of the definition of $T(D_2)$ we have
  \begin{equation*}
T(D_2) 
 \leq 4 \|D_2\|_2 
\leq  \frac{8 \|  U_{h^2} - \sigma_h^2\|_2}{\sigma_h^2} 
\leq \frac{8 \sqrt{m \bE[(h^2 - \sigma_h^2)^2]}}{\sigma_h^2 \sqrt{n}},
\end{equation*}
 where the second last inequality comes from $\bar{U}_{h^2} \geq \sigma_h^2/2$ and the  last inequality uses the standard  bound on U-statistic variance  \citep[Theorem 5.2]{hoeffding1948} .  This has proved the bound for $T(D_2)$ in \lemref{TD2_bdd}$(iii)$ under the stronger finite fourth moment assumption $\bE[h^4] < \infty$. 
 \end{proof}

The bounds for  $T(D_2)$ in \lemref{TD2_bdd}$(i)$ and $(ii)$ respectively for $\frakz_\cX$ in \configsref{Nlln}$(i)$ and  \configssref{Nasympn}$(i)$ under only the  finite third moment assumption $\bE[|h|^3] < \infty$ are much more involved to establish.  To prove them, 
we first upper-censor the $\eta_i$'s in \eqref{eta_def} as
\begin{equation*} \label{etai_upper_censored}
\bareta_i \equiv \eta_i I(\eta_i \leq 1) + 1 I(\eta_i > 1),
\end{equation*}
 from which, with the Hoeffding decomposition remainder term in \eqref{H_decomp_remainder},  we can further define
\begin{equation} \label{Pi1_def}
\Pi_1 \equiv \sum_{i=1}^n ( \bareta_i - \bE[\bareta_i]),
\end{equation}
\begin{align} \label{Pi2_def}
\Pi_2 &\equiv   \cR
- \sum_{i=1}^n \bE[(\eta_i - 1)I(|\eta_i| >1)],
\end{align}
\begin{equation}\label{D3_def}
D_3 \equiv - \frac{\sigma_h^2 (\Pi_1 + \Pi_2)}{\bar{U}_{h^2}}
\end{equation}
and 
\[
\bar{D}_3 \equiv D_3 I\bigg(-\frac{1}{2}\leq D_3 \leq \frac{1}{2}\bigg)- \frac{1}{2} I\bigg(D_3  < -\frac{1}{2}\bigg)+ \frac{1}{2 }I\bigg(D_3 >\frac{1}{2}\bigg).
\]
 
 The following moment bounds on \eqref{Pi1_def} and \eqref{Pi2_def} will be useful in the sequel.
 \begin{lemma} [Moment estimates on $\Pi_1$ and $\Pi_2$] \label{lem:Rosenthal_general}
If $\bE[|h|^3] < \infty$, the following moment bounds hold for $\Pi_1$ and  $\Pi_2$ defined in  \eqref{Pi1_def} and  \eqref{Pi2_def}:
\[
 \bE\big[ |\Pi_1|^\ell\big] \leq \begin{cases}
  C(\ell) \bigg\{  \Big(\frac{m^{3/2}\bE[\Psi_1^{3/2}]}{ n^{1/2}\sigma_h^3}\Big)^{\ell/2}
+ \frac{m^{3/2}\bE[\Psi_1^{3/2}]}{ n^{1/2}\sigma_h^3} \bigg\}  & \text{ for } \;\ 2 < \ell  < \infty\\
 C\frac{m^{3/2}\bE[\Psi_1^{3/2}]}{ n^{1/2}\sigma_h^3} & \text{ for } \;\ \frac{3}{2}\leq  \ell \leq2
\end{cases},
\]
where $C(\ell) > 0$ is a constant depending only on $\ell$, 
and 
\item

\[
\|\Pi_2\|_{3/2} \leq 
\frac{m^{3/2}\bE[\Psi_1^{3/2}]}{ n^{1/2}\sigma_h^3} +  \|\cR\|_{3/2}.
\]
\end{lemma}
\begin{proof}[Proof of \lemref{Rosenthal_general}]
We first prove the bound for $\Pi_1$, and note that the constant $C(\ell)$ below, while depending only on $\ell$, can change in values at different occurrences. 
For $\ell > 2$, by 
\citet[Theorem 3]{rosenthal1970subspaces}'s inequality, we have
\begin{align*}
&\bE\Big[ \Big|\sum_{i=1}^n ( \bareta_i - \bE[\bareta_i]) \Big|^\ell  \Big] \\
&\leq C(\ell)\Bigg\{ \Big( \sum_{i=1}^n \bE[  ( \bareta_i - \bE[\bareta_i])^2]\Big)^{\ell/2} + \sum_{i=1}^n \bE[|\bareta_i - \bE[\bareta_i]|^\ell ] \Bigg\} \\
&\leq   C(\ell)\Bigg\{ \Big( \sum_{i=1}^n \bE[  | \bareta_i - \bE[\bareta_i]|^{3/2}]\Big)^{\ell/2} + \sum_{i=1}^n \bE[| \bareta_i - \bE[\bareta_i]|^{3/2}  ] \Bigg\} \text{ since }  |\bar{\eta}_i - \bE[\bar{\eta}_i ]| \leq 1 \\
&\leq  C(\ell)\Bigg\{ \Big( \sum_{i=1}^n \bE[  \bareta_i^{3/2}]\Big)^{\ell/2} + \sum_{i=1}^n \bE[ \bareta_i^{3/2}  ] \Bigg\}  \text{ since } \bE[\bareta_i^{3/2} ] \geq (\bE[\bareta_i] )^{3/2}> 0\\
&\leq C(\ell) \bigg\{  \bigg(\frac{m^{3/2}\bE[\Psi_1^{3/2}]}{ n^{1/2}\sigma_h^3}\bigg)^{\ell/2}
+ \frac{m^{3/2}\bE[\Psi_1^{3/2}]}{ n^{1/2}\sigma_h^3} \bigg\}.
\end{align*}
Moreover,  by 
\citet[Lemma 1]{chatterji1969p},
we also have, for $\ell \in [\frac{3}{2}, 2]$,
\begin{align*}
\bE\bigg[ \Big|\sum_{i=1}^n ( \bareta_i - \bE[\bareta_i]) \Big|^\ell  \bigg] 
 &\leq 
  2\sum_{i =1}^n  \bE[  | \bareta_i - \bE[\bareta_i]|^\ell]\\
   &\leq 
2\sum_{i =1}^n  \bE[  ( \bareta_i + \bE[\bareta_i])^\ell] \text{ since }   \bar{\eta}_i \geq 0 \\
  &\leq    C\frac{m^{3/2}\bE[\Psi_1^{3/2}]}{ n^{1/2}\sigma_h^3} \text{ since } \bareta_i^\ell \leq \bareta_i^{3/2}.
\end{align*}
Hence we have proved the bound for $\Pi_1$.

Now we prove the bound for $\Pi_2$. First, 
\begin{equation} \label{remnant_expectation_bdd}
\sum_{i=1}^n \bE[(\eta_i - 1)I(|\eta_i| >1)] \leq \sum_{i=1}^n  \bE[\eta_i I(\eta_i >1)] \leq \sum_{i=1}^n  \bE[\eta_i^{3/2} ] = \frac{m^{3/2}\bE[\Psi_1^{3/2}]}{ n^{1/2}\sigma_h^3}.
\end{equation}
Next, by a standard moment inequality for U-statistics \citep[Theorem 2.1.3]{KB1994Ustat}, we have
\begin{multline*}
\bE\Bigg[\bigg| \sum_{r = 2}^m {m \choose r} {n \choose r}^{-1} \sum_{1 \leq i_1 < \dots < i_r \leq n}\frac{\pi_r(h^2) (X_{i_1}, \dots, X_{i_r})}{\sigma_h^2}\bigg|^{3/2}\Bigg] \\
\leq  (m - 1)^{1/2} \sum_{r = 2}^m  {m \choose r}^{3/2} {n \choose r}^{-1/2}  \frac{2^{r/2}\bE[|\pi_r(h^2)|^{3/2}]}{\sigma_h^{3}}.
\end{multline*}
Collecting these inequalities gives the bound for $\|\Pi_2\|_{3/2}$ in  \lemref{Rosenthal_general}.

\end{proof}

Since 
\[
m= \sum_{i=1}^n \bE[\eta_i] = \sum_{i=1}^n \bE[\bareta_i]+ \sum_{i=1}^n \bE[(\eta_i - 1)I(|\eta_i| >1)],
\]
in light of  the Hoeffding decomposition in \eqref{H_decomp_h2_ustat}, one can see that
\begin{equation} \label{D2_equal_D3}
D_2 = D_3 \text{ when  } \max_{1 \leq i \leq n} \eta_i \leq 1,
\end{equation}
by comparing the definitions of $D_2$ and $D_3$ in \eqref{D2_def} and \eqref{D3_def}. Now, because 
 $T(D_2)$'s component terms satisfy the bounds $P(|D_2|>1/2) \leq 1$, $|\bar{D}_2^2 - \bar{D}_3^2|  \leq 1/2$ and  $| \frakz_\cX (\bar{D}_2 - \bar{D}_3) f_{\frakz_\cX} (\sum_{\boldi \in I_{n, m}} \barxi_\boldi )| \leq 1$ by \lemref{fxfxprime_bdd}$(iv)$, 
 it is easy to see from   \eqref{D2_equal_D3} that
  \begin{equation} \label{TD2_prelim_bdd_D3}
T(D_2) 
 \leq P( |D_3|>1/2)+ \bE[\bar{D}_3^2 ] +  \Big| \bE\Big[\frakz_\cX \bar{D}_3  f_{\frakz_\cX} \Big(\sum_{\boldi \in I_{n, m}} \barxi_\boldi \Big) \Big]\Big| + \frac{5}{2} P(\Omega \backslash  \cE_{2, \cX}),
\end{equation}
where we have defined the event
  \[
 \cE_{2, \cX} \equiv \bigg\{ \max_{1 \leq i \leq n}\eta_i \leq 1\bigg\}.
 \]
In \appsref{Proof_of_tail_prob_bdd_for_D2}-\appssref{proof_of_x_exp_barD2_f_indicator_bdd1},  we will derive that one can bound each of the first three terms on the right hand side of \eqref{TD2_prelim_bdd_D3}  in  moments of $\Pi_1$ and $\Pi_2$ as:
\begin{equation}\label{tail_prob_bdd_for_D2}
P(|D_3|>1/2)  \leq C(\bE[\Pi_1^2] +  \|\Pi_2 \|_{3/2})
\end{equation}
\begin{equation} \label{something1}
\bE[\bar{D}_3^2 ] \leq  C( \bE[\Pi_1^2] + \|\Pi_2 \|_{3/2} ) 
\end{equation}
\begin{multline} \label{x_exp_barD2_f_indicator_bdd_final}
\Big| \bE\Big[ \frakz_\cX \bar{D}_3 f_{\frakz_\cZ} \Big(\sum_{\boldi \in I_{n, m}} \barxi_\boldi \Big)  \Big]\Big| 
 \leq 
 \Big| \bE\Big[  \frakz_\cX \Pi_1  f_{ \frakz_\cX} \Big(\sum_{\boldj \in I_{n, m}} \xi_\boldj \Big)  \Big] \Big| +   \frac{2 \sqrt{2}\bE[|h|^3] (1 - 2p + 2p^2)}{\sigma_h^3 (N(1-p))^{1/2}}
\\
+ C\bigg\{ \frac{m^{3/2} \bE[\Psi_1^{3/2}]}{\sigma_h^3 n^{1/2}} \cdot \big(1 +   \|\Pi_1\|_{3/2}  \big) + \bE[|\Pi_1|^{3/2}] + \bE[\Pi_1^2] + \bE[|\Pi_1|^3] +  \|\Pi_2\|_{3/2} \bigg\} + P(\Omega \backslash \cE_{1, \cX}).
\end{multline}
Hence, upon applying the moment bound for $\bE[|\Pi_1|^{3/2}]$, $\bE[\Pi_1^2]$ and $\bE[|\Pi_1|^3]$ in  \lemref{Rosenthal_general} to the right hand sides of \eqref{tail_prob_bdd_for_D2}-\eqref{x_exp_barD2_f_indicator_bdd_final}, one can further get
\begin{equation}\label{tail_prob_bdd_for_D2_modified}
P(|D_3|>1/2)  \leq C \bigg(\frac{m^{3/2}\bE[\Psi_1^{3/2}]}{ n^{1/2}\sigma_h^3} +  \|\Pi_2 \|_{3/2}\bigg)
\end{equation}
\begin{equation} \label{something1_modified}
\bE[\bar{D}_3^2 ] \leq  C \bigg(\frac{m^{3/2}\bE[\Psi_1^{3/2}]}{ n^{1/2}\sigma_h^3} + \|\Pi_2 \|_{3/2}  \bigg)
\end{equation}
\begin{multline} \label{x_exp_barD2_f_indicator_bdd_final_modified}
\Big| \bE\big[  \frakz_\cX\bar{D}_3 f_{ \frakz_\cX} \big(\sum_{\boldi \in I_{n, m}} \barxi_\boldi \big)  \big]\Big| 
 \leq 
 \Big| \bE\big[   \frakz_\cX\Pi_1  f_{ \frakz_\cX} \big(\sum_{\boldj \in I_{n, m}} \xi_\boldj \big)  \big] \Big| 
 \\
  + C \bigg(  \frac{\bE[|h|^3] (1 - 2p + 2p^2)}{\sigma_h^3 (N(1-p))^{1/2}}
+  \frac{m^{3/2}\bE[\Psi_1^{3/2}]}{ n^{1/2}\sigma_h^3}
+ \|\Pi_2\|_{3/2}\bigg) +   P(\Omega \backslash \cE_{1, \cX} )  ,
\end{multline}
where we have absorbed the terms $\big(\frac{m^{3/2}\bE[\Psi_1^{3/2}]}{ n^{1/2}\sigma_h^3} \big)^{3/2}$ and $\big(\frac{m^{3/2}\bE[\Psi_1^{3/2}]}{ n^{1/2}\sigma_h^3} \big)^{5/3}$  showing up in the bounds for $\bE[|\Pi_1|^3]$ and $\frac{m^{3/2} \bE[\Psi_1^{3/2}]}{\sigma_h^3 n^{1/2}}  \|\Pi_1\|_{3/2}$ into $\frac{m^{3/2}\bE[\Psi_1^{3/2}]}{ n^{1/2}\sigma_h^3}$. The latter are permissible,  because 
\[
| \bE[  \frakz_\cX\bar{D}_3 f_{ \frakz_\cX} (\sum_{\boldi \in I_{n, m}} \barxi_\boldi )  ]| \leq 1/2
\]
 by virtue of \lemref{fxfxprime_bdd}$(iv)$ and $\bar{D}_3$ being censored  within the interval $[-1/2, 1/2]$; for the bound on  $| \bE[  \frakz_\cX\bar{D}_3 f_{ \frakz_\cX} (\sum_{\boldi \in I_{n, m}} \barxi_\boldi )  ]|$ to be non-vacuous, one can hence without loss of generality assume  $\frac{m^{3/2}\bE[\Psi_1^{3/2}]}{ n^{1/2}\sigma_h^3} \leq 1$, which implies 
 \[
 \bigg(\frac{m^{3/2}\bE[\Psi_1^{3/2}]}{ n^{1/2}\sigma_h^3}\bigg)^{5/3} \vee  \bigg(\frac{m^{3/2}\bE[\Psi_1^{3/2}]}{ n^{1/2}\sigma_h^3}\bigg)^{3/2} \leq \frac{m^{3/2}\bE[\Psi_1^{3/2}]}{ n^{1/2}\sigma_h^3}
 \]
  to justify the absorption. By recalling \eqref{TD2_prelim_bdd_D3}, \eqref{tail_prob_bdd_for_D2_modified}-\eqref{x_exp_barD2_f_indicator_bdd_final_modified} give
\begin{multline}\label{TD2_premature_bdd}
T(D_2) \leq  \Big| \bE\big[  \frakz_\cX  \Pi_1  f_{ \frakz_\cX} \big(\sum_{\boldj \in I_{n, m}} \xi_\boldj \big)  \big] \Big|  \\
+ C\bigg(\frac{\bE[|h|^3] (1 - 2p + 2p^2)}{\sigma_h^3 (N(1-p))^{1/2}} + \frac{m^{3/2}\bE[\Psi_1^{3/2}]}{ n^{1/2}\sigma_h^3} + \|\Pi_2\|_{3/2}\bigg) +   P(\Omega \backslash \cE_{1, \cX} ) + \frac{5}{2} P(\Omega \backslash  \cE_{2, \cX}).
\end{multline}
For  the rest of this section, we will show the estimate 
\begin{multline} \label{paranthesis_2}
 \Big| \bE\big[  \frakz_\cX \Pi_1  f_{ \frakz_\cX} \big(\sum_{\boldj \in I_{n, m}} \xi_\boldj \big)  \big] \Big|\\
  \leq 
 \begin{cases}
 C   \frac{m^{3/2}\bE[\Psi_1^{3/2}]}{ n^{1/2}\sigma_h^3} \text{ if $\frakz_\cX$ is as in \configsref{Nlln}$(i)$}  \\
C \bigg(1 + \frac{\sqrt{mN}}{\sqrt{n(1-p)}}\bigg)  \frac{m^{3/2}\bE[\Psi_1^{3/2}]}{ n^{1/2}\sigma_h^3} \text{if $\frakz_\cX$ is as in \configsref{Nasympn}$(i)$}  
 \end{cases}.
\end{multline}
 Combining \eqref{TD2_premature_bdd}, \eqref{paranthesis_2}, the union bound
\begin{equation*}
P(\Omega \backslash \cE_{2, \cX}) \leq \sum_{i=1}^n P( \eta_i> 1)\leq \sum_{i=1}^n \bE[\eta_i^{3/2}] = \frac{m^{3/2} \bE[\Psi_1^{3/2}]}{n^{1/2}\sigma_h^3}
\end{equation*}
 and the bound for $\|\Pi_2\|_{3/2}$ in \lemref{Rosenthal_general} will then have proven \lemref{TD2_bdd}$(i)$ and $(ii)$ under the  finite third moment assumption $\bE[|h|^3] < \infty$.

For the proof of  \eqref{paranthesis_2} to follow, we will use the symbol  ``$\sum_{\boldj \in I_{n,m}, i \not\in \boldj}$'' and ``$\sum_{\boldj \in I_{n,m}, i \in \boldj}$'' introduced at the beginning of \secref{pf_Nggn_case}. 
Let 
\begin{equation} \label{Uh2_leave_one_out}
U_{h^2}^{(i)} \equiv {n \choose m}^{-1} \sum_{\boldj \in I_{n,m}, i \not\in \boldj} h^2 (X_\boldj),
\end{equation}
be the ``leave-one-out'' version of $U_{h^2}$ eliminating terms involving $X_i$,  and define
\begin{multline} \label{Uh2_leave_one_out_bar}
\bar{U}^{(i)}_{h^2} \equiv  \frac{\sigma_h^2}{2}I \bigg(U^{(i)}_{h^2} \leq \frac{\sigma_h^2}{2}\bigg) + U_{h^2}^{(i)} I\bigg(U_{h^2}^{(i)} > \frac{\sigma_h^2}{2}\bigg) \quad \text{ and }\quad
 {_i}\xi_{\boldj} =  \frac{ (Z_\boldj - p)h(\bX_\boldj)}{\sqrt{\bar{U}_{h^2}^{(i)}N(1 - p)}},
\end{multline}
analogously to $\bar{U}_{h^2}$ and $\xi_\boldj$. Moreover, 
 let 
 \begin{equation} \label{leave_one_out_frakz}
 \frakz_{\cX^{(i)} }  \equiv 
   \begin{cases} 
   z& \text{if $\frakz_\cX$ is as in \configsref{Nlln}$(i)$} \\
    \frac{z \sigma}{\sigma_h\sqrt{ \alpha(1 -p)}} -   \frac{\sqrt{N}U_n^{(i)}}{\sigma_h \sqrt{1 - p}}     &  \text{if $\frakz_\cX$ is as in \configsref{Nasympn}$(i)$}  \end{cases},
 \end{equation}
 where
  \[
 U_n^{(i)} \equiv {n \choose m}^{-1} \sum_{\boldj \in I_{n,m}, i \not\in \boldj} h (X_\boldj).
 \]
 is a leave-one-out version of $U_n$ for any given $i \in [n]$.
  Then write
\begin{align}
\bE\Big[  \frakz_\cX \Pi_1     f_{\frakz_\cX } \Big(\sum_{\boldj \in I_{n, m}} \xi_\boldj \Big)  \Big]
 = 
  &\sum_{i=1}^n\bE\Big[  ( \bareta_i - \bE[\bareta_i]) \cdot \Big( \frakz_\cX  f_{\frakz_\cX } \big(\sum_{\boldj \in I_{n, m}} \xi_\boldj \big) -\frakz_{\cX^{(i)}}   f_{\frakz_{\cX^{(i)} }} \big(\sum_{\boldj \in I_{n, m}} \xi_\boldj \big) \Big) \Big] \notag \\
 + & \sum_{i=1}^n\bE\Big[  \frakz_{\cX^{(i)} } \cdot ( \bareta_i - \bE[\bareta_i]) \cdot \Big( f_{\frakz_{\cX^{(i)} } } \big(\sum_{\boldj \in I_{n, m}} \xi_\boldj \big) - f_{ \frakz_{\cX^{(i)} } } \big(\sum_{\boldj \in I_{n, m}} {_i}\xi_\boldj \big) \Big) \Big] 
\notag\\
+  &\sum_{i=1}^n\bE\Big[  \frakz_{\cX^{(i)} } \cdot ( \bareta_i - \bE[\bareta_i]) \cdot f_{ \frakz_{\cX^{(i)} } } \Big(\sum_{\boldj \in I_{n, m}} {_i}\xi_\boldj \Big) \Big]. \label{cool_argument_prep}
\end{align}
We will next establish the bounds
\begin{multline} \label{coolest_ineq}
\Big|\sum_{i=1}^n\bE\Big[  ( \bareta_i - \bE[\bareta_i]) \cdot \Big( \frakz_\cX  f_{\frakz_\cX } \big(\sum_{\boldj \in I_{n, m}} \xi_\boldj \big) -\frakz_{\cX^{(i)}}   f_{\frakz_{\cX^{(i)} }} \big(\sum_{\boldj \in I_{n, m}} \xi_\boldj \big) \Big) \Big] \Big| \\
  \leq  \begin{cases} 
0 &  \text{ if }\quad \frakz_{\cX^{(i)}} = z \text{ as in \configsref{Nlln}$(i)$ }\\
 \frac{14 m^2 \sqrt{N} \bE[ \Psi_1^{3/2}]}{\sigma_h^3 n  \sqrt{1 - p}}    & \text{ if } \quad  \frakz_{\cX^{(i)}} = \frac{z \sigma}{\sigma_h\sqrt{ \alpha(1 -p)}} -   \frac{\sqrt{N}U_n^{(i)}}{\sigma_h \sqrt{1 - p}}  \text{ as in \configsref{Nasympn}$(i)$ }
     \end{cases} ,
\end{multline}
\begin{equation}\label{cool_argument_prep_1}
\Big| \sum_{i=1}^n\bE\Big[ \frakz_{\cX^{(i)} } \cdot ( \bareta_i - \bE[\bareta_i]) \cdot \Big( f_{\frakz_{\cX^{(i)} } } \big(\sum_{\boldj \in I_{n, m}} \xi_\boldj \big) - f_{ \frakz_{\cX^{(i)} }} \big(\sum_{\boldj \in I_{n, m}} {_i}\xi_\boldj \big) \Big) \Big]\Big|
 \leq \frac{(8 + 8\sqrt{2})m^{3/2} \bE[\Psi_1^{3/2}]}{n^{1/2} \sigma_h^3}
\end{equation}
and
\begin{equation}\label{cool_argument_prep_2}
\Big| \sum_{i=1}^n\bE\Big[ \frakz_{\cX^{(i)} } \cdot  ( \bareta_i - \bE[\bareta_i]) \cdot  f_{\frakz_{\cX^{(i)} } } \Big(\sum_{\boldj \in I_{n, m}} {_i}\xi_\boldj \Big) \Big]  \Big|\leq  \frac{  (12\sqrt{2} +4)  m^{3/2} \bE[\Psi_1^{3/2}]}{n^{1/2} \sigma_h^3};
\end{equation}
we remark that \eqref{cool_argument_prep_1} and \eqref{cool_argument_prep_2} are true regardless of the configurations of $\frakz_{\cX^{(i)}}$ in \eqref{leave_one_out_frakz}.
Combining \eqref{cool_argument_prep}-\eqref{cool_argument_prep_2} will then finish proving \eqref{paranthesis_2}. 

\subsection{Bound on the term in \eqref{coolest_ineq}}  \label{app:subsectionD1}
We note that the bound is completely trivial if $\frakz_{\cX^{(i)}} = z$ as in the first configuration in  \eqref{leave_one_out_frakz}, so for the rest of this \appref{subsectionD1} we focus on proving the bound for $\frakz_{\cX^{(i)}} =   \frac{z \sigma}{\sigma_h\sqrt{ \alpha(1 -p)}} -   \frac{\sqrt{N}U_n^{(i)}}{\sigma_h \sqrt{1 - p}}$.

For the Stein equation solution in \eqref{Stein_solution}, 
let 
$
s(z; w) \equiv \frac{\partial{(zf_z(w))}}{\partial z}
$
be the derivative of $zf_z(w)$ with respect to $z$, treated as a function in $z$. Since $zf_z(w)$ is not differentiable at $z = w$, we simply define $s(w; w) \equiv 0$. 
By the fundamental theorem of calculus, one can write
\[
 \frakz_\cX  f_{\frakz_\cX } \bigg(\sum_{\boldj \in I_{n, m}} \xi_\boldj \bigg) -\frakz_{\cX^{(i)}}   f_{\frakz_{\cX^{(i)} }} \bigg(\sum_{\boldj \in I_{n, m}} \xi_\boldj \bigg) =  \int^{\frakz_\cX}_{\frakz_{\cX^{(i)}}}  s\Big(z; \sum_{\boldj \in I_{n, m}} \xi_\boldj  \Big)dz,
\]
and hence, with Fubini's theorem,
\begin{multline} \label{after_FTC}
\bigg|\bE\Big[  ( \bareta_i - \bE[\bareta_i]) \cdot \Big( \frakz_\cX  f_{\frakz_\cX } \big(\sum_{\boldj \in I_{n, m}} \xi_\boldj \big) -\frakz_{\cX^{(i)}}   f_{\frakz_{\cX^{(i)} }} \big(\sum_{\boldj \in I_{n, m}} \xi_\boldj \big) \Big) \Big] \bigg|\\
\leq  \bE \Bigg[ |\bareta_i - \bE[\bareta_i]| \cdot \int^{\frakz_\cX \vee  \frakz_{\cX^{(i)}}}_{ \frakz_\cX \wedge \frakz_{\cX^{(i)}}}   \bE \bigg[\Big|s\Big(z; \sum_{\boldj \in I_{n, m}} \xi_\boldj  \Big) \Big|  \mid \cX\bigg]dz \Bigg].
\end{multline}
We will now bound the inner conditional expectation $\bE [|s(z; \sum_{\boldj \in I_{n, m}} \xi_\boldj  ) |  \mid \cX]$ in the preceding display. 
Given the form of $f_z$  in \eqref{Stein_solution},  we have
\begin{equation} \label{form_of_s_function}
s(z; w)
= 
\begin{cases}
\sqrt{2 \pi} e^{w^2/2} \Phi(w)( \barPhi(z) - z  \phi(z))  & \text{ if } w < z \\
 \sqrt{2 \pi} e^{w^2/2}  \bar{\Phi}(w) (\Phi(z) + z \phi(z)) & \text{ if } z < w 
\end{cases}.
\end{equation}
Since  (see \citet[p.16]{chen2010normal},  for instance)
\begin{equation} \label{mills_ratio_bdd}
\Phi(-w)  = \barPhi(w) \leq \min \bigg(\frac{1}{2} ,  \frac{1}{w \sqrt{2 \pi}} \bigg) e^{-w^2/2} 
\text{ for any } w > 0,
\end{equation}
one can see that, if $w <z$,  we have
\begin{equation*}
\Big|\sqrt{2 \pi} e^{ w^2/2} \Phi(w) \barPhi(z)\Big| \leq  \sqrt{\frac{\pi}{2}} \text{ and } 
\Big|\sqrt{2 \pi} e^{ w^2/2} \Phi(w) z \phi(z) \Big| =  |z| e^{(w^2 - z^2)/2} \Phi(w),
\end{equation*}
and if $z < w$, we have
\begin{multline*}
\Big|\sqrt{2 \pi} e^{ w^2/2} \barPhi(w) \Phi(z)\Big| \leq  \sqrt{\frac{\pi}{2}} \text{ and } 
\Big|\sqrt{2 \pi} e^{ w^2/2} \barPhi(w) z \phi(z)\Big|  = |z| e^{(w^2 - z^2)/2} \barPhi(w).
\end{multline*}
In light of \eqref{form_of_s_function}, 
these imply
\begin{multline} \label{after_the_two_if}
\Big|s \Big(z; \sum_{\boldj \in I_{n, m}} \xi_\boldj \Big)\Big| \leq  \sqrt{\frac{\pi}{2}}
  +   |z|   \times \\
e^{\frac{(\sum_{\boldj \in I_{n, m}} \xi_\boldj )^2-z^2}{2}} \Bigg\{
 I\Big(\sum_{\boldj \in I_{n, m}} \xi_\boldj  > z \Big) \cdot\ \barPhi \Big(\sum_{\boldj \in I_{n, m}} \xi_\boldj  \Big
) +
 I\Big(\sum_{\boldj \in I_{n, m}} \xi_\boldj  < z \Big) \cdot\Phi \Big(\sum_{\boldj \in I_{n, m}} \xi_\boldj  \Big)
\Bigg\},
\end{multline}
and we will respectively bound the conditional expectation of 
\begin{multline*}
  e^{\frac{(\sum_{\boldj \in I_{n, m}} \xi_\boldj )^2-z^2}{2}}  \cdot I\Big(\sum_{\boldj \in I_{n, m}} \xi_\boldj  >z \Big) \cdot\barPhi \Big(\sum_{\boldj \in I_{n, m}} \xi_\boldj  \Big)
   \text{ and }\\
   e^{\frac{(\sum_{\boldj \in I_{n, m}} \xi_\boldj )^2-z^2}{2}}  \cdot I\Big(\sum_{\boldj \in I_{n, m}} \xi_\boldj  < z \Big) \cdot\Phi \Big(\sum_{\boldj \in I_{n, m}} \xi_\boldj  \Big)
\end{multline*}
  on the right hand side: Write
\begin{align}
&\bE\bigg[e^{\frac{(\sum_{\boldj \in I_{n, m}} \xi_\boldj )^2-z^2}{2}} \cdot  I\Big(\sum_{\boldj \in I_{n, m}} \xi_\boldj  >z \Big) \cdot \barPhi \Big(\sum_{\boldj \in I_{n, m}} \xi_\boldj \Big) \mid \cX \bigg]  \notag\\
=&\quad \bE\bigg[e^{\frac{(\sum_{\boldj \in I_{n, m}} \xi_\boldj )^2-z^2}{2}} \cdot I\Big(\sum_{\boldj \in I_{n, m}} \xi_\boldj  >z \Big) \cdot  \barPhi \Big(\sum_{\boldj \in I_{n, m}} \xi_\boldj \Big) \cdot I\Big(|\sum_{\boldj \in I_{n, m}} \xi_\boldj | > \frac{|z|}{\sqrt{2}}\Big)  \mid \cX \bigg] \notag\\
& + \bE\bigg[e^{\frac{(\sum_{\boldj \in I_{n, m}} \xi_\boldj )^2-z^2}{2}} \cdot I\Big(\sum_{\boldj \in I_{n, m}} \xi_\boldj  >z \Big) \cdot  \barPhi\Big(\sum_{\boldj \in I_{n, m}} \xi_\boldj \Big) \cdot I\Big(|\sum_{\boldj \in I_{n, m}} \xi_\boldj | \leq \frac{|z|}{\sqrt{2}}\Big)  \mid \cX \bigg] \notag\\
\leq &\quad \bE\bigg[e^{\frac{(\sum_{\boldj \in I_{n, m}} \xi_\boldj )^2-z^2}{2}} \cdot  I\Big(\sum_{\boldj \in I_{n, m}} \xi_\boldj  >z \Big) \cdot  \barPhi \Big(\sum_{\boldj \in I_{n, m}} \xi_\boldj \Big) \cdot \frac{\sqrt{2} |\sum_{\boldj \in I_{n, m} } \xi_\boldj|}{|z|}  \mid \cX \bigg]  + e^{- z^2/4} \notag\\
\leq& \quad   \frac{ \sqrt{2} \bE[|\sum_{\boldj \in I_{n, m}} \xi_\boldj | \mid \cX]}{|z|} +  e^{-z^2/4}  \label{smartsmart}
\end{align}
where the second last inequality uses $I(|\sum_{\boldj \in I_{n, m}} \xi_\boldj | > \frac{|z|}{\sqrt{2}} ) \leq \frac{\sqrt{2}|\sum_{\boldj \in I_{n, m}} \xi_\boldj |}{|z|}$, and the last inequality uses
\begin{multline*}
e^{\frac{(\sum_{\boldj \in I_{n, m}} \xi_\boldj )^2-z^2}{2}} \cdot I\Big(\sum_{\boldj \in I_{n, m}} \xi_\boldj  >z \Big)\cdot \barPhi (\sum_{\boldj \in I_{n, m}} \xi_\boldj  ) \leq \\
I\Big(\sum_{\boldj \in I_{n, m}} \xi_\boldj  >z, \sum_{\boldj \in I_{n, m}} \xi_\boldj  \leq 0 \Big) 
 + I\Big(\sum_{\boldj \in I_{n, m}} \xi_\boldj  >z, \sum_{\boldj \in I_{n, m}} \xi_\boldj  >0 \Big) \cdot \frac{e^{- z^2/2}}{2}  
 \leq 1
\end{multline*}
 as a result of \eqref{mills_ratio_bdd}.  Since $\bE[ (\sum_{\boldj \in I_{n, m}} \xi_\boldj)^2 \mid \cX] \leq 1$ based on the conditional second moment calculation in \eqref{cX_cond_moment_properties_of_xi}, from \eqref{smartsmart} we further get
\begin{equation}  \label{smartsmart_a}
\bE\Big[e^{\frac{(\sum_{\boldj \in I_{n, m}} \xi_\boldj )^2-z^2}{2}} \cdot  I\Big(\sum_{\boldj \in I_{n, m}} \xi_\boldj  >z \Big) \cdot \barPhi \Big(\sum_{\boldj \in I_{n, m}} \xi_\boldj \Big) \mid \cX \Big]   \leq  \frac{ \sqrt{2}}{|z|} +  e^{-z^2/4}.
\end{equation}
Similarly, one can show that
 \begin{equation} \label{smartsmart_b}
 \bE\Big[ e^{\frac{(\sum_{\boldj \in I_{n, m}} \xi_\boldj )^2-z^2}{2}}  \cdot I\Big(\sum_{\boldj \in I_{n, m}} \xi_\boldj  < z \Big) \cdot\Phi \Big(\sum_{\boldj \in I_{n, m}} \xi_\boldj  \Big)  \mid \cX \Big] \leq   \frac{ \sqrt{2} }{|z|} +  e^{-z^2/4} .
 \end{equation}
Combining \eqref{after_the_two_if}, \eqref{smartsmart_a}, \eqref{smartsmart_b} and $\sup_{z \in \bR}|z| e^{-z^2/4} 
 \leq 1$ implies
\begin{equation} \label{smartsmart2}
\bE\Big[\Big|s \Big(z; \sum_{\boldj \in I_{n, m}} \xi_\boldj \Big)\Big| \mid \cX \Big]
\leq 
\sqrt{\frac{\pi}{2}}+
2 (\sqrt{2} + |z|e^{-z^2/4})
< 7.
\end{equation}

Plugging \eqref{smartsmart2} into \eqref{after_FTC} and using
\[
\frakz_\cX - \frakz_{\cX^{(i)}} = \frac{\sqrt{N}}{\sigma_h \sqrt{1 - p}} (U_n - U_n^{(i)})= \frac{m\sqrt{N}}{\sigma_h n\sqrt{1 - p}} \frac{\sum_{\boldj \in \Inm, i \in \boldj} h(X_\boldj)}{{n -1 \choose m -1}},
\] 
we get
\begin{align}
&\bigg|\bE\Big[  ( \bareta_i - \bE[\bareta_i]) \cdot \Big( \frakz_\cX  f_{\frakz_\cX } \big(\sum_{\boldj \in I_{n, m}} \xi_\boldj \big) -\frakz_{\cX^{(i)}}   f_{\frakz_{\cX^{(i)} }} \big(\sum_{\boldj \in I_{n, m}} \xi_\boldj \big) \Big) \Big] \bigg| \notag\\
& \leq   \frac{7m\sqrt{N}}{\sigma_h n\sqrt{1 - p}}\bE \Bigg[ |\bareta_i - \bE[\bareta_i]| \cdot  \bigg| \frac{\sum_{\boldj \in \Inm, i \in \boldj} h(X_\boldj)}{{n -1 \choose m -1}} \bigg| \Bigg] \notag \\
&\leq \frac{7 m \sqrt{N}}{\sigma_h n \sqrt{1 - p}}  \bE\Bigg[ | \bareta_i - \bE[\bareta_i]|  \cdot  \bigg(\frac{\sum_{\boldj \in \Inm, i \in \boldj} h^2(X_\boldj)}{{n -1 \choose m -1}} \bigg)^{1/2}  \Bigg] \notag\\
&= \frac{7m \sqrt{N}}{\sigma_h n \sqrt{1 - p}}  \bE\Bigg[ | \bareta_i - \bE[\bareta_i]|  \cdot  \bE\bigg[\bigg(\frac{\sum_{\boldj \in \Inm, i \in \boldj} h^2(X_\boldj)}{{n -1 \choose m -1}} \bigg)^{1/2}  \mid X_i\bigg]\Bigg]  \notag\\
&\leq \frac{7 m \sqrt{N}}{\sigma_h n \sqrt{1 - p}}  \bE\big[ | \bareta_i - \bE[\bareta_i]|  \cdot  \Psi_1^{1/2}(X_i)\big]  \notag \\
&\leq  \frac{7 m \sqrt{N}}{\sigma_h n \sqrt{1 - p}}  \bE\big[ ( \bareta_i + \bE[\bareta_i])  \cdot  \Psi_1^{1/2}(X_i)\big]  \notag \\
&= \frac{7 m \sqrt{N}}{\sigma_h n \sqrt{1 - p}}  \bigg( \bE[ \bareta_i  \Psi_1^{1/2}(X_i)]+ \bE[\bareta_i] \bE[ \Psi_1^{1/2}] \bigg) \notag \\
&\leq  \frac{7 m \sqrt{N}}{\sigma_h n \sqrt{1 - p}}  \bigg( \frac{m \bE[ \Psi_1^{3/2}]}{n \sigma_h^2} + \frac{m \bE[\Psi_1] \bE[\Psi_1^{1/2}] }{n \sigma_h^2} \bigg) \notag
\end{align}
which gives \eqref{coolest_ineq} by $\bE[\Psi_1] \bE[\Psi_1^{1/2}] \leq \bE[ \Psi_1^{3/2}]$.

\subsection{Bound on the term  in \eqref{cool_argument_prep_1}}
We will treat both configurations in \eqref{leave_one_out_frakz} together. We first have to separately consider $|\frakz_{\cX^{(i)} } | <2$ and $|\frakz_{\cX^{(i)} } | \geq 2$. Since $\sup_{z,w \in \bR} |\frac{\partial f_z}{\partial z}(w)| \leq 1$ by \lemref{fxfxprime_bdd}$(i)$, 
\begin{multline} \label{cool_argument_prep_1prime}
 \bigg|  \bE\bigg[ \frakz_{\cX^{(i)} }  ( \bareta_i - \bE[\bareta_i]) \cdot \Big( f_{\frakz_{\cX^{(i)} } } \Big(\sum_{\boldj \in I_{n, m}} \xi_\boldj \Big) - f_{\frakz_{\cX^{(i)} } } \Big(\sum_{\boldj \in I_{n, m}} {_i}\xi_\boldj\Big) \Big)\cdot I\Big(|\frakz_{\cX^{(i)} } | <2\Big)  \bigg] \bigg|\\
\leq 2 \bE\bigg[| \bareta_i - \bE[\bareta_i]| \cdot\Big|\sum_{\boldj \in \Inm}(\xi_\boldj - {_i}\xi_\boldj) \Big| \cdot I\Big(|\frakz_{\cX^{(i)} } | <2\Big) \bigg]. 
\end{multline}
Since
\[
 I\Big(\Omega \backslash   \Big\{ |\sum_{\boldj \in \Inm} {_i}\xi_\boldj | \leq |\frakz_{\cX^{(i)} } | - 1 \Big\}, \quad  |\frakz_{\cX^{(i)} } | \geq 2 \Big)  
 \leq  \frac{ | \sum_{\boldj \in \Inm} {_i}\xi_\boldj | }{ |\frakz_{\cX^{(i)} } | - 1} \cdot I( |\frakz_{\cX^{(i)} } | \geq 2),
\]
we can form the bound
\begin{align}
&\bigg|  \bE\bigg[ \frakz_{\cX^{(i)} }  ( \bareta_i - \bE[\bareta_i])  \Big( f_{ \frakz_{\cX^{(i)} } } \big(\sum_{\boldj \in I_{n, m}} \xi_\boldj \big) - f_{ \frakz_{\cX^{(i)} } } \big(\sum_{\boldj \in I_{n, m}} {_i}\xi_\boldj \big) \Big) I(|\frakz_{\cX^{(i)} } | \geq 2)\bigg] \bigg|   \notag\\
&\leq   \bigg|\bE\bigg[ \frakz_{\cX^{(i)} } ( \bareta_i - \bE[\bareta_i]) \Big( f_{\frakz_{\cX^{(i)} } } \big(\sum_{\boldj \in I_{n, m}} \xi_\boldj \big) - f_{ \frakz_{\cX^{(i)} }  } \big(\sum_{\boldj \in I_{n, m}} {_i}\xi_\boldj \big) \Big)  I\Big( \big|\sum_{\boldj \in \Inm} {_i}\xi_\boldj \big| \leq |\frakz_{\cX^{(i)} } | - 1, |\frakz_{\cX^{(i)} } | \geq 2\Big)  \bigg] \bigg|  \notag\\
  &+   \bigg| \bE\bigg[ \frakz_{\cX^{(i)} }  \Big( \bareta_i - \bE[\bareta_i]\Big) \cdot \Big( f_{\frakz_{\cX^{(i)} } } \Big(\sum_{\boldj \in I_{n, m}} \xi_\boldj \Big) - f_{\frakz_{\cX^{(i)} } } \Big(\sum_{\boldj \in I_{n, m}} {_i}\xi_\boldj \Big) \Big) \cdot \frac{ | \sum_{\boldj \in \Inm} {_i}\xi_\boldj | }{|\frakz_{\cX^{(i)} } | - 1}  \cdot I( |\frakz_{\cX^{(i)} } | \geq 2) \bigg] \bigg| \notag\\
  &\leq   \bE\bigg[ |\frakz_{\cX^{(i)} } | e^{1/2- |\frakz_{\cX^{(i)} } |} \cdot  \Big| \bareta_i - \bE[\bareta_i]\Big| \cdot \Big|\sum_{\boldj \in \Inm}(\xi_\boldj - {_i}\xi_\boldj) \Big| \cdot I(|\frakz_{\cX^{(i)} } | \geq 2)  \bigg]  \notag\\
&  \qquad \qquad + \bE\bigg[  \frac{|\frakz_{\cX^{(i)} } |}{|\frakz_{\cX^{(i)} } |-1} \cdot \Big| \bareta_i - \bE[\bareta_i]\Big| \cdot \Big|\sum_{\boldj \in \Inm}(\xi_\boldj - {_i}\xi_\boldj) \Big| \cdot  \Big| \sum_{\boldj \in \Inm} {_i}\xi_\boldj \Big|     \cdot I( |\frakz_{\cX^{(i)} } | \geq 2) \bigg] \label{first_and_second_pieces}
\end{align}   
where  the first term in \eqref{first_and_second_pieces} comes from facts:
\begin{enumerate}[(i)]
\item For any $z \in \bR$ with $|z| \geq 2$, $|\frac{\partial f_z}{\partial w} (w)| \leq e^{1/2 - |z|}$ if $|w| \leq |z|- 1$   by the bounds for $\frac{\partial f_z}{\partial w}$ in   \lemref{fxfxprime_bdd}$(ii)$ and $(iii)$;
\item For any any $z \in \bR$ with $|z| \geq 2$, $
 |\sum_{\boldj \in \Inm} {_i}\xi_\boldj|  \leq |z| - 1$ implies  $|\sum_{\boldj \in \Inm} \xi_\boldj | \leq |z| - 1$ because $\bar{U}^{(i)}_{h^2} \leq \bar{U}_{h^2}$ by definition,
\end{enumerate}
and the second term in \eqref{first_and_second_pieces} comes from  $ \sup_{w, z \in \bR}|\frac{\partial f_z}{\partial w}(w)|\leq 1$ by  \lemref{fxfxprime_bdd}$(i)$.

As 
\begin{equation} \label{some_sup_bdds}
\sup_{|z| \geq 2} |z| e^{\frac{1}{2}-|z|}= 2 e^{-3/2} \text{ and } \sup_{|z| \geq 2} \frac{|z|}{|z|-1}= 2,
\end{equation}
 \eqref{cool_argument_prep_1prime} and \eqref{first_and_second_pieces} imply   
\begin{multline} \label{first_and_second_pieces_pieces}
 \bigg|\bE\bigg[ \frakz_{\cX^{(i)} } \Big( \bareta_i - \bE[\bareta_i]\Big) \cdot \Big( f_{\frakz_{\cX^{(i)} } } \big(\sum_{\boldj \in I_{n, m}} \xi_\boldj \big) - f_{\frakz_{\cX^{(i)} } } \big(\sum_{\boldj \in I_{n, m}} {_i}\xi_\boldj \big) \Big) \bigg] \bigg| \leq\\
 2 \bigg\{\underbrace{\bE\Big[|\bareta_i - \bE[\bareta_i]| \cdot\Big|\sum_{\boldj \in \Inm}(\xi_\boldj - {_i}\xi_\boldj) \Big| \Big]}_{\equiv (1
 )} 
 + 
 \underbrace{\bE\Big[  |\bareta_i - \bE[\bareta_i]|  \cdot \Big|\sum_{\boldj \in \Inm}(\xi_\boldj - {_i}\xi_\boldj) \Big| \cdot  | \sum_{\boldj \in \Inm} {_i}\xi_\boldj |   \Big]}_{\equiv (2)}\bigg\}.
\end{multline}
Now we  estimate the  term $(1)$ in \eqref{first_and_second_pieces_pieces} as
\begin{align}
(1)&=  \bE\bigg[  |\bareta_i - \bE[\bareta_i]| \cdot \bigg| \frac{\sqrt{\bar{U}_{h^2}} - \sqrt{\bar{U}_{h^2}^{(i)}}}{\sqrt{\bar{U}_{h^2}^{(i)}   }}\bigg| \cdot \big|\sum_{\boldj \in \Inm} \xi_\boldj \big|  \bigg] \notag\\
&\leq \frac{\sqrt{2}}{\sigma_h}\bE\bigg[ |\bareta_i - \bE[\bareta_i]| \cdot \sqrt{\bar{U}_{h^2}- \bar{U}_{h^2}^{(i)}} \cdot 
\Big(\underbrace{\bE\Big[ \big(\sum_{\boldj \in \Inm} \xi_\boldj\big)^2 \mid \cX\Big]}_{ = U_{h^2}/\bar{U}_{h^2}\leq 1}\Big)^{1/2}   \bigg] \notag\\
&\leq \frac{\sqrt{2}}{\sigma_h}   \bE\bigg[ |\bareta_i - \bE[\bareta_i]|\cdot \bigg( {n \choose m}^{-1}\sum_{\boldj \in \Inm, i \in \boldj}h^2(\bX_\boldj )   \bigg)^{1/2} \bigg] 
\text{ by \propertyref{censoring_property}$(i)$ of censoring} \notag\\
&\leq  \frac{\sqrt{2}}{\sigma_h} \bE\bigg[ |\bareta_i - \bE[\bareta_i]| \cdot \bigg({n \choose m}^{-1}\sum_{\boldj \in \Inm, i \in \boldj} \bE[h^2(\bX_\boldj )  \mid X_i] \bigg)^{1/2} \bigg] \notag\\
&= \frac{\sqrt{2m} \bE[|\bareta_i - \bE[\bareta_i]|  \Psi_1^{1/2}(X_i)]}{\sigma_h\sqrt{n}} \leq 
\frac{\sqrt{2m} \bE[(\bareta_i + \bE[\bareta_i])  \Psi_1^{1/2}(X_i)]}{\sigma_h\sqrt{n}} \leq \frac{2 \sqrt{2}m^{3/2} \bE[\Psi_1^{3/2}]}{n^{3/2}\sigma_h^3 }
\label{(3)_bdd}.
\end{align}
To estimate the term $(2)$ in \eqref{first_and_second_pieces_pieces}, we first write
\begin{align}
(2) &=  \bE\bigg[  |\bareta_i - \bE[\bareta_i]| \cdot \Bigg|  \frac{\sqrt{\bar{U}_{h^2}} - \sqrt{\bar{U}_{h^2}^{(i)}}}{\sqrt{\bar{U}_{h^2}^{(i)}}}\Bigg|  \cdot  \sqrt{\frac{\bar{U}_{h^2}}{\bar{U}_{h^2}^{(i)}}}\cdot
\underbrace{\bE\bigg[\bigg(\sum_{\boldj \in \Inm}\frac{ (Z_\boldj - p) h(\bX_\boldj)}{\sqrt{\bar{U}_{h^2}N(1- p)}} \bigg)^2 \mid \cX \bigg] }_{= U_{h^2}/\bar{U}_{h^2}\leq  1}
 \bigg] \notag\\
&\leq  \frac{\sqrt{2}}{\sigma_h} \bE\bigg[ |\bareta_i - \bE[\bareta_i]|  \cdot \sqrt{ ( \bar{U}_{h^2}- \bar{U}_{h^2}^{(i)} ) \cdot \bigg(1 + \frac{\bar{U}_{h^2} - \bar{U}_{h^2}^{(i)}}{\bar{U}_{h^2}^{(i)}}\bigg)}   \bigg]   \text{ by } \sqrt{\bar{U}_{h^2}} - \sqrt{\bar{U}_{h^2}^{(i)}} \leq  \sqrt{\bar{U}_{h^2} -  \bar{U}_{h^2}^{(i)}} \notag\\
&\leq  \frac{\sqrt{2}}{\sigma_h} \bE\bigg[ |\bareta_i - \bE[\bareta_i]|  \cdot \bigg(\sqrt{  \bar{U}_{h^2}- \bar{U}_{h^2}^{(i)} } + \frac{\bar{U}_{h^2} - \bar{U}_{h^2}^{(i)}}{\sqrt{\bar{U}_{h^2}^{(i)}}} \bigg) \bigg], \label{intermediate_for_paranthesis2}
\end{align}
where the last inequality uses the concavity of the $\sqrt{\cdot}$ function.
Since $U_{h^2} - U_{h^2}^{(i)} = {n \choose m}^{-1} \sum_{\boldj \in \Inm, i \in \boldj} h^2(\bX_\boldj)$, continuing from \eqref{intermediate_for_paranthesis2}, we have
\begin{align}
(2)&\leq \frac{\sqrt{2}}{\sigma_h}   \bE\bigg[ |\bareta_i - \bE[\bareta_i]|  \cdot  \bigg(\sqrt{ \frac{\sum_{\boldj \in \Inm, i \in \boldj}h^2(\bX_\boldj )}{{n \choose m}}   } + \frac{\sqrt{2}\sum_{\boldj \in \Inm, i \in \boldj}h^2(\bX_\boldj )}{ {n \choose m}\sigma_h}\bigg)  \bigg] \text{ by \propertyref{censoring_property}$(i)$}  \notag  \\
&\leq \frac{\sqrt{2}}{\sigma_h}   \bE\bigg[ |\bareta_i - \bE[\bareta_i]| \cdot  \bigg(\sqrt{ \frac{\sum_{\boldj \in \Inm, i \in \boldj}\bE[h^2(\bX_\boldj )\mid X_i]}{{n \choose m}}   } + \frac{\sqrt{2}\sum_{\boldj \in \Inm, i \in \boldj}\bE[h^2(\bX_\boldj ) \mid X_i]}{ {n \choose m}\sigma_h}\bigg)  \bigg] \notag\\
&\leq  \frac{\sqrt{2}}{\sigma_h}   \bE\bigg[ (\bareta_i + \bE[\bareta_i]) \cdot  \bigg(\sqrt{ \frac{m\Psi_1(X_i)}{n}   } + \frac{\sqrt{2}m\Psi_1(X_i )}{ n \sigma_h}\bigg)  \bigg] \notag\\
&\leq   \frac{\sqrt{2}}{\sigma_h} \bigg\{   \bE\bigg[ \bareta_i    \bigg(\sqrt{ \frac{m\Psi_1(X_i)}{n}   } + \frac{\sqrt{2}m\Psi_1(X_i )}{ n \sigma_h}\bigg)  \bigg]+   \frac{(1 + \sqrt{2})m^{3/2}\bE[\Psi_1] \bE[\Psi_1^{1/2}]}{n^{3/2} \sigma_h^2} 
\bigg\} \text{ by } \bE[\bareta_i] \leq \bE[\bareta_i^{1/2}]\notag\\
&\leq  \frac{\sqrt{2}}{\sigma_h} \bigg\{   \frac{(1 + \sqrt{2})m^{3/2} \bE[\Psi_1^{3/2}]}{n^{3/2} \sigma_h^2}+   \frac{(1 + \sqrt{2})m^{3/2}\bE[\Psi_1] \bE[\Psi_1^{1/2}]}{n^{3/2} \sigma_h^2} 
\bigg\}  \text{ by } \bareta_i \leq \bareta_i^{1/2}\notag\\
&\leq \frac{(4 + 2\sqrt{2})m^{3/2} \bE[\Psi_1^{3/2}]}{n^{3/2} \sigma_h^3}
\label{(4)_bdd} \text{ by } \bE[\Psi_1] \bE[\Psi_1^{1/2}] \leq \bE[\Psi_1^{3/2}] .
\end{align}
Combining \eqref{first_and_second_pieces_pieces}-\eqref{(4)_bdd}, we have proved the bound in \eqref{cool_argument_prep_1}.
\subsection{Bound on the term in \eqref{cool_argument_prep_2}}
We will treat both configurations in \eqref{leave_one_out_frakz} together. 
By the independence between $X_i$ and $\{X_j\}_{j \neq i}$, we have 
\begin{equation*}
\bE\bigg[ \frakz_{\cX^{(i)}} \Big(\bareta_i - \bE[\bareta_i] \Big)  f_{\frakz_{\cX^{(i)}}}\Big(\sum_{\boldj \in \Inm, i \not\in\boldj} {_i}\xi_\boldj\Big) \bigg] 
= \underbrace{\bE\big[\bareta_i - \bE[\bareta_i]\big]}_{=0}  \bE\bigg[\frakz_{\cX^{(i)}} f_{\frakz_{\cX^{(i)}}}\Big(\sum_{\boldj \in \Inm, i \not\in\boldj} {_i}\xi_\boldj\Big) \bigg] ,
\end{equation*}
 so  we can  write 
 \begin{align}
&\bE\Big[ \frakz_{\cX^{(i)} }  ( \bareta_i - \bE[\bareta_i])  f_{\frakz_{\cX^{(i)} } } \Big(\sum_{\boldj \in I_{n, m}} {_i}\xi_\boldj \Big) \Big] \notag\\
&= \bE\bigg[ \frakz_{\cX^{(i)} } ( \bareta_i - \bE[\bareta_i])    \bigg( f_{\frakz_{\cX^{(i)} }}\Big(\sum_{\boldj \in \Inm}  {_i}\xi_\boldj\Big) -  f_{\frakz_{\cX^{(i)} }}\Big(\sum_{\boldj \in \Inm, i \not\in\boldj} {_i}\xi_\boldj\Big) \bigg) \bigg]  \notag\\
& \hspace{4cm} \text{ by the independence between $X_i$ and $\{X_j\}_{j \neq i}$} \notag\\
&=   \bE\bigg[\frakz_{\cX^{(i)} } ( \bareta_i - \bE[\bareta_i])   \int_0^{\sum_{\boldj \in \Inm, i \in \boldj} {_i}\xi_\boldj} \bE\bigg[\frac{\partial f_{\frakz_{\cX^{(i)} } }}{\partial w}\Big( \sum_{\boldj \in \Inm, i \not\in\boldj} {_i}\xi_\boldj+ t\Big)\mid \cX\bigg]d t 
\bigg]  \notag \\
& \hspace{1cm} \text{  by the independence of the $Z_\boldj$'s and the fundamental theorem of calculus}. \notag
\end{align}
From here, with both $\bareta_i$ and $\bE[\bareta_i]$ being non-negative, we can immediately write
\begin{multline}\label{first_blue_bdd}
\Big|\bE\Big[ \frakz_{\cX^{(i)} }  ( \bareta_i - \bE[\bareta_i])  f_{\frakz_{\cX^{(i)} } } \Big(\sum_{\boldj \in I_{n, m}} {_i}\xi_\boldj \Big) \Big] \Big|
 \\
\leq \bE\Bigg[ ( \bareta_i + \bE[\bareta_i])   \int_0^{|\sum_{\boldj \in \Inm, i \in \boldj} {_i}\xi_\boldj|} \bE \bigg[ \Big|\frakz_{\cX^{(i)}} \frac{\partial f_{\frakz_{\cX^{(i)}}}}{\partial w}\Big( \sum_{\boldj \in \Inm, i \not\in\boldj} {_i}\xi_\boldj+ t\Big)\Big| \mid \cX \bigg] dt.
\Bigg]
\end{multline}
At this point, we will bound the integrand $\bE [ |\frakz_{\cX^{(i)}} \frac{\partial f_{\frakz_{\cX^{(i)}}} }{\partial w}( \sum_{\boldj \in \Inm, i \not\in\boldj} {_i}\xi_\boldj+ t)| \mid \cX]$ in \eqref{first_blue_bdd}. By $\sup_{w, z \in \bR}|\frac{\partial f_z}{\partial w}(w)| \leq 1$ from \lemref{fxfxprime_bdd}$(i)$, 
\begin{equation}\label{integrand_bdd_one}
\bE \bigg[ \Big| \frakz_{\cX^{(i)}} \frac{\partial f_{\frakz_{\cX^{(i)}}} }{\partial w}\Big( \sum_{\boldj \in \Inm, i \not\in\boldj} {_i}\xi_\boldj+ t\Big) \Big|\cdot I\big(|\frakz_{\cX^{(i)}}| \leq 2\big) \mid \cX \bigg]  \leq 2 
.
\end{equation}
Since $\sup_{z, w \in \bR }| \frac{\partial f_z}{\partial w}(w)| \leq 1$  by  \lemref{fxfxprime_bdd}$(i)$ and  $|\frac{\partial f_z}{\partial w}(w)| \leq e^{1/2 - |z|}$ for any $|z| >2$ and $|w| \leq |z| -1$ by  \lemref{fxfxprime_bdd}$(ii)$ and $(iii)$,  for any $t \geq 0$, we have
\begin{align}
&\bE \bigg[ \Big| \frakz_{\cX^{(i)}} \frac{\partial f_{ \frakz_{\cX^{(i)}}} }{\partial w} \big( \sum_{\boldj \in \Inm, i \not\in\boldj} {_i}\xi_\boldj+ t\big)\Big|  \cdot I\big(|\frakz_{\cX^{(i)}}| > 2\big) \mid \cX \bigg]  \notag\\
&\leq  \bE\Big[|\frakz_{\cX^{(i)}}| e^{1/2-|\frakz_{\cX^{(i)}}|}\cdot I \Big(\big|\sum_{\boldj \in \Inm, i \not\in\boldj} {_i}\xi_\boldj+ t \big| \leq |\frakz_{\cX^{(i)}}| -1,|\frakz_{\cX^{(i)}}| >2 \Big)\mid \cX\Big]  \notag\\
& \hspace{3cm}+ \bE\Big[|\frakz_{\cX^{(i)}}|  \cdot I \Big( \big|\sum_{\boldj \in \Inm, i \not\in\boldj} {_i}\xi_\boldj+ t \big| > |\frakz_{\cX^{(i)}}| -1, |\frakz_{\cX^{(i)}}| >2\Big) \mid \cX\Big] \notag\\
&\leq  \bE\Big[|\frakz_{\cX^{(i)}}| e^{1/2-|\frakz_{\cX^{(i)}}|}\cdot I \Big(\big|\sum_{\boldj \in \Inm, i \not\in\boldj} {_i}\xi_\boldj+ t \big| \leq |\frakz_{\cX^{(i)}}| -1,|\frakz_{\cX^{(i)}}| >2 \Big)\mid \cX\Big]  \notag\\
& \hspace{3cm} + \bE\bigg[|\frakz_{\cX^{(i)}}|  \cdot  \frac{ \big|\sum_{\boldj \in \Inm, i \not\in\boldj} {_i}\xi_\boldj|+ t }{|\frakz_{\cX^{(i)}}| -1} \cdot I\big( |\frakz_{\cX^{(i)}}| >2\big) \mid \cX\bigg] \notag\\
&\leq 2 + 2\Big( \bE\Big[ \Big|\sum_{\boldj \in \Inm, i \not\in\boldj} {_i}\xi_\boldj \Big| \mid \cX\Big]+t\Big) \text{ by \eqref{some_sup_bdds}} \notag \\
&\leq 2  + 2\bigg( \Big(\bE\Big[ \sum_{\boldj \in \Inm, i \not\in\boldj} {_i}\xi_\boldj^2\mid \cX\Big]\Big)^{1/2} + t\bigg) 
\text{ by } \bE\Big[ ( \sum_{\boldj \in \Inm, i \not\in\boldj} {_i}\xi_\boldj )^2\mid \cX\Big]=  \bE\Big[ \sum_{\boldj \in \Inm, i \not\in\boldj} {_i}\xi_\boldj^2\mid \cX\Big]\notag\\
&\leq  2+ 2(1+t)   \quad \text{ by } \quad  \bE\Big[ \sum_{\boldj \in \Inm, i \not\in\boldj} {_i}\xi_\boldj^2\mid \cX\Big] =   U_{h^2}^{(i)} / \bar{U}_{h^2}^{(i)} \leq 1 .    \label{integrand_bdd_two}
\end{align}
Hence, from \eqref{integrand_bdd_one} and \eqref{integrand_bdd_two}, we have 
\[
\bE \bigg[ \Big| \frakz_{\cX^{(i)}} \frac{ \partial f_{ \frakz_{\cX^{(i)}}}  }{\partial w}\big( \sum_{\boldj \in \Inm, i \not\in\boldj} {_i}\xi_\boldj+ t\big)\Big|  \mid \cX \bigg] \leq 6 + 2t;
\]
plugging this latter fact into \eqref{first_blue_bdd}, we get 
\begin{align}
&\Big| \bE\Big[ \frakz_{\cX^{(i)} }  ( \bareta_i - \bE[\bareta_i])  f_{\frakz_{\cX^{(i)} } } \Big(\sum_{\boldj \in I_{n, m}} {_i}\xi_\boldj \Big) \Big] \Big| \notag\\
&\leq \bE\bigg[ ( \bareta_i + \bE[\bareta_i])   \int_0^{|\sum_{\boldj \in \Inm, i \in \boldj} {_i}\xi_\boldj|} 6 + 2t \;\ dt
\bigg] \notag\\
&= \bE\bigg[ ( \bareta_i + \bE[\bareta_i])  \bigg( 6 \Big|\sum_{\boldj \in \Inm,i \in \boldj} {_i}\xi_\boldj \Big| + \Big(\sum_{\boldj \in \Inm, i \in \boldj} {_i}\xi_\boldj \Big)^2 \bigg)
\bigg] \notag\\
&\leq \bE\bigg[ ( \bareta_i + \bE[\bareta_i])  \bigg( 6 \Big(\bE \Big[ \Big(\sum_{\boldj \in \Inm,i \in \boldj} {_i}\xi_\boldj \Big)^2  \mid X_i \Big] \Big)^{1/2}+ \bE\Big[\Big(\sum_{\boldj \in \Inm, i \in \boldj} {_i}\xi_\boldj\Big)^2  \mid X_i \Big] \bigg)\bigg] \notag\\
&\leq \bE\bigg[ ( \bareta_i + \bE[\bareta_i])  \bigg( 6 \bigg( \frac{2m \Psi_1(X_i)}{n \sigma_h^2}
 \bigg)^{1/2}+ \frac{2m \Psi_1(X_i)}{n \sigma_h^2}\bigg)
\bigg], \label{second_blue_bdd}
\end{align}
where \eqref{second_blue_bdd} comes from evaluating and bounding the conditional expectation as
\begin{multline*}
\bE\Big[ \Big(\sum_{ \boldj \in \Inm, i \in \boldj} {_i}\xi_\boldj \Big)^2 \mid X_i \Big] 
= 
\bE\Big[\bE\Big[ \Big(\sum_{ \boldj \in \Inm, i \in \boldj} {_i}\xi_\boldj \Big)^2 \mid \cX \Big]\mid X_i \Big] \\
= \sum_{\boldj \in \Inm, i \in \boldj} \frac{\bE[h^2(\bX_\boldj) \mid X_i]}{{n \choose m} \bar{U}^{(i)}_{h^2}} = \frac{m \Psi_1(X_i)}{n \bar{U}^{(i)}_{h^2}} \leq \frac{2m \Psi_1(X_i)}{n \sigma_h^2}.
\end{multline*}
Since  $\bE[\bareta_i] \leq \bE[\bareta_i^{1/2}] $, we can then further continue from \eqref{second_blue_bdd} as
\begin{align}
&\Big|\bE\Big[ \frakz_{\cX^{(i)} }  ( \bareta_i - \bE[\bareta_i])  f_{\frakz_{\cX^{(i)} } } \Big(\sum_{\boldj \in I_{n, m}} {_i}\xi_\boldj \Big) \Big] \Big| \notag\\
&\leq  \bE\bigg[ \bareta_i   \bigg( 6\bigg( \frac{2m \Psi_1(X_i)}{n \sigma_h^2}
 \bigg)^{1/2}+ \frac{2m \Psi_1(X_i)}{n \sigma_h^2}\bigg)
\bigg] +\frac{(6\sqrt{2} +2)m^{3/2} \bE[\Psi_1] \bE[\Psi_1^{1/2}]}{n^{3/2} \sigma_h^3}\notag\\
 &\leq 
  \frac{(6\sqrt{2} + 2)m^{3/2} \bE[\Psi_1^{3/2}]}{n^{3/2} \sigma_h^3}
+   \frac{(6\sqrt{2} +2)m^{3/2} \bE[\Psi_1] \bE[\Psi_1^{1/2}]}{n^{3/2} \sigma_h^3}  \text{ by }  \bareta_i \leq \bareta_i^{1/2} \notag \\
&\leq   \frac{  (12\sqrt{2} +4)  m^{3/2} \bE[\Psi_1^{3/2}]}{n^{3/2} \sigma_h^3} \label{3_final_bdd}, 
\end{align}
where the last inequality  is true via $\bE[\Psi_1] \bE[\Psi_1^{1/2}] \leq \bE[\Psi_1^{3/2}]$. So \eqref{cool_argument_prep_2} is proved.

\section{Supplementary calculations for  \appref{bounding_T_D2} } \label{app:tediousproof}
This section proves  \eqref{tail_prob_bdd_for_D2}-\eqref{x_exp_barD2_f_indicator_bdd_final} in \appref{bounding_T_D2}.


\subsection{Proof of \eqref{tail_prob_bdd_for_D2}} \label{app:Proof_of_tail_prob_bdd_for_D2}
Since $\bar{U}_{h^2} > \sigma_h^2/2$, we have
\begin{align*}
P(|D_3| > 1/2) &\leq P( |\Pi_1 + \Pi_2|> 1/4)\\
&\leq P( |\Pi_1| > 1/8) + P( |\Pi_2|> 1/8) \\
&\leq 64 \bE[\Pi_1^2] + 8 \|\Pi_2 \|_{3/2}
\end{align*}

\subsection{Proof of  \eqref{something1}} \label{app:proof_of_something1}

Recognizing that  $|\bar{D}_3|$ amounts to the non-negative random variable $|D_3|$   upper-censored at $1/2$, one can use \propertyref{censoring_property}$(ii)$ and $\bar{U}_{h^2} > \sigma_h^2/2$ to get
\begin{equation}\label{first_bdd_barD2}
|\bar{D}_3| \leq  |D_3|\leq  2(|\Pi_1| + |\Pi_2|)
\end{equation}
Since  the value of $|\bar{D}_3|$ must be capped at $1/2$, 
\eqref{first_bdd_barD2} must  also necessarily implies
\begin{equation} \label{D2bar_bdd_Pi1_Pi2}
 |\bar{D}_3| \leq 2 (|\bar{\Pi}_1| + |\bar{\Pi}_2|)  
\end{equation}
where $\bar{\Pi}_1$ and $\bar{\Pi}_2$ represent  $\Pi_1$ and $\Pi_2$ censored within the interval $[-1, 1]$ as
\[
\bar{\Pi}_k \equiv  \Pi_k I(|\Pi_k| \leq 1) + I(\Pi_k \geq 1) - I(\Pi_k < -1) \text{ for } k = 1, 2.
\]
Hence, \eqref{D2bar_bdd_Pi1_Pi2} gives \eqref{something1} by \propertyref{censoring_property}$(ii)$ and $|\bar{\Pi}_2|^2 \leq |\bar{\Pi}_2|$.

\subsection{Proof of \eqref{x_exp_barD2_f_indicator_bdd_final} } \label{app:proof_of_x_exp_barD2_f_indicator_bdd1}
By the definition of $D_3$ in  \eqref{D3_def}, one can first write 
\begin{multline} \label{x_exp_barD2_f_indicator_bdd1a}
 \bE\Big[ \frakz_\cX \bar{D}_3 f_{\frakz_\cX} \Big(\sum_{\boldi \in I_{n, m}} \barxi_\boldi \Big)  \Big]
 =  
 - \bE\Big[ \frakz_\cX\frac{\sigma_h^2 \Pi_1}{\bar{U}_{h^2}} f_{\frakz_\cX} \Big(\sum_{\boldi \in I_{n, m}} \barxi_\boldi \Big)  \Big]
 - \bE\Big[ \frakz_\cX\frac{\sigma_h^2 \Pi_2}{\bar{U}_{h^2}} f_{\frakz_\cX} \Big(\sum_{\boldi \in I_{n, m}} \barxi_\boldi \Big)  \Big]\\
 -\bE\Big[ \frakz_\cX (D_3 - 1/2) f_{\frakz_\cX} \Big(\sum_{\boldi \in I_{n, m}} \barxi_\boldi \Big) I(D_3 > 1/2) \Big] 
  -\bE\Big[ \frakz_\cX (D_3 + 1/2) f_{\frakz_\cX}\Big(\sum_{\boldi \in I_{n, m}} \barxi_\boldi \Big) I(D_3 < - 1/2) \Big].
\end{multline}
Given  $\sup_{z, w \in \bR} |z f_z(w) |\leq 1$ from \lemref{fxfxprime_bdd}$(iv)$, from   \eqref{D3_def} it is easy to see that 
\begin{align}
\bigg| \bE\Big[ \frac{\frakz_\cX}{2} f_{\frakz_\cX} \Big(\sum_{\boldi \in I_{n, m}} \barxi_\boldi \Big)  I(D_3 > 1/2)\Big] \bigg| 
 &\leq \frac{1}{2} \bE\bigg[I\bigg( |\Pi_1 + \Pi_2|>\frac{1}{4}\bigg) \bigg] \notag\\
&\leq \frac{1}{2} \bE\bigg[ I\bigg(|\Pi_1|>\frac{1}{8}\bigg) + I\bigg(|\Pi_2|>\frac{1}{8}\bigg)     \bigg] \notag\\
&\leq 32  \bE[\Pi_1^2] + 4 \|\Pi_2 \|_{3/2}.  \label{x_exp_barD2_f_indicator_bdd1b}
\end{align}
Similarly, we also have
\begin{equation}
\bigg| \bE\Big[ \frac{\frakz_\cX}{2} f_{\frakz_\cX} \Big(\sum_{\boldi \in I_{n, m}} \barxi_\boldi \Big) I(D_3 <- 1/2)\Big] \bigg| \leq 
32  \bE[\Pi_1^2] + 4 \|\Pi_2 \|_{3/2}.  \label{x_exp_barD2_f_indicator_bdd1b}
\end{equation}
Next, using  \lemref{fxfxprime_bdd}$(iv)$ and    \eqref{D3_def} again, 
\begin{align}
&\bigg| \bE\Big[ \frakz_\cX D_3 f_{\frakz_\cX}\Big(\sum_{\boldi \in I_{n, m}} \barxi_\boldi \Big)I(|D_3 |>1/2)\Big] \bigg| \notag\\
&\leq  \bE\bigg[|D_3| I\bigg(|D_3|>\frac{1}{2}\bigg)\bigg] \notag\\
&\leq  2 \bE\bigg[ |\Pi_1 + \Pi_2| I\bigg( |\Pi_1 + \Pi_2|>\frac{1}{4}\bigg)\bigg] \notag\\
&\leq   8\bE\big[ \big|\Pi_1 | \cdot |\Pi_1 + \Pi_2|\big] + 2 \|\Pi_2\|_{3/2}\notag\\
&\leq  8( \bE[\Pi_1^2] + \| \Pi_1\|_3 \|\Pi_2\|_{3/2} )+ 2 \|\Pi_2\|_{3/2}.\label{x_exp_barD2_f_indicator_bdd1c}
\end{align}
Lastly, using \lemref{fxfxprime_bdd}$(iv)$ again, 
\begin{equation} \label{running_out_of_label}
\Big| \bE\Big[ \frakz_\cX \frac{\sigma_h^2 \Pi_2}{\bar{U}_{h^2}} f_{\frakz_\cX} \Big(\sum_{\boldi \in I_{n, m}} \barxi_\boldi \Big)  \Big] \Big|
\leq 2 \bE[|\Pi_2|] \leq 2 \|\Pi_2\|_{3/2}.
\end{equation}
Combining \eqref{x_exp_barD2_f_indicator_bdd1a}-\eqref{running_out_of_label} gives
\begin{multline} \label{x_exp_barD2_f_indicator_bdd1_not_yet}
\Big| \bE\Big[ \frakz_\cX \bar{D}_3 f_{\frakz_\cX} \Big(\sum_{\boldi \in I_{n, m}} \barxi_\boldi \Big)  \Big]\Big| \leq 
 \bigg|\bE\bigg[ \frakz_\cX    \frac{\sigma_h^2 \Pi_1}{\bar{U}_{h^2}} f_{\frakz_\cX}\bigg(\sum_{\boldj \in \Inm} \barxi_\boldj\bigg) \bigg] \bigg|  
\\+ C\Big(\bE[\Pi_1^2] + \| \Pi_1\|_3 \|\Pi_2\|_{3/2} + \|\Pi_2\|_{3/2}  \Big)
\end{multline}

We will further bound the term $\big|\bE\big[ \frakz_\cX    \frac{\sigma_h^2 \Pi_1}{\bar{U}_{h^2}} f_{\frakz_\cX  }\big(\sum_{\boldj \in \Inm} \barxi_\boldj\big) \big] \big|$ in \eqref{x_exp_barD2_f_indicator_bdd1_not_yet}; first write
\begin{align}
&\bE\bigg[   \frakz_\cX   \frac{\sigma_h^2  \Pi_1 }{\bar{U}_{h^2}}  f_{\frakz_\cX } \Big(\sum_{\boldj \in I_{n, m}} \barxi_\boldj \Big)  \bigg] \notag\\
&=  \bE\bigg[   \frakz_\cX  \cdot  \Pi_1 \cdot \frac{\sigma_h^2}{\bar{U}_{h^2}} \cdot \Big( 1 -\frac{\bar{U}_{h^2}}{\sigma_h^2} \Big)\cdot f_{\frakz_\cX} \Big(\sum_{\boldj \in I_{n, m}} \barxi_\boldj \Big)  \bigg] 
+  \bE\Big[  \frakz_\cX    \Pi_1  f_{\frakz_\cX} \Big(\sum_{\boldj \in I_{n, m}} \barxi_\boldj \Big)  \Big]. \label{paranthesis_12_break}
\end{align}
Using  \lemref{fxfxprime_bdd}$(iv)$,  we have
\begin{align}
&\bigg| \bE\bigg[  \frakz_\cX \cdot \Pi_1 \cdot \frac{\sigma_h^2}{\bar{U}_{h^2}} \cdot \Big( 1 -\frac{\bar{U}_{h^2}}{\sigma_h^2} \Big)\cdot f_{ \frakz_\cX} \Big(\sum_{\boldj \in I_{n, m}} \barxi_\boldj \Big)  \bigg]\bigg| \notag\\
&\leq \bE\Big[ \Big|\Pi_1 \cdot \frac{\sigma_h^2}{\bar{U}_{h^2}}  \cdot \Big( 1 -\frac{\bar{U}_{h^2}}{\sigma_h^2}\Big)  \Big|\Big] \notag\\
& \leq 2\bE\Big[ \Big|\Pi_1   \Big( 1 -\frac{\bar{U}_{h^2}}{\sigma_h^2}\Big) \Big|\Big]  \notag\\
&=2\bE\Big[ \Big|\Pi_1   \Big( 1 -\frac{U_{h^2}}{\sigma_h^2}\Big) \Big|  \cdot I(\cE_{1, \cX} )\Big]  
+ 2\bE\Big[ |\Pi_1   |  \cdot I\big(\Omega \backslash \cE_{1, \cX} \big)\Big]
\notag\\
&\leq 2 \|\Pi_1\|_3 \bigg( \Big\|\sum_{i=1}^n\frac{m (\Psi_1(X_i) - \sigma_h^2)}{n \sigma_h^2} \Big\|_{3/2}+ \|\Pi_2\|_{3/2} \bigg) \notag \\
& \hspace{2cm}+ 2 \|\Pi_1\|_{3/2} \sum_{i=1}^n \bE[ (\eta_i - 1) I(\eta_i >1)]  +2 \|\Pi_1\|_3 P(\Omega \backslash \cE_{1, \cX})^{2/3}   \text{ by } \eqref{H_decomp_h2_ustat} \notag\\
&\leq 2 \|\Pi_1\|_3 \bigg(  \frac{ 2^{5/3}m \|\Psi_1\|_{3/2}}{\sigma_h^2 n^{1/3}}  + \frac{m^{3/2}\bE[\Psi_1^{3/2}]}{ n^{1/2}\sigma_h^3} \bigg) + \frac{2 \|\Pi_1\|_{3/2} m^{3/2}\bE[\Psi_1^{3/2}]}{ n^{1/2}\sigma_h^3}  + 2  \|\Pi_1\|_3 P(\Omega \backslash \cE_{1, \cX})^{2/3}, \label{paranthesis_1}
 \end{align}
where the last inequality uses  \citet[Lemma 1]{chatterji1969p}'s inequality
\begin{multline*}
\bE\Big[ \Big(\sum_{i=1}^n\frac{m (\Psi_1(X_i) - \sigma_h^2)}{n \sigma_h^2}  \Big)^{3/2} \Big] 
\leq  \frac{2m^{3/2} \bE[| \Psi_1(X_1) - \sigma_h^2 |^{3/2}]}{\sigma_h^3\sqrt{n}} \\
\leq  \frac{2m^{3/2} \bE[| \Psi_1(X_1) + \sigma_h^2 |^{3/2}]}{\sigma_h^3\sqrt{n}} \leq \frac{(2m)^{3/2} (\bE[ \Psi_1^{3/2} ] + \sigma_h^3)}{\sigma_h^3\sqrt{n}}  \leq  \frac{2(2m)^{3/2} \bE[ \Psi_1^{3/2} ] }{\sigma_h^3\sqrt{n}}
 \end{multline*}
 and \eqref{remnant_expectation_bdd}. 
 Since 
 \begin{multline*}
 \|\Pi_1\|_3  \frac{m \|\Psi_1\|_{3/2}}{\sigma_h^2 n^{1/3}} \leq \bE[|\Pi_1|^3] +  \frac{m^{3/2} \bE[\Psi_1^{3/2}]}{\sigma_h^3 n^{1/2}}, \\
  \|\Pi_1\|_3 \|\Pi_2\|_{3/2} \leq \bE[|\Pi_1|^3] + \bE[|\Pi_2|^{3/2}] \text{ and }\\
    \|\Pi_1\|_3 P(\Omega \backslash \cE_{1, \cX})^{2/3} \leq \bE[|\Pi_1|^3] + P(\Omega \backslash \cE_{1, \cX}),
 \end{multline*} 
   \eqref{x_exp_barD2_f_indicator_bdd1_not_yet}-\eqref{paranthesis_1} together imply
\begin{multline*} 
\Big| \bE\Big[ \frakz_\cX\bar{D}_3 f_{\frakz_\cX} \Big(\sum_{\boldi \in I_{n, m}} \barxi_\boldi \Big)  \Big]\Big| \leq 
 \Big| \bE\Big[ \frakz_\cX  \Pi_1  f_{\frakz_\cX} \Big(\sum_{\boldj \in I_{n, m}} \barxi_\boldj \Big)  \Big] \Big| +
\\
 C\bigg\{ \frac{m^{3/2} \bE[\Psi_1^{3/2}]}{\sigma_h^3 n^{1/2}} \cdot \big(1 +   \|\Pi_1\|_{3/2}  \big)+  \bE[\Pi_1^2] + \bE[|\Pi_1|^3] + \bE[|\Pi_2|^{3/2}]+ \|\Pi_2\|_{3/2}  \bigg\}+ P(\Omega \backslash \cE_{1, \cX}) ,
\end{multline*}
which can be slightly simplified as 
\begin{multline} \label{paranthesis_1_prime}
\Big| \bE\Big[  \frakz_\cX \bar{D}_3 f_{\frakz_\cX} \Big(\sum_{\boldi \in I_{n, m}} \barxi_\boldi \Big)  \Big]\Big| \leq 
 \Big| \bE\Big[ \frakz_\cX \Pi_1  f_{\frakz_\cX} \Big(\sum_{\boldj \in I_{n, m}} \barxi_\boldj \Big)  \Big] \Big| +
\\
 C\bigg\{ \frac{m^{3/2} \bE[\Psi_1^{3/2}]}{\sigma_h^3 n^{1/2}} \cdot \big(1 +   \|\Pi_1\|_{3/2}  \big) + \bE[\Pi_1^2] + \bE[|\Pi_1|^3] +  \|\Pi_2\|_{3/2}  \bigg\} + P(\Omega \backslash \cE_{1, \cX})\end{multline}
by absorbing $\bE[|\Pi_2|^{3/2}]$ into $\|\Pi_2\|_{3/2}$; this is because
$\Big| \bE\Big[\frakz_\cX \bar{D}_3 f_{\frakz_\cX} \Big(\sum_{\boldi \in I_{n, m}} \barxi_\boldi \Big)  \Big]\Big| \leq 1/2$ by \lemref{fxfxprime_bdd}$(iv)$, one can without loss of generality assume that $\|\Pi_2\|_{3/2} \leq 1$ (otherwise the bound is vacuous), so $\bE[|\Pi_2|^{3/2}] \leq \|\Pi_2\|_{3/2}$.

Lastly, we will bound the term $ \Big| \bE\Big[ \frakz_\cX  \Pi_1  f_{\frakz_\cX} \Big(\sum_{\boldj \in I_{n, m}} \barxi_\boldj \Big)  \Big] \Big| $ on the right hand side of \eqref{paranthesis_1_prime}. Recall the definition of $\xi_\boldj$ in \eqref{xi_def}, and define the event
\begin{equation*} \label{cE_3_cX_def}
\scE_{\cX, \cZ}  = \Big\{\max_{\boldj \in I_{n, m}}  |\xi_\boldj|  \leq 1\Big\}
\end{equation*}
which depends on both $\cX$ and $\cZ$; it will be useful to note that 
\begin{equation} \label{cE_3_bdd}
P\Big(\Omega \backslash \scE_{\cX, \cZ} \mid  \cX\Big) \leq \sum_{\boldj \in \Inm} \bE[ |\xi_\boldj|^3 \mid \cX] = \frac{ U_{|h|^3} (1 - 2p + 2p^2)}{\bar{U}_{h^2}^{3/2} \sqrt{N(1 - p)}} \leq  
\frac{2^{3/2}U_{|h|^3} (1 - 2p + 2p^2)}{\sigma_h^3 \sqrt{N(1 - p)}},
\end{equation}
by the computation of $\bE [ |\xi_\boldj|^3 \mid \cX ] $ in \eqref{cX_cond_moment_properties_of_xi}.
One can then bound   it  as
\begin{align}
&\bigg| \bE\Big[  \frakz_\cX \Pi_1 f_{\frakz_\cX} \Big(\sum_{\boldj \in I_{n, m}} \barxi_\boldj \Big)  \Big]\bigg|\notag\\
&\leq \bigg| \bE\bigg[  \frakz_\cX \Pi_1     f_{ \frakz_\cX}\Big(\sum_{\boldj \in I_{n, m}} \xi_\boldj \Big) I(\scE_{\cX, \cZ}) \bigg] \bigg| +
 \bigg|\bE\bigg[  \frakz_\cX \Pi_1     f_{ \frakz_\cX} \Big(\sum_{\boldj \in I_{n, m}} \barxi_\boldj \Big) I(\Omega \backslash \scE_{\cX, \cZ}  ) \bigg] \bigg|\notag\\
 &\leq  \bigg| \bE\bigg[   \frakz_\cX\Pi_1     f_{ \frakz_\cX} \Big(\sum_{\boldj \in I_{n, m}} \xi_\boldj \Big) I(\scE_{\cX, \cZ}) \bigg] \bigg| +
 \bE\big[  \big|\Pi_1      I(\Omega \backslash \scE_{\cX, \cZ}  ) \big| \big]  \text{ by  \lemref{fxfxprime_bdd}$(iv)$}\notag\\
 &\leq  \bigg| \bE\bigg[  \frakz_\cX \Pi_1     f_{ \frakz_\cX} \Big(\sum_{\boldj \in I_{n, m}} \xi_\boldj \Big)I(\scE_{\cX, \cZ}) \bigg] \bigg| +\frac{\sqrt{2}\|\Pi_1\|_{3/2}\|h\|_3 (1 - 2p + 2p^2)^{1/3}}{  \sigma_h (N(1-p))^{1/6}}, \label{2_breakdown}
\end{align}
where the last inequality comes from applying H\"older's inequality and \eqref{cE_3_bdd} on the right hand side of
\[
\bE\big[  \big|\Pi_1      I(\Omega \backslash \scE_{\cX, \cZ}  ) \big| \big]  = \bE[|\Pi_1| P(\Omega \backslash \scE_{\cX, \cZ} \mid  \cX)] \leq \bE\Big[|\Pi_1| \big(P(\Omega \backslash \scE_{\cX, \cZ} \mid  \cX)\big)^{1/3}\Big] .
\]
 Further,  $\bE\Big[  \frakz_\cX \Pi_1     f_{ \frakz_\cX} \Big(\sum_{\boldj \in I_{n, m}} \xi_\boldj \Big) I(\scE_{\cX, \cZ} ) \Big]$ in \eqref{2_breakdown} can be written as
\begin{multline}\label{2_breakdown_first_piece_bdd}
\bE\Big[ \frakz_\cX \Pi_1     f_{ \frakz_\cX} \Big(\sum_{\boldj \in I_{n, m}} \xi_\boldj \Big)  I(\scE_{\cX, \cZ} )   \Big] = \\
   \bE\Big[  \frakz_\cX\Pi_1 f_{ \frakz_\cX} \Big(\sum_{\boldj \in I_{n, m}} \xi_\boldj \Big)\Big] - 
 \bE\Big[  \frakz_\cX  \Pi_1f_{ \frakz_\cX} \Big(\sum_{\boldj \in I_{n, m}} \xi_\boldj \Big)I(\Omega \backslash  \scE_{\cX, \cZ})\Big];
\end{multline}
on the other hand,   \lemref{fxfxprime_bdd}$(iv)$ implies
\begin{align}
&\bigg|  \bE\Big[  \frakz_\cX \Pi_1f_{ \frakz_\cX} \Big(\sum_{\boldj \in I_{n, m}} \xi_\boldj \Big)I(\Omega \backslash \scE_{\cX, \cZ})\Big]\bigg| \notag\\
&\leq \|\Pi_1\|_{3/2} \| I(\Omega \backslash \scE_{\cX, \cZ})\|_3 \notag\\
&= \|\Pi_1\|_{3/2} \Big( \bE[ P(\Omega \backslash \scE_{\cX, \cZ} \mid \cX)]\Big)^{1/3} \notag\\
&\leq \|\Pi_1\|_{3/2} \Bigg( \frac{ 2^{3/2} \bE[|h|^3] (1 - 2p + 2p^2)}{\sigma_h^3 \sqrt{N(1 - p)}}\Bigg)^{1/3}  \text{ by \eqref{cE_3_bdd} } \notag\\
&\leq  \frac{\sqrt{2}  \|\Pi_1\|_{3/2} \|h\|_3 (1 - 2p + 2p^2)^{1/3}}{ \sigma_h (N(1-p))^{1/6}}.  
\label{2_breakdown_first_piece_bdd_bdd}
\end{align}
Putting \eqref{2_breakdown}-\eqref{2_breakdown_first_piece_bdd_bdd} together gives  
\begin{align} 
 &\Big| \bE\Big[  \frakz_\cX \Pi_1  f_{ \frakz_\cX}\Big(\sum_{\boldj \in I_{n, m}} \barxi_\boldj \Big)  \Big] \Big| \notag\\
 &\leq 
  \frac{2\sqrt{2} \|\Pi_1\|_{3/2} \|h\|_3 (1 - 2p + 2p^2)^{1/3}}{ \sigma_h (N(1-p))^{1/6}} 
+ \bigg|   \bE\bigg[ \frakz_\cX \Pi_1 f_{ \frakz_\cX} \Big(\sum_{\boldj \in I_{n, m}} \xi_\boldj \Big)\bigg]\bigg| \notag\\
&\leq 2\sqrt{2} \bigg(\bE[|\Pi_1|^{3/2}] + \frac{\bE[|h|^3] (1 - 2p + 2p^2)}{\sigma_h^3 (N(1-p))^{1/2}} \bigg)
+ \bigg|   \bE\bigg[  \frakz_\cX \Pi_1 f_{ \frakz_\cX} \Big(\sum_{\boldj \in I_{n, m}} \xi_\boldj \Big)\bigg]\bigg|. \label{2_breakdown_further}
\end{align} 
Combining \eqref{paranthesis_1_prime} and \eqref{2_breakdown_further} gives \eqref{x_exp_barD2_f_indicator_bdd_final}.

\bibliographystyle{plainnat}
\bibliography{icu_BE_ejp}


\begin{acks}
Dennis Leung is grateful to Mr. Liqian Zhang at the University of Manchester for his careful proofreading and fruitful discussion,  and to an anonymous referee whose comments have led to tremendous improvements of this paper.
\end{acks}

\end{document}